\documentclass[a4paper]{amsart}
\usepackage{amsmath,amssymb,amsfonts}
\usepackage{mathrsfs,latexsym,amsthm,stmaryrd,enumerate,epic,mathrsfs,dsfont}
\usepackage{amscd}
\usepackage[all]{xy}
\begin{document}
\title[Categorification of cell modules]{Categorification
of (induced) cell modules and the rough structure of generalized
Verma modules}

\author{Volodymyr Mazorchuk}
\address{V. M.: Department of Mathematics, Uppsala University (Sweden).}
\email{mazor\symbol{64}math.uu.se}
\thanks{The first author was supported by STINT, the Royal Swedish Academy
of Sciences, and the Swedish Research Council, the second author was
supported by EPSRC grant 32199}
\author{Catharina Stroppel}
\address{C. S.: Department of Mathematics, University of Glasgow
(United Kingdom).} \email{c.stroppel\symbol{64}maths.gla.ac.uk}

\numberwithin{equation}{section}

\newtheorem{proposition}{Proposition}
\newtheorem{lemma}[proposition]{Lemma}
\newtheorem{corollary}[proposition]{Corollary}
\newtheorem{theorem}[proposition]{Theorem}
\newtheorem{definition}[proposition]{Definition}
\newtheorem{conjecture}[proposition]{Conjecture}
\newtheorem{example}[proposition]{Example}
\newtheorem{remark}[proposition]{Remark}
\newcommand{\oplusop}[1]{{\mathop{\oplus}\limits_{#1}}}
\newcommand{\oplusoop}[2]{{\mathop{\oplus}\limits_{#1}^{#2}}}
\newcommand{\ff}{\footnote}
\renewcommand{\to}{\rightarrow}
\newtheorem{theoremintro}{Theorem}
\renewcommand{\thetheoremintro}{\Roman{theoremintro}}

\font\sc=rsfs10 at 12 pt
\font\scs=rsfs10 at 10 pt
\font\scb=rsfs10 at 16 pt
\font\scbb=rsfs10 at 18 pt

\newcommand{\ccC}{\mathscr{C}}
\newcommand{\ccS}{\mathscr{S}}

\def\la{\lambda}
\def\op{\operatorname}
\def\C{\mathbb C}
\def\R{\mathbb R}
\def\N{\mathbb N}
\def\Z{\mathbb Z}
\def\Q{\mathbb Q}
\def\g{\mathfrak g}
\def\p{\mathfrak p}
\def\h{\mathfrak h}
\def\n{\mathfrak n}
\def\mmm{\mathbf m}
\def\mmn{\mathbf n}
\newcommand{\mC}{\mathbb{C}}
\newcommand{\Ext}{\operatorname{Ext}}
\newcommand{\End}{\operatorname{End}}

\newcommand{\add}{\operatorname{add}}
\newcommand{\Ann}{\operatorname{Ann}}

\def\F{\mathbb F}
\def\S{\mathbb S}
\def\l{\lbrace}
\def\r{\rbrace}
\def\o{\otimes}
\def\lra{\longrightarrow}
\newcommand{\ba}{\mathbf{a}}
\newcommand{\cA}{\mathcal{A}}
\newcommand{\cB}{\mathcal{B}}
\newcommand{\cL}{\mathcal{L}}
\newcommand{\dmod}{\mathrm{-mod}}
\newcommand{\gmod}{\mathrm{-gmod}}
\newcommand{\mc}{\mathcal}
\newcommand{\mZ}{\mathbb{Z}}
\newcommand{\tto}{\twoheadrightarrow}
\newcommand{\mg}{\mathfrak{g}}
\newcommand{\mh}{\mathfrak{h}}
\newcommand{\ma}{\mathfrak{a}}
\newcommand{\mb}{\mathfrak{b}}
\newcommand{\cU}{\mathcal{U}}
\def\cF{\mathcal{F}}
\def\Hom{\textrm{Hom}}
\def\drawing#1{\begin{center} \epsfig{file=#1} \end{center}}
\def\mc{\mathcal}
\def\mf{\mathfrak}
\def\mb{\mathbb}

\def\yesnocases#1#2#3#4{\left\{ \begin{array}{ll} #1 & #2 \\ #3 & #4
\end{array} \right. }

\newcommand{\define}{\stackrel{\mbox{\scriptsize{def}}}{=}}
\def\hsm{\hspace{0.05in}}

\def\cO{\mathcal{O}}   
\def\cC{\mathscr{C}}
\def\sln{\mathfrak{sl}(n)}

\begin{abstract}
This paper presents categorifications of (right) cell modules and induced
cell modules for Hecke algebras of finite Weyl groups. In type $A$ we
show that these categorifications depend only on the isomorphism class of
the cell module, not on the cell itself. Our main application is multiplicity
formulas for parabolically induced modules over a reductive Lie algebra of
type $A$, which finally determines the so-called rough structure of generalized
Verma modules. On the way we present several categorification results and give
the positive answer to Kostant's problem from \cite{Jo} in many cases. We
also give a general setup of decategorification, precategorification and
categorification.
\end{abstract}
\maketitle

\tableofcontents

\section{Introduction}

The Weyl group acts via (exact) translation functors on the principal block
of the Bernstein-Gelfand-Gelfand category $\mathcal{O}$ associated with a
semi-simple complex finite-dimensional Lie algebra, see \cite{BG}. On the
level of the Grothendieck group, this becomes the regular representation of
the Weyl group. The nature of translation functors is such that they obviously
preserve several classes of modules - for example projective, injective or
tilting modules.

This naturally leads to the question whether the isomorphism classes of such
modules, considered as elements of the Grothendieck group, can be interpreted
in terms of the representation theory of the Weyl group, in particular in terms
of the regular representation of the Weyl group.

One of the most remarkable breakthrough results in the theory of semisimple
complex Lie algebras is that such an interpretation actually exists. The
connection is given by the so-called {\em Kazhdan-Lusztig theory}, which first
`upgrades' the Weyl group to the corresponding Hecke algebra, and also the
corresponding category $\mathcal{O}$ to its graded version, and then says that
the isomorphism classes of the graded indecomposable projective modules in the
regular block of the category $\mathcal{O}$ descend (on the level of the
Grothendieck group) to what is now known as the {\em Kazhdan-Lusztig basis} of
the Hecke algebra. The introduction of this Kazhdan-Lusztig basis together with
the Kazhdan-Lusztig conjecture (\cite[Conjecture~1.5]{KLCoxeter}) was a
milestone in combinatorial representation theory which finally turned the
computation of the character of any simple highest weight module for a complex
semisimple Lie algebra into a purely combinatorial task.

One main idea in this combinatorial representation theory showed up already
before \cite{KLCoxeter}, namely the idea of (left or right) {\it cells} for
finite Weyl group, in particular for the symmetric group. The latter was first
studied by combinatorialists (see e.g. \cite{Knuth}) and afterwards introduced
into representation theory (\cite{Jopreprint}, \cite{V}).

A natural consequence of the theory of cells is the definition of a special
class of modules for the Hecke algebra, namely the {\em cell modules}. In type
$A$ these modules contribute an exhaustive list of all irreducible modules. For
other types however, they are not irreducible in general.

The first objective of the present paper is to give a categorical version of
(right) cell modules. To each cell in the Weyl group $W$ we associate a certain
quotient category of some subcategory of the category $\mathcal{O}$ (of the
corresponding semisimple Lie algebra $\mg$) which is stable under the action of
translation functors. The categories used for this categorification are
indecomposable. When passing to the Grothendieck group we obtain the cell
module corresponding to our chosen cell. In other words: we categorify cell
modules for the Hecke algebra. Note that two different cells might have
isomorphic cell modules. In type $A$ the isomorphism classes of cell modules are
exactly the isomorphism classes of irreducible modules. We show that the
categorical picture is the same:

\begin{theoremintro}[Uniqueness theorem]\label{thmintro1}
Assume that $W$ is of type $A$. Then if two cell modules are isomorphic
then their categorifications are equivalent.
\end{theoremintro}

We will make this equivalence concrete by giving an explicit functor which
naturally commutes with the functorial action of the Hecke algebra. This is
what we call the `uniqueness' of categorifications (Theorem~\ref{thm6}).
As a result, we therefore have to each right cell $\mathbf{R}$ a
categorification $\mathscr{C}_{\mathbf{R}}$ together with an equivalence
$\Phi:\mathscr{C}_{\mathbf{R}_1}\to\mathscr{C}_{\mathbf{R}_2}$ whenever the
cell modules corresponding to $\mathbf{R}_1$  and $\mathbf{R}_2$  are isomorphic
(i.e. $\mathbf{R}_1$  and $\mathbf{R}_2$ are in the same double cell).
The Kazhdan-Lusztig cell theory equips the cell modules with a distinct basis
which corresponds in the categorification to the  isomorphism classes of
indecomposable projective modules.

Given a parabolic subgroup $W'$ of $W=S_n$, a right cell $\mathbf{R}'$ of
$W'$ and the corresponding cell module $\mathscr{C}_{\mathbf{R}'}$ of its
Hecke algebra $\mathds{H}(W')$ there is the induced cell module
$\mathscr{C}_{\mathbf{R}'}\otimes_{\mathds{H}(W')} \mathds{H}(W)$. To these
data we associate a certain category $\mathscr{X}=\mathscr{X}
(W,W',\mathbf{R}')$ of $\mg$-modules (in fact a subcategory of the category
$\cO$) such that the  following holds (for details see Theorem~\ref{thm53},
Proposition~\ref{cunique}, Theorem~\ref{combinatorics}):

\begin{theoremintro}\label{thm2}
\begin{enumerate}[(i)]
\item\label{thm2.1} The category $\mathscr{X}$ is a categorification of
$\mathscr{C}_{\mathbf{R}'}\otimes_{\mathds{H}(W')}\mathds{H}(W)$, with the
$\mathds{H}$-action given by translation functors.
\item\label{thm2.2} Up to equivalence $\mathscr{X}$ only depends on the
isomorphism class of $\mathscr{C}_{\mathbf{R}'}$, not on $\mathbf{R}'$ itself.
\item\label{thm2.3} There is a combinatorial description of $\mathscr{X}$ in
terms of Kazhdan-Lusztig polynomials in the following sense: the module
$\mathscr{C}_{\mathbf{R}'}\otimes_{\mathds{H}(W')}\mathds{H}(W)$ is
equipped with four natural bases corresponding to four natural classes of
modules in $\mathscr{X}$.
\end{enumerate}
\end{theoremintro}

A consequence of the (now proved) \cite[Conjecture 1.5]{KLCoxeter} is that
the Kazhdan-Lusztig basis of the Hecke algebra turns the problem of finding
multiplicities of composition factors of Verma modules into a purely
combinatorial statement: the multiplicities are given by evaluating the
corresponding Kazhdan-Lusztig polynomials. Verma modules are a special sort of
induced modules obtained by inducing one-dimensional (irreducible) modules over
a Borel subalgebra. In general, one would like to understand the structure of
modules obtained by inducing from an arbitrary irreducible module over a
parabolic subalgebra, ideally with a combinatorial description similar to the
case of Verma modules. This is however a very difficult task because of at
least two reasons: Firstly, there is no classification or reasonable
understanding of simple modules for finite dimensional complex Lie algebras
available (except for the Lie algebra $\mathfrak{sl}_2$, see \cite{Bl}), hence
the starting point for the induction process is not understood at all.
Secondly, it might happen that the induced modules are of infinite length (due
to a result of Stafford on existence of non-holonomic simple modules over the
Weyl algebra and $U(\mathfrak{sl}_2\times \mathfrak{sl}_2)$, see
\cite{Stafford}).

Nevertheless, our paper goes a big step further in solving these problems. The
principal idea is that we realize the induced module we are interested in, as a
(proper) standard object in some category which is equivalent to some
$\mathscr{X}$ as above. Then the Kazhdan-Lusztig theory together with
Theorem~\ref{thm2}\eqref{thm2.3} provides the necessary combinatorics
and as a result we can describe the so-called {\it rough structure} of
parabolically induced arbitrary simple modules.

One of the difficulties is actually to give a precise definition of what is
meant by {\it rough structure} (this is the topic of the last section of the
articles). In this introduction we just try to give the main idea. To do so let
$\mathfrak{g}$ be a Lie algebra with triangular decomposition. Let
$\mathfrak{p}$ be a parabolic subalgebra of $\mathfrak{g}$, and $V$ a simple
module over the reductive part of $\mathfrak{p}$. Then $V$ trivially extends
to a simple $\mathfrak{p}$-module, and the corresponding induced module
\begin{displaymath}
\Delta(\mathfrak{p},V)=U(\mathfrak{g})\otimes_{U(\mathfrak{p})}V
\end{displaymath}
is called a {\em generalized Verma module}. We want to describe the composition
factors of generalized Verma modules. Our main result is the following
statement (for details and notation see Section~\ref{s9},
in particular Theorem~\ref{s9.5-cor1}):

\begin{theoremintro}\label{thmintro3}
Assume that the reductive part of $\mathfrak{p}$ is of type $A$.
For $X,Y\in \op{Irr}^{\mathfrak{g}}\big(\mathcal{O}\{\mathfrak{p},\op{Coker}
(\overline{N}\otimes E)\}_{\op{int}}\big)$ we have the
following multiplicity formula in the category of $\mg$-modules:
\begin{equation}\label{blabla}
[\Delta(\mathfrak{p},V_X):L(\mathfrak{p},V_Y)]=
[\Delta(\mathfrak{p},V_{\hat{\xi}(X)}):
L(\mathfrak{p},V_{\hat{\xi}(Y)})].
\end{equation}
\end{theoremintro}

Here, the generalized Verma module $\Delta(\mathfrak{p},V_X)$ is the
one we are interested in, that means we want to describe the
multiplicities of the left hand side of the equation \eqref{blabla}.
On the other hand, $\Delta(\mathfrak{p},V_{\hat{\xi}(X)})$ is a generalized
Verma module induced from a simple {\em highest weight} module, and hence
is easier to understand. We will prove that the multiplicity on the right
hand side is given by Kazhdan-Lusztig combinatorics of a certain
category $\mathscr{X}$ as in Theorem~\ref{thm2} above. In fact, the
module $\Delta(\mathfrak{p},V_{\hat{\xi}(X)})$ belongs to one of the
four classes from Theorem~\ref{thm2}\eqref{thm2.3}. Therefore, it
becomes in principle possible to compute the multiplicities
completely. The only problem here is that  simple subquotients of
the form $L(\mathfrak{p},V_Y)$, $Y\in
\op{Irr}^{\mathfrak{g}}(\mathcal{O}\{\mathfrak{p},\op{Coker}
(\overline{N}\otimes E)\}_{\op{int}})$, do not exhaust all simple
subquotients of $\Delta(\mathfrak{p},V_X)$. Roughly speaking,
Theorem~\ref{thmintro3} gives information only about simple subquotients
having small enough annihilator. It turns out that the number and
multiplicities of such subquotients are always finite. {\em Knowing all the
multiplicities for these `allowed' simple subquotients is what we
call `knowing the rough structure of a generalized Verma module'}.
To prove Theorem~\ref{thmintro3} we use the approach of \cite{MiSo}
together with \cite{KM2}, which associates to simple modules of the
form $L(\mathfrak{p},V_Y)$ certain simple objects of some
$\op{Coker}$-category. Without any restriction on the simple module
$L$ to start with, our result seems to be the best possible, since
all we know in general about $L$ is its annihilator. The complete
(i.e. fine) structure of $\Delta(\mathfrak{p},V_X)$ depends  heavily
on $V_X$, not just on its annihilator. This becomes transparent by
comparing for instance the structure of generalized Verma modules
induced from Gelfand-Zetlin-modules on the one hand with generalized
Verma modules induced from simple Verma modules on the other hand.
In the first case the rough structure always coincides with the fine
structure (see for example \cite{MO-0}), whereas in the second  case
the fine structure is different from the rough structure already in
the case of the algebra $\mathfrak{sl}_3$ (this follows for example
from \cite[Theorem~7.6.23]{Di}).

Let now $\ma$ be the semisimple part of $\p$ and $W'$ the corresponding Weyl
group. As a consequence of Theorem~\ref{thmintro3}, we are able to deduce a
criterion for the irreducibility of the generalized Verma module
$\Delta(\mathfrak{p},L)$, where $L$ is an {\it arbitrary} simple $\ma$-module
(we formulate the statement in the case when  $L$ has trivial central
character, however, standard arguments extend this to the arbitrary type $A$
case, see Remark~\ref{simplicity}): We first associate in a combinatorial way
to $L$ a pair $(x,w)\in W'\times W$ (see Section~\ref{s9}) and then deduce the
following result:

\begin{theoremintro}\label{thmintro4}
The module $\Delta(\mathfrak{p},L)$ is irreducible if and only if
$w$ belongs to the same coset in $W'\setminus W$ as the longest element
$w_0$ of $W$.
\end{theoremintro}

An essential part of the approach from \cite{MiSo} is the study
(including an answer in special cases) of the so-called {\em
Kostant's problem} from \cite{Jo}. If $M$ is a $\mathfrak{g}$-module
and $\mathrm{Ann}(M)$ is the annihilator of $M$ in
$U(\mathfrak{g})$, then  the vector space
$U(\mathfrak{g})/\mathrm{Ann}(M)$ canonically embeds into the vector
space of all $\mathbb{C}$-linear automorphisms of $M$, which are
locally finite with respect to the adjoint action of $\mathfrak{g}$.
The question, which was called {\em Kostant's problem for $M$} in
\cite{Jo}, is to determine for which modules $M$ the canonical
injection above is in fact an isomorphism. We answer this question
for several modules $M$. In particular, we prove the following
statement:

\begin{theoremintro}\label{thmintro5}
Let $\mathfrak{g}$ be of type $A$, and $x,y\in W=S_n$ be elements in the same
left cell. Then Kostant's problem has a positive answer for the simple highest
weight module $L(x\cdot 0)$ if and only if it has a positive answer for the
simple highest weight module $L(y\cdot 0)$.
\end{theoremintro}

We believe that our approach will finally lead to a complete positive answer
of this problem in type $A$. For other types our approach fails, but the answer
to Kostant's problem is negative as well (as was shown in \cite{Jo}). A
detailed analysis of the problem, our partial solutions, and the
obstacles for other types are given in Subsection~\ref{s9.2}.

On the way to our main results we also obtain several categorification results
which we think are of interest on their own. We also obtain some unexpected
applications of the categorification procedure, in particular we define
a canonical filtration on integral permutation and induced cell modules
for the symmetric group $S_n$.
\smallskip

\noindent
{\bf A structural overview.}
The paper starts with a general discussion on the notion of {\em
precategorification} and {\em categorification} in Section~\ref{s2}.
In Section~\ref{s25} we give a brief summary of the known categorifications
of the regular representations of
the Hecke algebra. A categorification of the Kazhdan-Lusztig (right) cell
modules is given in Section~\ref{s3}. The categories appearing in this
categorification are not very well understood. They are defined as quotients
of certain subcategories of the category $\cO$. In general they are not
highest weight categories and have infinite homological dimension
(see Subsection~\ref{s4.4}).

From our uniqueness result it follows that the categorifications of the cell
modules are certain module categories over (in general non-commutative)
symmetric algebras including as a special case Khovanov's algebra
$\mathcal{H}^n$ (from \cite{Kh}, \cite{St3}); and the uniqueness result
together with \cite{Br}, \cite{St2} shows that the centres of these categories
are isomorphic to the cohomology ring of a certain Springer fiber, that means
the fixed point variety of the flag variety $GL(n,\mC)/B$ under a nilpotent
matrix $N$ (in Jordan normal form). For the standard examples of induced
modules, namely the induced sign or induced trivial module, categorifications
in terms of parabolic category $\cO$ (see e.g. \cite{SoKipp}, \cite{StDuke})
and Harish-Chandra bimodules (see e.g. \cite{MS}) are well-known and
well-studied. Our categories corresponding to induced modules are
generalizations of both, parabolic  category $\cO$ and certain categories of
Harish-Chandra bimodules. A short summary can be found in Section~\ref{s5}.

In an induced cell module for the Iwahori-Hecke algebra, we have four special
bases which will have a very natural categorical interpretation
(Theorem~\ref{combinatorics}) in terms of isomorphism classes of projective
modules, simple modules, standard modules (which are induced projective modules
from the categorification of the cell module) and proper standard modules
(which are induced simple modules). The categorifications for induced cell
modules are stratified in the sense of \cite{CPS}, and even weakly properly
stratified in the sense of \cite{Fr}. The latter structure plays
an important role in several parts of the paper.

In Section~\ref{s8} we study properties of the categories used to categorify
induced cell modules (in type $A$). Generalizing Irving's results from
\cite{Irself}, we classify all projective modules which are also injective
(Theorem~\ref{irving}) and then deduce a double centralizer property
(Theorem~\ref{pr991}) which generalizes Soergel's original Struktursatz
from \cite{Sperv} highly non-trivially. Maybe the most surprising result here
is the description of the center of these induced categories
(Theorem~\ref{pcentre}): the center is isomorphic to the center of a certain
parabolic category $\cO$. Therefore, we again have the explicit description of
the center as given in \cite{Br}. Moreover, the categories categorifying
induced cell modules are all Ringel self-dual (Theorem~\ref{trsd}),
which means that there is an equivalence between the additive subcategory of
all projective modules and the additive subcategory of all tilting modules.

The categorifications of induced cell modules will finally be used to describe
the best possible general result about generalized Verma modules, that means
parabolically induced {\it arbitrary} simple modules. The generalized Verma
modules as briefly explained above appear as the so-called (proper) standard
objects in our categorifications. Our combinatorial description can then be
used to deduce at least the multiplicity of certain composition factors
(namely the one which can be seen in our categories), and leads to what is
called the `rough structure' of generalized Verma modules. In this rough
structure all the multiplicities become finite. A very special case of our
setup was already considered in \cite{MiSo} and \cite{KM2}.
\smallskip

\noindent
{\bf General terminology.} A {\it ring} always means an associative
unitary ring. {\it Graded} always means $\mathbb{Z}$-graded. For a ring $R$
we denote by $R\operatorname{-mod}$ and $\operatorname{mod-}R$ the categories
of finitely generated left and right $R$-modules respectively. If $R$ is graded
then we denote by $R\gmod$ and $\operatorname{gmod-}R$ the categories of
finitely generated graded left and right $R$-modules respectively. Inclusions
are denoted by $\subset$. If it is necessary to point out that some inclusion
is proper we use the symbol $\subsetneq$.

Let $\mathbb{F}$ be a commutative ring. We denote by $\mathbb{F}[v,v^{-1}]$
and $\mathbb{F}((v))$ the rings of Laurent polynomials and formal Laurent
series in the variable $v$ with coefficients in $\mathbb{F}$ respectively.
In the paper we usually work over $\mathbb{Z}$  or over $\mathbb{C}$. We
abbreviate $\otimes_{\mathbb{C}}$ as $\otimes$.
\smallskip

\noindent {\bf Acknowledgments.} We would like to thank Henning
Haahr Andersen, Roman Bezrukavnikov, Oleksandr Khomenko, Ryszard
Rubinsztein and David Vogan for useful and stimulating discussions.

\section{Decategorification, precategorification and
categorification}\label{s2}

In this section we define a general algebraic notion of
categorification. The definition is based on and further develops
the ideas of \cite{KMS}, \cite{KMS2}, \cite{MS}.

\subsection{Ordinary setup}\label{s2.1}

Let $\ccC$ be a category. If $\ccC$ is abelian or triangulated, we denote by
$\op{Gr}(\ccC)$ the {\it Grothendieck group} of $\ccC$. The latter one is
by definition the free abelian group generated by the isomorphism classes
$[M]$ of objects $M$ of $\ccC$ modulo the relation $[C]=[A]+[B]$ whenever
there is a short exact sequence $A\hookrightarrow C\tto B$ if $\ccC$ is abelian;
and whenever there is a triangle $(A,C,B,f,g,h)$ if $\ccC$ is triangulated.
If $\ccC$ is additive, we denote by $\op{Gr}(\ccC)_{\oplus}$
the {\em split Grothendieck group} of $\ccC$, which is by definition the free
abelian group generated by the isomorphism classes $[M]$ of objects $[M]$ of
$\ccC$ modulo the relation $[C]=[A]+[B]$ whenever  $C\cong A\oplus B$. For
$M\in\ccC$ we denote by $[M]$ the image of $M$ in the (split) Grothendieck
group.  Let $\mathbb{F}$ be  a commutative ring with $1$.

\begin{definition}{\rm
Let $\ccC$ be an abelian or triangulated, respectively additive, category.
Then the  {\em $\mathbb{F}$-decategorification} of $\ccC$ is the
$\mathbb{F}$-module $[\ccC]^{\mathbb{F}}:= \mathbb{F}\otimes_{\mathbb{Z}}
\op{Gr}(\ccC)$ (resp. $[\ccC]^{\mathbb{F}}_{\oplus}:=
\mathbb{F}\otimes_{\mathbb{Z}}  \op{Gr}(\ccC)_{\oplus}$). }
\end{definition}

The element $1\otimes [M]$ of the $\mathbb{F}$-decategorification is
abbreviated as $[M]$ as well. We set $[\ccC]:=[\ccC]^{\mathbb{Z}}$ and
$[\ccC]_{\oplus}:= [\ccC]^{\mathbb{Z}}_{\oplus}$.

\begin{definition}{\rm
Let $V$ be an $\mathbb{F}$-module. An {\em $\mathbb{F}$-precategorification}
$(\ccC,\varphi)$ of $V$ is an abelian (resp. triangulated or additive)
category $\ccC$ with a fixed monomorphism $\varphi$ from $V$ to the $\mathbb{F}$-decategorification of $\ccC$.  If $\varphi$ is an isomorphism,
then $(\ccC,\varphi)$ is called an {\em $\mathbb{F}$-categorification}
of $V$.
}
\end{definition}

Hence categorification is in some sense the `inverse' of
decategorification. Whereas the latter is uniquely defined, there
are usually several different categorifications. In case
$\mathbb{F}=\mathbb{Z}$ and $V$ is torsion-free there is always the
(trivial) categorification given by a semisimple category of the
appropriate size.

\begin{definition}{\rm
Let $V$ be an $\mathbb{F}$-module and $f:V\to V$ be an
$\mathbb{F}$-endomorphism. Given an $\mathbb{F}$-precategorification
$(\ccC,\varphi)$ of $V$, an {\em $\mathbb{F}$-categorification} of $f$ is an
exact (resp. triangulated or additive) functor $F:\ccC\to \ccC$ such that
$[F]\circ \varphi=\varphi\circ f$, where
$[F]$ denotes the endomorphism of $[\ccC]^{\mathbb{F}}$
(or $[\ccC]_{\oplus}^{\mathbb{F}}$ if $\ccC$ is abelian)
induced by $F$. In other words, the following diagram commutes:
\begin{displaymath}
\xymatrix{
V\ar[rr]^{f}\ar[d]_{\varphi}&&
V\ar[d]^{\varphi}\\
[\ccC]_{(\oplus)}^{\mathbb{F}}\ar[rr]^{[F]} && [\ccC]_{(\oplus)}^{\mathbb{F}}
}
\end{displaymath}
}
\end{definition}

\begin{definition}
{\rm
Assume $A$ is some $\mathbb{F}$-algebra defined by generators $a_1,\dots,a_k$
and relations $R_j$, $j\in J$. Given an $A$-module $M$, each generator $a_i$ of
$A$ defines a linear endomorphism, $f_i$, of $M$. A {\em very weak
$\mathbb{F}$-(pre)categorification} of $M$ is a (pre)categorification
$(\ccC,\varphi)$ of the vector space $M$ together with a
categorification $F_i$, $i=1,\dots,k$, of each $f_i$.}
\end{definition}

If there is an `interpretation' of the relations $R_j$ between the
generators of $A$ in terms of isomorphisms of functors, we will call
$(\ccC,\varphi,F_1,\dots,F_k)$ a {\em (pre)ca\-te\-go\-ri\-fi\-ca\-ti\-on}
of the  $A$-module $M$. The interpretation of the relations will depend on
the example.

\begin{example}{\rm
Let $R=\mC[x]/(x^2)$ and $\ccC=R\dmod$. Then $\op{Gr}(\ccC)\cong\mZ$,
generated by the isomorphism class $[\mC]$ of the unique simple $R$-module,
and $[\ccC]\cong\mZ$. Thus $\ccC$ is a $\mZ$-categorification of $\mZ$.}
\end{example}

\subsection{Graded setup}\label{s2.2}

If $\ccC$ is equivalent to a category of modules over a graded ring,
then $\op{Gr}(\ccC)$ (or $\op{Gr}(\ccC)_{\oplus}$) becomes a
$\mathbb{Z}[v,v^{-1}]$-module via $v^{i}[M]=[M\langle i\rangle]$ for
any $M\in\ccC$, $i\in\mZ$, where $M\langle i\rangle$ is the module
$M$, but in the grading shifted by $i$ such that $(M\langle
i\rangle)_j=M_{j-i}$.

To define the notion of a decategorification for a category of
graded modules (or complexes of graded modules) let $\mathbb{F}$ be
a commutative ring with $1$ and $\iota:\mathbb{Z}[v,v^{-1}]\to
\mathbb{F}$ be a fixed homomorphism of unitary rings. Then $\iota$ defines on
$\mathbb{F}$ the structure of a (right)
$\mathbb{Z}[v,v^{-1}]$-module.

\begin{definition}{\rm
The {\em $(\mathbb{F},\iota)$-decategorification} of $\ccC$ is
the $\mathbb{F}$-module
\begin{displaymath}
[\ccC]^{(\mathbb{F},\iota)}:=
\mathbb{F}\otimes_{\mathbb{Z}[v,v^{-1}]} \op{Gr}(\ccC)\quad
(\text{resp. }[\ccC]^{(\mathbb{F},\iota)}_{\oplus}:=
\mathbb{F}\otimes_{\mathbb{Z}[v,v^{-1}]} \op{Gr}(\ccC)_{\oplus}).
\end{displaymath}
}
\end{definition}

In most of our examples the homomorphism $\iota:\mathbb{Z}[v,v^{-1}]\to
\mathbb{F}$ will be the obvious canonical inclusion. In such cases
we will omit $\iota$ in the notation. We set
\begin{displaymath}
[\ccC]:=[\ccC]^{(\mathbb{Z}[v,v^{-1}],\mathrm{id})},\quad
[\ccC]_{\oplus}:=[\ccC]^{(\mathbb{Z}[v,v^{-1}],\mathrm{id})}_{\oplus}.
\end{displaymath}

\begin{definition}{\rm
Let $V$ be an $\mathbb{F}$-module. A {\em $\iota$-precategorification}
$(\ccC,\varphi)$ of $V$ is an abelian or triangulated, respectively additive,
category $\ccC$ with a fixed free action of $\mathbb{Z}$ and a fixed
monomorphism $\varphi$ from $V$ to the $(\mathbb{F},\iota)$-decategorification
of $\ccC$. If $\varphi$ is an isomorphism, $(\ccC,\varphi)$ is called a
{\em $\iota$-categorification} of $V$.
}
\end{definition}

The definitions of a {\em $\iota$-categorification} of an endomorphism,
$f:V\to V$, and of a {\em $\iota$-(pre)categorification} of a module over
some $\mathbb{F}$-algebra are completely analogous to the corresponding
definitions from the previous subsection.

\begin{example}{\rm
Let $R=\mC[x]/(x^2)$. Consider $R$ as a graded ring (we usually consider it
as a cohomology ring and put $x$ in degree two), and take $\ccC=R\gmod$.
Then $[\ccC]\cong\mZ[v,v^{-1}]$ as a $\mZ[v,v^{-1}]$-module, hence the  graded
category $\ccC$ is a $(\mZ[v,v^{-1}],\mathrm{id})$-categorification  of
$\mZ[v,v^{-1}]$. Note that $\ccC$, considered just as an abelian category,
is also a $\mZ[v,v^{-1}]$-ca\-te\-go\-ri\-fi\-ca\-ti\-on of $\mZ[v,v^{-1}]$.
}
\end{example}

\section{The Hecke algebra as a bimodule over itself
and its categorifications}\label{s25}

In this section we recall the definition of Hecke algebras and give
several examples of categorifications of regular (bi)modules over
these algebras. We refer the reader to \cite{KMS2} for more examples
of categorifications.

 From now on we fix a finite Weyl group $W$
with identity element $e$, set of simple reflections $S$, and length
function $l$. Denote by $w_0$ the longest element of $W$. Let
further $\leq$ be the Bruhat order on $W$. With respect to this
order the element $e$ is the minimal and $w_0$ is  the maximal
element. Our main example will be $W=S_n$, the symmetric group  on
$n$ elements, and $S=\{(i,i+1),i=1,\dots,n-1\}$, the set of
elementary transpositions.

\subsection{The Hecke algebra}\label{s25.2}

Denote by $\mathds{H}=\mathds{H}(W,S)$ the
{\it Hecke algebra} associated  with $W$ and $S$; that is the
$\mathbb{Z}$-algebra which is a free $\mathbb{Z}[v,v^{-1}]$-module with
basis $\{H_x\mid x\in W\}$ and multiplication given by
\begin{equation} \label{eqhecke}
H_xH_y=H_{xy}\, \text{if $l(x)+l(y)=l(xy),\,\,$ and}\,\,
H_s^2=H_e+(v^{-1}-v)H_s\, \text{for $s\in S$}.
\end{equation}
The algebra $\mathds{H}$ is a deformation of the group  algebra
$\mathbb{Z}[W]$. As a $\mathbb{Z}[v,v^{-1}]$-algebra it is generated by
$\{H_s\mid s\in S\}$, or (which will turn out to be more convenient) by the set
$\{\underline{H}_s=H_s+vH_e\mid s\in S\}$. Note that $\underline{H}_s$ is fixed
under the involution ${}^-$, which maps $v\mapsto v^{-1}$ and  $H_s\mapsto
(H_s)^{-1}$, i.e, $\underline{H}_s$ is a Kazhdan-Lusztig basis element. More
general, for $w\in W$ we denote by $\underline{H}_w$ the corresponding element
from the Kazhdan-Lusztig bases for $\mathds{H}$ in the normalization of
\cite{SoKipp}. The Kazhdan-Lusztig polynomials $h_{x,y}\in \mathbb{Z}[v]$ are
defined via $\underline{H}_x=\sum_{y\in W}h_{y,x}H_x$. With respect to the
generators $\underline{H}_s$, $s\in S$, we have the following set of defining
relations (in the case $W=S_n$):
\begin{eqnarray} \label{eqhecke2}
\underline{H}_s^2&=&(v+v^{-1})\underline{H}_s;\\ \nonumber
\underline{H}_s\underline{H}_t&=&\underline{H}_t\underline{H}_s,
\quad\quad\quad\quad\, \text{ if } ts=st;\\ \nonumber
\underline{H}_s\underline{H}_t\underline{H}_s+
\underline{H}_t&=&\underline{H}_t\underline{H}_s\underline{H}_t
+\underline{H}_s, \,\, \text{ if }\, ts\neq st.
\end{eqnarray}

Let $\mathbb{F}$ be any commutative ring and $\iota:\mathbb{Z}[v,v^{-1}]\to
\mathbb{F}$ be a homomorphism of unitary rings. Then we have the
{\em specialized} Hecke algebra
$\mathds{H}^{(\mathbb{F},\iota)}=\mathbb{F}\otimes_{\mathbb{Z}[v,v^{-1}]}
\mathds{H}$. Again if $\iota$ is clear from the context (for
instance if $\iota$ is the natural inclusion), we will omit it in
the notation.

\begin{example}\label{ex6}
{\rm
Let again $R=\mC[x]/(x^2)$. Putting $x$ in degree two, induces
a grading on $R$ and $B_s:=R\langle-1\rangle$ becomes a graded $R$-bimodule.
Let $\ccS$ be the additive category generated by the graded left $R$-modules
$\mC\langle j\rangle$ and $B_s\langle j\rangle$, $j\in\mZ$. Then
$\op{Gr}(\ccS)_\oplus$ is a free $\mZ[v,v^{-1}]$-module of rank two, and is
isomorphic to $\mathds{H}(S_2, \{s\})$ via $[\mC\langle l\rangle]\mapsto
v^{-l}H_e$, $[B_s\langle l\rangle]\mapsto v^{-l}\underline{H}_s$. The functor
$F_s^l=B_s\otimes_R -$ satisfies the condition $F_s^l\circ F_s^l\cong
F_s^l\langle 1\rangle\oplus F_s^l\langle -1\rangle$ which is an interpretation
of the first relation in \eqref{eqhecke}. Hence we get a categorification of
the left regular $\mathds{H}(S_2, \{s\})$-module. This example generalizes to
arbitrary finite Weyl groups as we will describe in the next subsection.
}
\end{example}

\subsection{Special bimodules}\label{s25.3}

Associated with $W$ we have the additive category given by the so-called
{\em special bimodules} $B_w$, $w\in W$, introduced by Soergel in \cite{SHC},
see also \cite{SKLP}. To define these bimodules we consider the geometric
representation $(V_{\mathbb{R}},\varphi)$  of $W$ and its complexification
$(V,\varphi)$, see \cite[4.2]{BjBr}. Let $R$ be the ring of regular functions
on $V$ with its natural $W$-action. This ring becomes graded by putting $V^*$
in degree $2$. For any $s\in S$ let $R^s$ be the subring of $s$-invariants in
$R$. Note that this is in fact a graded subring of $R$. Given $w\in W$ with a
fixed reduced expression $[w]=s_1s_2\cdot\ldots\cdot s_k$ define the graded
$R$-bimodule $R_{[w]}$ as follows:
\begin{displaymath}
R_{[w]}=R\otimes_{R^{s_1}}R\otimes_{R^{s_2}}\dots
\otimes_{R^{s_k}}R \langle -l(w)\rangle.
\end{displaymath}
Following \cite{SKLP} we define $B_w$ as the unique indecomposable direct
summand of $R_{[w]}$, which is not isomorphic to a direct summand of any
$R_{[x]}$ with $l(x)<l(w)$. Let $\ccS$ be the smallest additive category which
contains all special bimodules, and is closed under taking direct sums and
graded shifts. There is a unique isomorphism $\mathcal{E}$ of
$\mathbb{Z}[v,v^{-1}]$-modules, which satisfies
\begin{eqnarray*}
\mathcal{E}: \quad \mathds{H} & \overset{\sim}{\longrightarrow}
& \left[\ccS\right]_{\oplus} \\ \underline{H}_w & \mapsto & \left[B_w\right].
\end{eqnarray*}
For any $s\in S$ we have the additive endofunctors $\mathrm{F}_s^l=
B_s\otimes_{R}{}_-$ and $\mathrm{F}_s^r={}_-\otimes_{R}B_s$ of $\ccS$.
Altogether we get a categorification of the regular Hecke module as follows
(see \cite[Theorem~1.10]{SKLP}, \cite[Satz~7.9]{Ha} and \cite[Theorem~1]{SHC}):

\begin{proposition}\label{prop1}
\begin{enumerate}[(i)]
\item\label{prop1.1} $(\ccS,\mathcal{E},\{\mathrm{F}_s^r\}_{s\in S})$ is
a categorification of the right regular representation of $\mathds{H}$
with respect to the generators $\underline{H}_s$, $s\in S$.
\item\label{prop1.2} $(\ccS,\mathcal{E},\{\mathrm{F}_t^l\}_{t\in S})$ is
a categorification of the left regular representation of $\mathds{H}$
with respect to the generators $\underline{H}_t$, $t\in S$.
\end{enumerate}
\end{proposition}

The interpretation of the relations \eqref{eqhecke2} is given by the existence
(see \cite[Theorem~1]{SHC}) of isomorphisms of functors as follows
(in case $W=S_n$):
\begin{eqnarray*}
(\mathrm{F}_s^{\sharp})^2&\cong&\mathrm{F}_s^{\sharp}\langle 1\rangle
\oplus\mathrm{F}_s^{\sharp}\langle-1\rangle;\\ \nonumber
\mathrm{F}_s^{\sharp}\mathrm{F}_t^{\sharp}&\cong&
\mathrm{F}_t^{\sharp}\mathrm{F}_s^{\sharp},
\quad\quad\quad\quad\, \text{ if } ts=st;\\ \nonumber
\mathrm{F}_s^{\sharp}\mathrm{F}_t^{\sharp}\mathrm{F}_s^{\sharp}\oplus
\mathrm{F}_t^{\sharp}&\cong&\mathrm{F}_t^{\sharp}
\mathrm{F}_s^{\sharp}\mathrm{F}_t^{\sharp}
\oplus\mathrm{F}_s^{\sharp}, \,\, \text{ if }\, ts\neq st,
\end{eqnarray*}
where $\sharp$ is either $l$ or $r$. For other types the interpretation is
similar.

\begin{remark}\label{rem1}
{\rm
\begin{enumerate}
\item\label{rem1.1} The functors $\mathrm{F}_s^r$ and $\mathrm{F}_t^l$ naturally
commute (with each other),  hence the parts \eqref{prop1.1} and \eqref{prop1.2}
of  Proposition~\ref{prop1} together give a categorification of the
regular Hecke {\bf bimodule}.
\item\label{rem1.2} The above categorification is not completely satisfactory,
mostly because it is given by an additive category which is not abelian.
As a consequence, we cannot see the standard basis of the Hecke module in
this categorification, hence we will present a categorification given by an
abelian category. This will be done in the next subsection.
\item\label{rem1.3} The proof of Proposition~\ref{prop1} given
in \cite[Theorem~1]{SHC} is quite involved and uses the full power of the
Kazhdan-Lusztig Theory (or the decomposition theorem
\cite[Theorem~6.2.5]{BBD}).
\item\label{rem1.4} If one prefers to work with finite-dimensional algebras and
modules, one could replace the polynomial ring $R$ with the coinvariant ring
$C$, which is the quotient of $R$ modulo the ideal generated by homogeneous
$W$-invariant polynomials of positive degree. One can define the {\em special
$C$-bimodules} $B_w\otimes_R C$ and obtains a completely analogous result
to Proposition~\ref{prop1}, and Remark (1), see \cite[Theorem~2]{SHC}.
\item\label{rem1.5}
We could also define the {\em special right $C$-modules}
$\overline{B}_w=\mathbb{C}\otimes_R B_w\otimes_R C$.
Since they are preserved by the functors $\mathrm{F}_s^r$, $s\in S$,
Proposition~\ref{prop1}\eqref{prop1.1} provides another
categorification of the regular right $\mathds{H}$-module,
\cite[Zerlegungssatz~1 and Section~2.6]{Sperv}.
\end{enumerate}
}
\end{remark}

\subsection{Harish-Chandra bimodules}\label{s25.4}

In this section we would like to improve  Proposition~\ref{prop1} and work
with abelian categories. We start with introducing the setup, which then will
also be used in the next subsection.

Let $\mathfrak{g}$ be a reductive finite-di\-men\-sio\-nal complex  Lie
algebra associated with the Weyl group $W$. Let $U(\g)$ be the universal
enveloping algebra of $\g$ with its center $Z(\g)$. Fix a  triangular decomposition $\g=\n_-\oplus\h\oplus\n_+$, where $\h$ is a fixed Cartan
subalgebra of $\g$ contained in the Borel subalgebra  $\mathfrak{b}=
\h\oplus\n_+$. For $\lambda\in \h^*$ we denote by $M(\lambda)$ the Verma
module with highest weight $\lambda$. Let $\rho$ be the half-sum of all
positive roots. Define $\h^*_{dom}:=\{\lambda\in\h^*\,:\,\lambda+\rho
\text{ is dominant}\}$, which is the dominant Weyl chamber with
respect to the {\em dot-action} of $W$ on $\h^*$ given by $w\cdot
\lambda=w(\lambda+\rho)-\rho$.

Denote by $\mathcal{H}$ the category of Harish-Chandra bimodules for $\g$,
that is the category of finitely generated $U(\g)$-bimodules of finite
length, which are locally finite with respect to the adjoint action of $\g$
(which is defined for a bimodule $M$ as $x.m=xm-mx$ for any $x\in\mg$ and
$m\in M$). The action of the center defines the following block decomposition
of $\mathcal{H}$:
\begin{displaymath}
\mathcal{H}=\bigoplus_{\mmm,\mmn\in \text{Max}Z(\g)}
{}_{\mmm}\mathcal{H}_{\mmn},
\quad\text{where}\quad
{}_{\mmm}\mathcal{H}_{\mmn}=
\left\{ M\in \mathcal{H}| \exists\, k\in \N:\mmm^k M=0=M\mmn^k
\right\}.
\end{displaymath}
Note that $Z(\g)\cong R$ (via the Harish-Chandra isomorphism and \cite[18-1]{Kane}) hence it is positively graded. Let $\mathbf{0}\in
\text{Max}Z(\g)$ denote the annihilator (in $Z(\g)$) of the trivial
$U(\g)$-module. Consider the block ${}_{\mathbf{0}}\mathcal{H}_{\mathbf{0}}$.
Tensoring with finite-dimensional left and right $U(\g)$-modules are
endofunctors on $\mathcal{H}$ and their direct sums and summands are
called {\em projective functors}. Indecomposable projective functors were
classified in  \cite[Theorem~3.3]{BG}. It turns out
that these summands are naturally labelled by the elements of $W$. For $w\in W$
we denote by $\theta_w^l$ the indecomposable projective endofunctor of
${}_{\mathbf{0}}\mathcal{H}_{\mathbf{0}}$ corresponding to $w$ and induced by tensoring with a finite dimensional {\it left} $\mg$-module (as the supindex $l$ indicates). Similarly, we can consider projective functors given by
tensoring with finite-dimensional {\it right} $U(\g)$-modules and obtain the
corresponding functors $\theta_w^r$. For two $\g$-modules $M$ and $N$ we
denote by $\mathscr{L}(M,N)$ the largest $\mathrm{ad}(\g)$-finite submodule of
$\mathrm{Hom}_{\mathbb{C}}(M,N)$, see \cite[1.7.9]{Di}. The classes
$[\mathscr{L}(M(0),M(w\cdot 0))]$, $w\in W$, form a basis of
$[{}_{\mathbf{0}}\mathcal{H}_{\mathbf{0}}]$, see \cite{BG}, \cite[6.15]{Ja2}.

Following \cite[Theorem~2]{SHC} we form the positively graded algebra
\begin{displaymath}
A^{\infty}=\mathrm{End}_{R-R}(\oplus_{w\in W}B_w)
\end{displaymath}
and we have an equivalence (see \cite[Theorem~3]{SHC}) of categories
\begin{displaymath}
{}_{\mathbf{0}}\mathcal{H}_{\mathbf{0}}\cong \mathrm{nil-}A^{\infty},
\end{displaymath}
where $\mathrm{nil-}A^{\infty}$ is the category of all finite dimensional
right $A^{\infty}$-modules $M$ satisfying $M\,A^{\infty}_i=0$ for all $i\gg 0$
(for example, this is obviously satisfied for any finite dimensional
{\em gradable} $A^{\infty}$-module).
We consider the category $\mathrm{gmod-}A^{\infty}$ of all finite-dimensional
graded right $A^{\infty}$-modules. The functors $\theta_w^l$ and $\theta_w^r$
lift to endofunctors of $\mathrm{gmod-}A^{\infty}$, see \cite[Appendix]{MO}.
The modules $\mathscr{L}(M(0),M(w\cdot 0))$ admit graded lifts as well
and we fix standard lifts $\tilde{M}_w$ such that their heads are
concentrated in degree $0$. Let $\tilde{\mathcal{E}}$ be the unique
isomorphism of the $\mathbb{Z}[v,v^{-1}]$-modules such that
\begin{eqnarray*}
\tilde{\mathcal{E}}: \quad \mathds{H} & \overset{\sim}{\longrightarrow} &
\left[\mathrm{gmod-}A^{\infty}\right] \\
H_w & \mapsto &\left[\tilde{M}_w\right] .
\end{eqnarray*}

\begin{proposition}\label{prop2}
\begin{enumerate}[(i)]
\item\label{prop2.1} $(\mathrm{gmod-}A^{\infty},\tilde{\mathcal{E}},
\{\theta_s^l\}_{s\in S})$ is a categorification of the right regular
representation of $\mathds{H}$ with respect to the generators
$\underline{H}_s$, $s\in S$.
\item\label{prop2.2} $(\mathrm{gmod-}A^{\infty},\tilde{\mathcal{E}},
\{\theta_t^r\}_{t\in S})$ is a categorification of the left regular
representation of $\mathds{H}$ with respect to the generators
$\underline{H}_t$, $t\in S$.
\end{enumerate}
\end{proposition}

This statement can be found for example in \cite{StThesis} and
\cite{Kh2}. Basically, it follows from \cite{SHC}. We would like to
emphasize the difference between
Proposition~\ref{prop2} and Proposition~\ref{prop1}: In
Proposition~\ref{prop1} the {\em right} regular representation of
$\mathds{H}$ was categorified using functors $\mathrm{F}_s^r$ of
tensor product from the {\em right}, while in
Proposition~\ref{prop1} the {\em right} regular representation of
$\mathds{H}$ was categorified using the {\em left} translation
functors $\theta_s^l$. The interpretation of the relations
\eqref{eqhecke2} is similar to the one given after
Proposition~\ref{prop1}. Again, the functors $\theta_s^l$ and
$\theta_t^r$ naturally commute with each other and hence the parts
\eqref{prop2.1} and \eqref{prop2.2} of  Proposition~\ref{prop2}
together give a categorification of the regular Hecke bimodule. The
connection to Subsection~\ref{s25.3} is given by
\cite[Section~3]{SHC}.

\subsection{Category $\mathcal{O}$}\label{s25.5}

We stick to the setup at the beginning of the previous subsection. Consider
the BGG category $\mathcal{O}=\mathcal{O}(\g,\mathfrak{b})$ (\cite{BGG2})
with its block decomposition
\begin{displaymath}
\mathcal{O}=\bigoplus_{\lambda\in \h^*_{dom}}\mathcal{O}_{\lambda},
\quad\text{where}\quad
\mathcal{O}_{\lambda}=\left\{ M\in \mathcal{O}| \exists\, k\in
\N:(\mathrm{Ann}_{Z(\g)}(M(\lambda)))^k M=0 \right\}.
\end{displaymath}
For $\lambda\in\h^*$ let $P(\lambda)$ be the projective cover of
$M(\lambda)$ and $L(\lambda)$ be the simple quotient of $P(\lambda)$.

Following \cite{Sperv} we form the graded algebra
$A=\mathrm{End}_{C}(\oplus_{w\in W}\overline{B}_w)$ (see
Remark~\ref{rem1}) and obtain an equivalence of categories between
$\mathcal{O}_0$ and $\mathrm{mod-}A$, the category of finite
dimensional right $A$-modules. We denote by
$\mathcal{O}_0^{\mathbb{Z}}$ the category of finite-dimensional {\em
graded} right $A$-modules. To connect this with the previous
subsection let ${}_{\mathbf{0}}\mathcal{H}_{\mathbf{0}}^1$ denote
the full subcategory of ${}_{\mathbf{0}}\mathcal{H}_{\mathbf{0}}$
consisting of all bimodules which are annihilated by $\mathbf{0}$
from the right hand side. Then there is an equivalence of categories
${}_{\mathbf{0}}\mathcal{H}_{\mathbf{0}}^1\cong \mathcal{O}_0$, see
\cite[Theorem~5.9]{BG}. Via this equivalence the functors $\theta_w^l$,
$w\in W$, restrict to exact endofunctors of $\mathcal{O}_0$, which admit
graded lifts. Unfortunately, the functors $\theta_w^r$ do not preserve
$\mathcal{O}_0$. However, for $s\in S$ there is a unique up to scalar natural
transformation $\mathrm{Id}\langle 1\rangle\to \theta_s^r$, whose cokernel we
denote by $\mathrm{T}_s$. These are the so-called {\em twisting functors} on
$\mathcal{O}_0$, see \cite{AS} and \cite{KM}. Each $\mathrm{T}_s$ preserves
$\mathcal{O}_0$ and has a graded lift by definition, but it is only right
exact. Therefore we consider $\mathcal{D}^b(\mathcal{O}_0^{\mathbb{Z}})$, the
bounded derived category of the category of finite-dimensional graded right
$A$-modules with shift functor $\llbracket\cdot \rrbracket$. Let
$\mathcal{L}\mathrm{T}_s$ be the left derived functor of $T_s$.

For $w\in W$ we abbreviate $\Delta(w)=M(w\cdot 0)$, $L(w)=L(w\cdot 0)$ and
$P(w)=P(w\cdot 0)$. All simple, standard and projective modules in
$\mathcal{O}_0$ have standard graded lifts (i.e. their heads are concentrated
in degree zero), which we will denote by the same symbols. We fix the unique
isomorphism of the  $\mathbb{Z}[v,v^{-1}]$-modules such that
\begin{eqnarray*}
\hat{\mathcal{E}}: \quad \mathds{H} & \overset{\sim}{\longrightarrow} & \left[\mathcal{D}^b(\mathcal{O}_0^{\mathbb{Z}})\right]
\\ H_w & \mapsto & \left[\Delta(w)\right]
\end{eqnarray*}
and obtain the following well-known result:

\begin{proposition}\label{prop3}
\begin{enumerate}[(i)]
\item\label{prop3.1} $(\mathcal{D}^b(\mathcal{O}_0^{\mathbb{Z}}),
\hat{\mathcal{E}},\{\theta_s^l\}_{s\in S})$ is a categorification of
the  right regular representation of $\mathds{H}$ with respect to the
generators $\underline H_s$, $s\in S$.
\item\label{prop3.2} $(\mathcal{D}^b(\mathcal{O}_0^{\mathbb{Z}}),
\hat{\mathcal{E}}, \{\mathcal{L}T_t\}_{t\in S})$  is a categorification of
the left regular representation of $\mathds{H}$ with respect to the
generators $H_t$, $t\in S$.
\item\label{prop3.3} $\hat{\mathcal{E}}(\underline{H}_w)=[P(w)]$
for all $w\in W$.
\end{enumerate}
\end{proposition}

\begin{proof}
The ungraded resp. graded cases of \eqref{prop3.1} are treated in
\cite[Theorem~3.4(iv)]{BG} and \cite[Theorem~7.1]{Stgrad}.
The ungraded resp. graded cases of \eqref{prop3.2} follow from
\cite[(2.3) and Theorem~3.2]{AS} and  \cite[Appendix]{MO}.
The claim \eqref{prop3.3} follows from
\cite[Theorem~3.11.4(i) and (iv)]{BGS}.
\end{proof}

For Proposition~\ref{prop3}\eqref{prop3.1}, the interpretation of the relations
\eqref{eqhecke2} is similar to the one given after Proposition~\ref{prop1}.
We note that the statements \eqref{prop3.1} and \eqref{prop3.3} of
Proposition~\ref{prop3} can be formulated entirely using the  underlying
abelian category,  whereas  the statement \eqref{prop3.2} can not. The
functors $\mathcal{L}T_t$, $t\in S$, satisfy braid relations (this can be
proved analogously to \cite[Proposition~11.1]{MS3} using \cite[Theorem~2.2]{AS}
and \cite[Section~6]{KM}), but we do not know any functorial interpretation
for the relation $H_s^2=H_e+(v^{-1}-v)H_s$. Hence (at least for the moment)
Proposition~\ref{prop3}\eqref{prop3.2} gives only a {\em very weak}
categorification of the left regular representation of $\mathds{H}$, but a
categorification (in the stronger sense) of the underlying representation
of the  braid group, see \cite{Rouquier}.

The functors $\theta_s^l$ and $\mathcal{L}T_t$ naturally commute with each
other and hence the parts \eqref{prop3.1} and \eqref{prop3.2} of
Proposition~\ref{prop3} together give a (very weak) categorification of the
regular Hecke bimodule. The connection to Remark~\ref{rem1}\eqref{rem1.5}
is given by Soergel's functor $\mathbb{V}$, see \cite{Sperv}.

\subsubsection{$\mathfrak{gl}_2$-example, the basis given by
standard modules}\label{s2.6}

Consider the case $W=S_2=\{e,s\}$. In this case the category $\mathcal{O}_0$
is equivalent to the category of finite-dimensional right $A$-modules, where
$A$ is the path algebra of the following quiver with relations:
\begin{displaymath}
\xymatrix{
s\ar@/^/[rr]^{\alpha} && e\ar@/^/[ll]^{\beta}
}, \quad\quad \alpha\beta=0.
\end{displaymath}
The algebra $A$ is graded with respect to the length of paths. The algebra
$A$ has a simple preserving duality and  hence the categories of
finite-dimensional right and left  $A$-modules are equivalent. Working with
left $A$-modules reflects better the natural $\mathfrak{gl}_2$-weight picture,
so we will use it. The category $A\mathrm{-mod}$ has $5$ indecomposable
objects,  namely,
\begin{displaymath}
\Delta(s)=L(s):\quad\xymatrix{
\mathbb{C}\ar@/^/[rr]^{0} && 0\ar@/^/[ll]^{0}
}\quad\quad\quad
L(e):\quad\xymatrix{
0\ar@/^/[rr]^{0} && \mathbb{C}\ar@/^/[ll]^{0}
}
\end{displaymath}
\begin{displaymath}
P(s):\quad\xymatrix{
\mathbb{C}\oplus\mathbb{C}
\ar@/^/[rr]^{
\text{\tiny $\left(\begin{array}{cc}1&0\end{array}\right)$}
} && \mathbb{C}
\ar@/^/[ll]^{
\text{\tiny $\left(\begin{array}{c}0\\1\end{array}\right)$}}
}\quad\quad\quad
\Delta(e)=P(e):\quad\xymatrix{
\mathbb{C}\ar@/^/[rr]^{0} && \mathbb{C}\ar@/^/[ll]^{\mathrm{id}}
}
\end{displaymath}
\begin{displaymath}
I(e):\quad\xymatrix{
\mathbb{C}\ar@/^/[rr]^{\mathrm{id}} && \mathbb{C}\ar@/^/[ll]^{0}
}.
\end{displaymath}
Let $f_s$ and $f_e$ denote the primitive idempotents of $A$ corresponding to the vertices $s$ and $e$. Then
the functor $\theta_s^l$ is given by tensoring with the
bimodule $Af_s\otimes_{\mathbb{C}} f_sA$, and the functor
$\mathrm{T}_s$ is given by tensoring with the bimodule $Af_sA$.
We have $\mathcal{L}_i\mathrm{T}_s=0$, $i>1$.
The values of $\theta_s^l$, $\mathrm{T}_s$ and
$\mathcal{L}_1\mathrm{T}_s$ on the indecomposable objects
from $\mathcal{O}_0^{\mathbb{Z}}$ are:
\begin{displaymath}
\begin{array}{|c||c|c|c|c|c|}
\hline
M               & L(s) & L(e) & P(s) & P(e) & I(e)\\
\hline
\hline
\theta_s^l M    & P(s)\langle -1\rangle& 0  & P(s)\langle -1\rangle\oplus
P(s)\langle 1\rangle & P(s) & P(s)\langle -2\rangle\\
\hline
\mathrm{T}_s M  & I(e) & 0    & P(s)\langle -1\rangle & L(s) &
I(e)\langle -1\rangle\\
\hline
\mathcal{L}_1\mathrm{T}_s M  & 0 & L(e)\langle 1\rangle & 0 & 0&
L(e)\langle 1\rangle\\
\hline
\end{array}
\end{displaymath}
There are several basis for the Grothendieck group, the {\it standard} choice is given by the isomorphism classes of the standard modules $\Delta(w)$, $w=e,s$. In this basis, the action of our functors is as follows:
\begin{displaymath}
\begin{array}{rcclccl}
\left[\theta_s^l \Delta(e)\right] &=&
v & \left[\Delta(e)\right]& +& &\left[\Delta(s)\right];\\
\left[\theta_s^l \Delta(s)\right] &=&
&\left[\Delta(e)\right]&+ &v^{-1}& \left[\Delta(s)\right];\\
\left[\mathcal{L}\mathrm{T}_s \Delta(e)\right] &=&
&&&&\left[\Delta(s)\right]; \\
\left[\mathcal{L}\mathrm{T}_s \Delta(s)\right] &=&
&\left[\Delta(e)\right]&+&(v^{-1}-v)& \left[\Delta(s)\right].\\
\end{array}
\end{displaymath}
This is a categorification of the regular Hecke bimodule in the standard basis.

\subsubsection{$\mathfrak{gl}_3$-example, the basis given by  simple modules}
\label{s2.7}

Consider the case $W=S_3=\{e,s,t,st,ts,sts=tst=w_0\}$. In this case the
category $\mathcal{O}_0$ has infinitely many indecomposable objects (and is
in fact wild, see \cite{FNP}). However one can still compute the actions of
$\theta_s^l$, $\theta_t^l$, $\mathrm{T}_s$ and $\mathrm{T}_t$ in various bases
using known properties of these functors. The easiest basis is given by
standard modules; here, however, we present the answer for $\theta_s^l$,
$\theta_t^l$ in the most natural basis, namely the one given by simple modules.
To shorten the notation we will denote our simple modules just by the
corresponding elements of the Weyl group. Here are the graded filtrations of
the values of the translation functors $\theta_s^l$ and $\theta_t^l$ applies
to simple modules:
\begin{displaymath}
\begin{array}{|c||c|c|c|c|c|c|}
\hline
M               & e & s & t & st & ts & w_0 \\
\hline
\hline
\theta_s^l M    & 0 & \begin{array}{ccc}&s& \\
st& & e\\ &s&\end{array}& 0 & 0 &
\begin{array}{c}ts \\ t \\ ts\end{array}&
\begin{array}{c}w_0 \\ st \\ w_0\end{array} \\
\hline
\theta_t^l M    & 0 &0 & \begin{array}{ccc}&t& \\
ts& & e\\ &t&\end{array} & \begin{array}{c}st \\ s \\ st\end{array}&
0& \begin{array}{c}w_0 \\ ts \\ w_0\end{array}
\\ \hline
\end{array}
\end{displaymath}
From this we can draw the following graph which shows all the non-zero
coefficients of the action of $\theta_s^l$ (indicated by solid arrows) and
$\theta_t^l$ (indicated by dotted arrows) in the bases of simple modules:
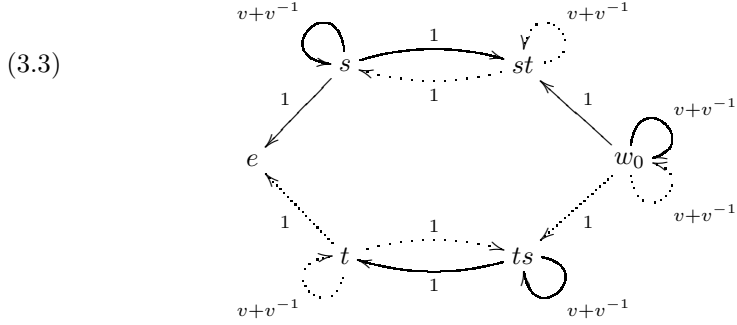
\hspace{2mm}

\begin{equation}\label{figuresl3}
\xymatrix{
  & s\ar[dl]_1\ar@/^/[rr]^1\ar@(u,l)[]_{v+v^{-1}} &&
  st\ar@{.>}@/^/[ll]^1\ar@(r,u)@{.>}[]_{v+v^{-1}} & \\
e & && & w_0\ar[ul]_1\ar@{.>}[dl]^1\ar@(u,r)[]^{v+v^{-1}}
\ar@(d,r)@{.>}[]_{v+v^{-1}}\\
  & t\ar@{.>}[ul]^1\ar@{.>}@/^/[rr]^1\ar@(d,l)@{.>}[]^{v+v^{-1}}
  && ts\ar@/^/[ll]^1\ar@(r,d)[]^{v+v^{-1}} & \\
}
\end{equation}
\hspace{2mm}

The graph \eqref{figuresl3} should be compared for
example with \cite[Figure~6.2]{BjBr} (in order to get
\cite[Figure~6.2]{BjBr} one should formally evaluate $v=1$ and
subtract the identity from $\theta_s^l$ and $\theta_t^l$).
From \eqref{figuresl3} one can deduce immediately the existence
of the following flag of $\mathds{H}$-submodules inside our regular
$\mathds{H}$-module:
\begin{displaymath}
\langle[e]\rangle\subset
\langle[e],[s],[st]\rangle\subset
\langle[e],[s],[st],[t],[ts]\rangle\subset
\langle[e],[s],[st],[t],[ts],[w_0]\rangle.
\end{displaymath}
The subquotients of this flag are the Kazhdan-Lusztig cell modules
for $\mathds{H}$. As we will show later on, this can be extended to
an explicit categorification of these cell modules via some
subcategories of $\mathcal{O}_0$. The definition of a left cell and
the categorification of cell modules is the topic of the next
section.

\section{Categorifications of Cell and Specht modules}\label{s3}

In this section we will introduce two categorifications of cell
modules - one which we believe is `the correct one' and one which
seems to be more canonical, easier, and straight forward on the
first sight, but turns out to be less natural at the end. We do not
know if the associated categories are in fact derived equivalent.

\subsection{Kazhdan-Lusztig's cell theory}\label{s3.1}

In this subsection we recall some facts from the Kazhdan-Lusztig cell theory.
Our main references here are \cite{KLCoxeter} and \cite{BjBr} and we refer
the reader to these papers for details. We will use the notation from
\cite{SoKipp}.

If $x\leq y$ then denote by $\mu(x,y)$ the coefficient of $v$ in the
Kazhdan-Lusztig polynomial $h_{x,y}$ and extend it to a symmetric function
$\mu:W\times W\to \mathbb{Z}$. In our normalization the formula
\cite[(1.0.a)]{KLCoxeter} reads then as follows:
\begin{equation}\label{formula1}
\underline{H}_x\underline{H}_s=
\begin{cases}
\underline{H}_{xs}+\sum_{y<x,ys<y}\mu(y,x)\underline{H}_y, & xs>x;\\
(v+v^{-1})\underline{H}_{x}, & xs<x.
\end{cases}
\end{equation}
In particular, $\mu(x,xs)=\mu(xs,x)=1$ for any $x\in W$ and $s\in S$.
For $w\in W$ define the {\em left} and the {\em right descent} sets of $w$
as follows:
\begin{displaymath}
D_\mathsf{L}(w):=\{s\in S\,:\, sw<w\},\quad
D_\mathsf{R}(w):=\{s\in S\,:\, ws<w\}.
\end{displaymath}
Now for $x,y\in W$ we write $x\rightarrow_\mathsf{L} y$ provided that
$\mu(x,y)\neq 0$ and there is some $s\in S$ such that $s\in D_{\mathsf{L}}(x)$
and $s\not\in D_\mathsf{L}(y)$. Denote by $\geq_{\mathsf{L}}$ the transitive
closure of the relation $\rightarrow_{\mathsf{L}}$. The relation
$\geq_\mathsf{L}$ is called the {\em left} pre-order on $W$. The equivalence
classes with respect to $\geq_\mathsf{L}$ are called the {\em left cells}.
The fact that $x,y\in W$ belong to the same left cell will be denoted
$x\sim_{\mathsf{L}}y$. The {\it right} versions $\geq_\mathsf{R}$ and
$\sim_{\mathsf{R}}$ of the above are obtained by applying the
involution $x\mapsto x^{-1}$, which  yields the notion of {\em right cells}.

Given  a right  cell $\mathbf{R}\subset W$, the $\mathbb{C}[v,v^{-1}]$-span
$X$ of $\underline{H}_x$, $x\geq_{\mathsf{R}} \mathbf{R}$, carries a natural
structure of a right $\mathds{H}$-module via \eqref{formula1}. The $\mathbb{C}[v,v^{-1}]$-span $Y$ of $\underline{H}_x$, $x>_\mathsf{R}
\mathbf{R}$, is a submodule of $X$. The $\mathds{H}$-module $X/Y$ is called the
{\em (right) cell module} associated with $\mathbf{R}$ and will be denoted by
$S(\mathbf{R})$. We leave it as an exercise to the reader to verify
that our definition of a cell module in fact agrees with the one from
\cite{KLCoxeter}.

\subsection{Presentable modules}\label{s3.105}

Here we would like to recall the construction of the category of presentable
modules from \cite{Au}, a basic construction which will be crucial in the
sequel. Let $\mathscr{A}$ be an abelian category and $\mathscr{B}$
be a full additive subcategory of $\mathscr{A}$. Denote by
$\overline{\mathscr{B}}$ the full subcategory of $\mathscr{A}$, which
consists of all $M\in \mathscr{A}$ for which there is an exact sequence
$N_1\to N_0\to M\to 0$ with $N_1,N_0\in \mathscr{B}$. This exact sequence
is called a $\mathscr{B}$-presentation of $M$. In the special case when
$\mathscr{B}=\mathrm{add}(P)$ for some projective object $P\in \mathscr{A}$
we have that $\overline{\mathscr{B}}$ is equivalent to the category of
right $\mathrm{End}_{\mathscr{A}}(P)$-modules, see
\cite[Section~5]{Au}. In particular, $\overline{\mathrm{add}(P)}$
is abelian.

\subsection{Categorification of cell modules}\label{s3.2}

Let $\mathbf{R}$ be a right cell of $W$. Set
\begin{displaymath}
\mathbf{\hat{R}}=\{w\in W\,:\,w\leq_\mathsf{R} x
\text{ for some }x\in \mathbf{R}\}.
\end{displaymath}
Let $\mathcal{O}_0^{\mathbf{\hat{R}}}$ denote the full subcategory
of $\mathcal{O}_0$, whose objects are all $M\in \mathcal{O}_0$ such that
each composition subquotient of $M$ has the form $L(w)$,
$w\in \mathbf{\hat{R}}$. For example if $\mathbf{R}=\{e\}$, the
category $\mathcal{O}_0^{\mathbf{\hat{R}}}$ contains only finite
direct sums of copies of the trivial $\mathfrak{g}$-module. In any case,
the inclusion functor $\mathfrak{i}^{\mathbf{\hat{R}}}:
\mathcal{O}_0^{\mathbf{\hat{R}}}\hookrightarrow \mathcal{O}_0$ is
exact and has as left adjoint the functor $\mathrm{Z}^{\mathbf{\hat{R}}}$
which picks out the maximal quotient contained in
$\mathcal{O}_0^{\mathbf{\hat{R}}}$. In particular, the indecomposable
projective modules in $\mathcal{O}_0^{\mathbf{\hat{R}}}$ are
the $P^{\mathbf{\hat{R}}}(w)=\mathrm{Z}^{\mathbf{\hat{R}}}P(w)$,
$w\in \mathbf{\hat{R}}$.

\begin{remark}\label{rem2}
{\rm
If $\mathbf{R}$ contains $w_0^{\mathfrak{p}}w_0$, where $w_0^{\mathfrak{p}}$
is the longest element in the parabolic subgroup of $W$ corresponding to a
parabolic subalgebra $\mathfrak{p}\supset\mathfrak{b}$ of $\mathfrak{g}$, then
$\mathcal{O}_0^{\mathbf{\hat{R}}}=\mathcal{O}_0^{\mathfrak{p}}$, the principal
block of the parabolic category $\mathcal{O}$ in the sense of \cite{RC}. This
follows from \cite[Proposition~6.2.7]{BjBr} and the fact that all simple
modules in $\mathcal{O}_0^{\mathfrak{p}}$ can be obtained as subquotients of
translations of the simple tilting module in $\mathcal{O}_0^{\mathfrak{p}}$
(as shown in \cite{CI}).
}
\end{remark}

\begin{proposition}\label{prop5}
Let $\mathbf{R}$ be a right cell of $W$.
\begin{enumerate}[(i)]
\item\label{prop5.1} The category $\mathcal{O}_0^{\mathbf{\hat{R}}}$ is
stable under $\theta_s^l$, $s\in S$.
\item\label{prop5.2} The additive category generated by
$P^{\mathbf{\hat{R}}}(w)$, $w\in \mathbf{R}$, is stable under $\theta_s^l$,
$s\in S$.
\end{enumerate}
\end{proposition}

\begin{proof}
To prove \eqref{prop5.1} it is enough to show that
$\theta_s^lL(w)\in \mathcal{O}_0^{\mathbf{\hat{R}}}$
for all $w\in \mathbf{\hat{R}}$. For $z\in W$ using
the self-adjointness of $\theta_s^l$, equation \eqref{formula1}
and Proposition~\ref{prop3}\eqref{prop3.3} we have:
\begin{eqnarray*}
\mathrm{Hom}_{\mathcal{O}}(P(z),\theta_s^l L(w))
&=&\mathrm{Hom}_{\mathcal{O}}(\theta_s^l P(z),L(w))\\
&=&
\begin{cases}
\mathrm{Hom}_{\mathcal{O}}(P(z)\oplus P(z),L(w)), & zs<z;\\
\mathrm{Hom}_{\mathcal{O}}(P(zs)\oplus
{\displaystyle\bigoplus_{y<z,ys<y}}P(y)^{\mu(y,z)},L(w)), & zs>z.
\end{cases}
\end{eqnarray*}
The latter space has the chance to be non-zero only in the following cases:
$z=w$ or $zs=w>z$, or, finally, $w<z$ where $ws<w$ and $\mu(w,z)\neq 0$.
In all these cases  $w\in \mathbf{\hat{R}}$ implies $z\in \mathbf{\hat{R}}$
and \eqref{prop5.1} follows.

To prove \eqref{prop5.2} we use \eqref{prop5.1} and note that $\theta_s^l$ maps
projectives from $\mathcal{O}_0^{\mathbf{\hat{R}}}$ to projectives from
$\mathcal{O}_0^{\mathbf{\hat{R}}}$ since it is self-adjoint. Now take $x\in
\mathbf{R}$. Then $\theta_s^lP^{\mathbf{\hat{R}}}(x)$ is a direct sum of some
$P^{\mathbf{\hat{R}}}(y)$'s. The possible $y$'s to occur are given by
\eqref{formula1}, hence either $y=x$ or $y\in \mathbf{R}$, or $y=xs>x$. In the
last case we have either $y\in \mathbf{R}$ or $y\not\in \mathbf{\hat{R}}$,
which is not possible since $\theta_s^l$ preserves
$\mathcal{O}_0^{\mathbf{\hat{R}}}$ by \eqref{prop5.1}. This completes the proof.
\end{proof}

We know already that the indecomposable projective module
$P(x)\in\mathcal{O}_0$ has a {\em standard} graded lift
$\mathtt{P}(x)$ for all $x\in W$ (for the definition of graded lift
we refer to \cite[Section 3]{Stgrad}; here and further a {\em
standard} graded lift of a projective or simple or standard module
is the lift in which the top of the module is concentrated in degree
zero). Now for $x\in \mathbf{\hat{R}}$ the module
$P^{\mathbf{\hat{R}}}(x)=\mathrm{Z}^{\mathbf{\hat{R}}}P(x)$ is the
quotient of $P(x)$ modulo the trace of all $P(y)$ such that $y\not
\leq_{\mathsf{R}}x$. The corresponding quotient
$\mathtt{P}^{\mathbf{\hat{R}}}(x)$ of $\mathtt{P}(x)$ is then a
standard graded lift of $P^{\mathbf{\hat{R}}}(x)$. Let
$\mathcal{P}^{\mathbf{R}}$ be the additive category, closed under
grading shifts, and generated by
$\mathtt{P}^{\mathbf{\hat{R}}}(w)$, $w\in \mathbf{R}$. This category
is the graded version of the additive category from
Proposition~\ref{prop5}\eqref{prop5.2}.  Set
$\mathscr{C}^{\mathbf{R}}=\overline{\mathcal{P}^{\mathbf{R}}}$ (see
Subsection~\ref{s3.105}), which is equivalent to the category of
graded finite-dimensional right modules over the algebra
$B^{\mathbf{R}}:= \mathrm{End}_{\mathcal{O}_0} (\oplus_{w\in
\mathbf{R}} P^{\mathbf{\hat{R}}}(w))$, which inherits a grading from
the algebra $A$ (Subsection~\ref{s25.5}). From
Proposition~\ref{prop5}\eqref{prop5.2} it follows that
$\mathscr{C}^{\mathbf{R}}$ is closed under $\theta_s^l$, $s\in S$.
Our first result is the following statement (we recall that
$\mathbb{Z}((v))$ denotes the ring of formal Laurent series in $v$
with integer coefficients):

\begin{theorem}[Categorification of cell modules]\label{thm5}
{\tiny .}

\begin{enumerate}[(i)]
\item\label{thm5.0}
There is a unique monomorphism of $\mathds{H}$-modules
such that
\begin{eqnarray*}
\mathcal{E}^{\mathbf{R}}:\quad\quad\quad S(\mathbf{R}) &
\longrightarrow & \left[ \mathscr{C}^{\mathbf{R}}
\right] \\ \underline{H}_w & \mapsto &
\left[\mathtt{P}^{\mathbf{\hat{R}}}(w)\right].
\end{eqnarray*}
\item\label{thm5.1}
The monomorphism  $\mathcal{E}^{\mathbf{R}}$ defines a precategorification
$(\mathscr{C}^{\mathbf{R}},\mathcal{E}^{\mathbf{R}},
\{\theta_s^l\}_{s\in S})$ and induces a categorification
$(\mathcal{P}^{\mathbf{R}},\mathcal{E}^{\mathbf{R}},
\{\theta_s^l\}_{s\in S})$ of the right cell $\mathds{H}$-module
$S(\mathbf{R})$ with respect to the generators $H_s$, $s\in S$.
\item\label{thm5.2}
The monomorphism $\mathcal{E}^{\mathbf{R}}$ from \eqref{thm5.0}
extends uniquely to a
$\mathbb{Z}((v))$-ca\-te\-go\-ri\-fi\-ca\-tion
$(\mathscr{C}^{\mathbf{R}},\mathcal{E}^{\mathbf{R}},\{\theta_s^l\}_{s\in
S})$ of the right  cell $\mathds{H}^{\mathbb{Z}((v))}$-module
$S(\mathbf{R})^{\mathbb{Z}((v))}$ with respect to the generators
$H_s$, $s\in S$.
\end{enumerate}
\end{theorem}

\begin{proof}
The statement \eqref{thm5.0} follows from
Proposition~\ref{prop3}\eqref{prop3.1}, Proposition~\ref{prop5}\eqref{prop5.2}
and the definitions. The statement \eqref{thm5.1} follows from \eqref{thm5.0}.
Note that $B^{\mathbf{R}}$ has infinite homological dimension in general.
Hence the statement \eqref{thm5.2} follows from \eqref{thm5.1} as the
extension of scalars from $\mathbb{Z}[v,v^{-1}]$ to $\mathbb{Z}((v))$ allows
one to work with infinite projective resolutions.
\end{proof}

\subsection{Remarks on another categorification of cell modules}
\label{s3.244}

Formula \eqref{formula1} suggests another way to categorify cell modules.
For a right cell $\mathbf{R}$ of $W$ set
\begin{displaymath}
\mathbf{\check{R}}=\{w\in W\,:\,x\leq_\mathsf{R} w
\text{ for some }x\in \mathbf{R}\}
\end{displaymath}
(note the difference to $\mathbf{\hat{R}}$). Let $\mathscr{A}$  denote the additive category, generated by $P(w)$, $w\in \mathbf{\check{R}}$. Denote
also by $\mathscr{A}'$  the additive category, generated by $P(w)$, $w\in \mathbf{\check{R}}\setminus \mathbf{R}$. Consider the categories
$\mathcal{O}_0^{\mathbf{\check{R}}}=\overline{\mathscr{A}}$
and $\tilde{\mathcal{O}}_0^{\mathbf{\check{R}}}=\overline{\mathscr{A}'}$.

Note that if $\mathbf{R}$ contains $w_0^{\mathfrak{p}}$, where
$w_0^{\mathfrak{p}}$ is the longest element in the parabolic subgroup
$W_{\mathfrak{p}}$ of $W$ corresponding to a parabolic subalgebra
$\mathfrak{p}\supset\mathfrak{b}$ of $\mathfrak{g}$, then
$\mathcal{O}_0^{\mathbf{\check{R}}}$ coincides with the category of
$\mathfrak{p}$-presentable modules in $\mathcal{O}_0$
(\cite[Section~2]{MS}) and is equivalent to
${}_{\mathbf{0}}\mathcal{H}^1_{\lambda}$, where $\lambda\in\h^*_{dom}$ is
integral and has stabilizer $W_{\mathfrak{p}}$ (\cite[Theorem~5.9(ii)]{BG}).

Formula \eqref{formula1} and Proposition~\ref{prop3}\eqref{prop3.3}
immediately imply that both, the category $\mathcal{O}_0^{\mathbf{\check{R}}}$
and the category $\tilde{\mathcal{O}}_0^{\mathbf{\check{R}}}$, are stable under
$\theta_s^l$, $s\in S$. And the `quotient' should be exactly the cell module.
To define this `quotient' we let $\mathcal{Q}^{\mathbf{R}}$ denote the additive
category, closed under grading shifts, and generated by $\mathtt{P}(w)$, $w\in
\mathbf{R}$. Set $\mathscr{D}^{\mathbf{R}}=
\overline{\mathcal{Q}^{\mathbf{R}}}$. The functors $\theta_s^l$, $s\in S$, do
not preserve $\mathcal{Q}^{\mathbf{R}}$ unless $\mathbf{R}=\{w_0\}$. However,
one can use them to define right exact functors $\tilde{\theta}_s^l$ on
$\mathscr{D}^{\mathbf{R}}$ as follows: First we define the functor
$\tilde{\theta}_s^l$ on the indecomposable projective module
$\mathtt{P}(x)$. Let $s\in S$ and $x\in \mathbf{R}$. If
$\theta_s^l\mathtt{P}(x)\in \mathcal{Q}^{\mathbf{R}}$, we set
$\tilde{\theta}_s^l\mathtt{P}(x)=\theta_s^l\mathtt{P}(x)$, otherwise
\eqref{formula1} gives
\begin{displaymath}
\theta_s^l\mathtt{P}(x)=\mathtt{P}(xs)\oplus
\bigoplus_{y<x,ys<y}\mathtt{P}(y)^{\mu(y,x)}.
\end{displaymath}
This decomposition into two summands is unique since the first summand
coincides with the trace of the module $\mathtt{P}(xs)$ in
$\theta_s^l\mathtt{P}(x)$ and the second summand coincides with the trace of
the module $\oplus_{w\in \mathbf{R}}\mathtt{P}(w)$ in $\theta_s^l\mathtt{P}(x)$.
Hence we can define $\tilde{\theta}_s^l\mathtt{P}(x)=
\bigoplus_{y<x,ys<y}\mathtt{P}(y)^{\mu(y,x)}$ and define
$\tilde{\theta}_s^l$ on morphisms via restriction. In the standard way
$\tilde{\theta}_s^l$ extends uniquely to a right exact endofunctor on
$\mathscr{D}^{\mathbf{R}}$. We do not know if $\tilde{\theta}_s^l$ is exact.
By \eqref{formula1}, the action of  $\tilde{\theta}_s^l$, $s\in S$, on the
Grothendieck group of $\mathcal{D}^b(\mathscr{D}^{\mathbf{R}})$ coincides with
the action of $\underline{H}_s$ on $S(\mathbf{R})$ and hence we obtain another
categorification of the cell module $S(\mathbf{R})$. We do not know whether
this categorification is (derived) equivalent to the one constructed in
Theorem~\ref{thm5} or not. The principal disadvantage with this
categorification is that we do not know to which extend our uniqueness result
from Subsection~\ref{s4.1} holds in this setup.

\subsection{$\mathfrak{gl}_3$-example}\label{s3.245}

\begin{figure}[tbh]
\special{em:linewidth 0.4pt} \unitlength 0.80mm
\linethickness{1pt}
\begin{picture}(150.00,40.00)
\special{em:linewidth 0.4pt} \unitlength 0.80mm
\linethickness{1pt}
\drawline(00.00,00.00)(150.00,00.00)
\drawline(00.00,00.00)(00.00,40.00)
\drawline(150.00,40.00)(00.00,40.00)
\drawline(150.00,40.00)(150.00,00.00)
\drawline(150.00,30.00)(00.00,30.00)
\drawline(15.00,40.00)(15.00,00.00)
\drawline(35.00,40.00)(35.00,00.00)
\drawline(55.00,40.00)(55.00,00.00)
\drawline(85.00,40.00)(85.00,00.00)
\drawline(115.00,40.00)(115.00,00.00)
\put(7.50,35.00){\makebox(0,0)[cc]{\tiny{$P(e)$}}}
\put(7.50,25.00){\makebox(0,0)[cc]{\tiny{$e$}}}
\put(5.00,22.00){\makebox(0,0)[cc]{\tiny{$s$}}}
\put(10.00,22.00){\makebox(0,0)[cc]{\tiny{$t$}}}
\put(5.00,19.00){\makebox(0,0)[cc]{\tiny{$st$}}}
\put(10.00,19.00){\makebox(0,0)[cc]{\tiny{$ts$}}}
\put(7.50,16.00){\makebox(0,0)[cc]{\tiny{$w_0$}}}
\linethickness{0.2pt}
\dashline{1}(7.50,28.00)(2.50,23.00)
\dashline{1}(7.50,28.00)(12.50,23.00)
\dashline{1}(2.50,23.00)(2.50,18.00)
\dashline{1}(12.50,18.00)(12.50,23.00)
\dashline{1}(12.50,18.00)(7.50,13.00)
\dashline{1}(7.50,13.00)(2.50,18.00)
\put(25.00,35.00){\makebox(0,0)[cc]{\tiny{$P(s)$}}}
\put(29.50,22.00){\makebox(0,0)[cc]{\tiny{$e$}}}
\put(27.00,19.00){\makebox(0,0)[cc]{\tiny{$s$}}}
\put(32.00,19.00){\makebox(0,0)[cc]{\tiny{$t$}}}
\put(27.00,16.00){\makebox(0,0)[cc]{\tiny{$st$}}}
\put(32.00,16.00){\makebox(0,0)[cc]{\tiny{$ts$}}}
\put(29.50,13.00){\makebox(0,0)[cc]{\tiny{$w_0$}}}
\dashline{1}(29.50,23.50)(25.50,19.00)
\dashline{1}(29.50,23.50)(33.50,19.00)
\dashline{1}(25.50,19.00)(25.50,14.00)
\dashline{1}(33.50,14.00)(33.50,19.00)
\dashline{1}(33.50,14.00)(29.50,11.00)
\dashline{1}(29.50,11.00)(25.50,14.00)
\put(21.50,25.00){\makebox(0,0)[cc]{\tiny{$s$}}}
\put(19.50,22.00){\makebox(0,0)[cc]{\tiny{$st$}}}
\put(24.50,22.00){\makebox(0,0)[cc]{\tiny{$ts$}}}
\put(21.50,19.00){\makebox(0,0)[cc]{\tiny{$w_0$}}}
\dashline{1}(21.50,27.00)(17.00,22.00)
\dashline{1}(21.50,27.00)(26.50,22.00)
\dashline{1}(21.50,16.50)(17.00,22.00)
\dashline{1}(21.50,16.50)(26.50,22.00)
\put(45.00,35.00){\makebox(0,0)[cc]{\tiny{$P(t)$}}}
\put(49.50,22.00){\makebox(0,0)[cc]{\tiny{$e$}}}
\put(47.00,19.00){\makebox(0,0)[cc]{\tiny{$s$}}}
\put(52.00,19.00){\makebox(0,0)[cc]{\tiny{$t$}}}
\put(47.00,16.00){\makebox(0,0)[cc]{\tiny{$st$}}}
\put(52.00,16.00){\makebox(0,0)[cc]{\tiny{$ts$}}}
\put(49.50,13.00){\makebox(0,0)[cc]{\tiny{$w_0$}}}
\dashline{1}(49.50,23.50)(45.50,19.00)
\dashline{1}(49.50,23.50)(53.50,19.00)
\dashline{1}(45.50,19.00)(45.50,14.00)
\dashline{1}(53.50,14.00)(53.50,19.00)
\dashline{1}(53.50,14.00)(49.50,11.00)
\dashline{1}(49.50,11.00)(45.50,14.00)
\put(41.50,25.00){\makebox(0,0)[cc]{\tiny{$t$}}}
\put(39.50,22.00){\makebox(0,0)[cc]{\tiny{$st$}}}
\put(44.50,22.00){\makebox(0,0)[cc]{\tiny{$ts$}}}
\put(41.50,19.00){\makebox(0,0)[cc]{\tiny{$w_0$}}}
\dashline{1}(41.50,27.00)(37.00,22.00)
\dashline{1}(41.50,27.00)(46.50,22.00)
\dashline{1}(41.50,16.50)(37.00,22.00)
\dashline{1}(41.50,16.50)(46.50,22.00)
\put(70.00,35.00){\makebox(0,0)[cc]{\tiny{$P(st)$}}}
\put(70.00,19.00){\makebox(0,0)[cc]{\tiny{$e$}}}
\put(68.00,16.00){\makebox(0,0)[cc]{\tiny{$s$}}}
\put(72.50,16.00){\makebox(0,0)[cc]{\tiny{$t$}}}
\put(68.00,13.00){\makebox(0,0)[cc]{\tiny{$st$}}}
\put(72.50,13.00){\makebox(0,0)[cc]{\tiny{$ts$}}}
\put(70.00,10.00){\makebox(0,0)[cc]{\tiny{$w_0$}}}
\dashline{1}(70.00,20.50)(66.00,16.00)
\dashline{1}(70.00,20.50)(74.00,16.00)
\dashline{1}(66.00,16.00)(66.00,11.00)
\dashline{1}(74.00,11.00)(74.00,16.00)
\dashline{1}(74.00,11.00)(70.00,8.00)
\dashline{1}(70.00,8.00)(66.00,11.00)
\put(62.00,22.00){\makebox(0,0)[cc]{\tiny{$s$}}}
\put(60.00,19.00){\makebox(0,0)[cc]{\tiny{$st$}}}
\put(65.00,19.00){\makebox(0,0)[cc]{\tiny{$ts$}}}
\put(62.00,16.00){\makebox(0,0)[cc]{\tiny{$w_0$}}}
\dashline{1}(62.00,24.00)(57.50,19.00)
\dashline{1}(62.00,24.00)(67.00,19.00)
\dashline{1}(62.00,13.50)(57.50,19.00)
\dashline{1}(62.00,13.50)(67.00,19.00)
\put(78.00,22.00){\makebox(0,0)[cc]{\tiny{$t$}}}
\put(76.00,19.00){\makebox(0,0)[cc]{\tiny{$st$}}}
\put(81.00,19.00){\makebox(0,0)[cc]{\tiny{$ts$}}}
\put(78.00,16.00){\makebox(0,0)[cc]{\tiny{$w_0$}}}
\dashline{1}(78.00,24.00)(73.50,19.00)
\dashline{1}(78.00,24.00)(83.00,19.00)
\dashline{1}(78.00,13.50)(73.50,19.00)
\dashline{1}(78.00,13.50)(83.00,19.00)
\put(70.00,25.00){\makebox(0,0)[cc]{\tiny{$st$}}}
\put(70.00,22.00){\makebox(0,0)[cc]{\tiny{$w_0$}}}
\dashline{1}(68.00,27.00)(72.00,27.00)
\dashline{1}(68.00,27.00)(68.00,20.00)
\dashline{1}(68.00,20.00)(72.00,20.00)
\dashline{1}(72.00,27.00)(72.00,20.00)
\put(100.00,35.00){\makebox(0,0)[cc]{\tiny{$P(ts)$}}}
\put(100.00,19.00){\makebox(0,0)[cc]{\tiny{$e$}}}
\put(98.00,16.00){\makebox(0,0)[cc]{\tiny{$s$}}}
\put(102.50,16.00){\makebox(0,0)[cc]{\tiny{$t$}}}
\put(98.00,13.00){\makebox(0,0)[cc]{\tiny{$st$}}}
\put(102.50,13.00){\makebox(0,0)[cc]{\tiny{$ts$}}}
\put(100.00,10.00){\makebox(0,0)[cc]{\tiny{$w_0$}}}
\dashline{1}(100.00,20.50)(96.00,16.00)
\dashline{1}(100.00,20.50)(104.00,16.00)
\dashline{1}(96.00,16.00)(96.00,11.00)
\dashline{1}(104.00,11.00)(104.00,16.00)
\dashline{1}(104.00,11.00)(100.00,8.00)
\dashline{1}(100.00,8.00)(96.00,11.00)
\put(92.00,22.00){\makebox(0,0)[cc]{\tiny{$s$}}}
\put(90.00,19.00){\makebox(0,0)[cc]{\tiny{$st$}}}
\put(95.00,19.00){\makebox(0,0)[cc]{\tiny{$ts$}}}
\put(92.00,16.00){\makebox(0,0)[cc]{\tiny{$w_0$}}}
\dashline{1}(92.00,24.00)(87.50,19.00)
\dashline{1}(92.00,24.00)(97.00,19.00)
\dashline{1}(92.00,13.50)(87.50,19.00)
\dashline{1}(92.00,13.50)(97.00,19.00)
\put(108.00,22.00){\makebox(0,0)[cc]{\tiny{$t$}}}
\put(106.00,19.00){\makebox(0,0)[cc]{\tiny{$st$}}}
\put(111.00,19.00){\makebox(0,0)[cc]{\tiny{$ts$}}}
\put(108.00,16.00){\makebox(0,0)[cc]{\tiny{$w_0$}}}
\dashline{1}(108.00,24.00)(103.50,19.00)
\dashline{1}(108.00,24.00)(113.00,19.00)
\dashline{1}(108.00,13.50)(103.50,19.00)
\dashline{1}(108.00,13.50)(113.00,19.00)
\put(100.00,25.00){\makebox(0,0)[cc]{\tiny{$ts$}}}
\put(100.00,22.00){\makebox(0,0)[cc]{\tiny{$w_0$}}}
\dashline{1}(98.00,27.00)(102.00,27.00)
\dashline{1}(98.00,27.00)(98.00,20.00)
\dashline{1}(98.00,20.00)(102.00,20.00)
\dashline{1}(102.00,27.00)(102.00,20.00)
\put(132.50,35.00){\makebox(0,0)[cc]{\tiny{$P(w_0)$}}}
\put(132.50,16.00){\makebox(0,0)[cc]{\tiny{$e$}}}
\put(130.50,13.00){\makebox(0,0)[cc]{\tiny{$s$}}}
\put(135.00,13.00){\makebox(0,0)[cc]{\tiny{$t$}}}
\put(130.50,10.00){\makebox(0,0)[cc]{\tiny{$st$}}}
\put(135.00,10.00){\makebox(0,0)[cc]{\tiny{$ts$}}}
\put(132.50,7.00){\makebox(0,0)[cc]{\tiny{$w_0$}}}
\dashline{1}(132.50,17.50)(128.50,13.00)
\dashline{1}(132.50,17.50)(136.50,13.00)
\dashline{1}(128.50,13.00)(128.50,8.00)
\dashline{1}(136.50,8.00)(136.50,13.00)
\dashline{1}(136.50,8.00)(132.50,5.00)
\dashline{1}(132.50,5.00)(128.50,8.00)
\put(122.50,19.00){\makebox(0,0)[cc]{\tiny{$s$}}}
\put(120.50,16.00){\makebox(0,0)[cc]{\tiny{$st$}}}
\put(125.50,16.00){\makebox(0,0)[cc]{\tiny{$ts$}}}
\put(122.50,13.00){\makebox(0,0)[cc]{\tiny{$w_0$}}}
\dashline{1}(122.50,21.00)(118.00,16.00)
\dashline{1}(122.50,21.00)(127.50,16.00)
\dashline{1}(122.50,10.50)(118.00,16.00)
\dashline{1}(122.50,10.50)(127.50,16.00)
\put(142.50,19.00){\makebox(0,0)[cc]{\tiny{$t$}}}
\put(140.50,16.00){\makebox(0,0)[cc]{\tiny{$st$}}}
\put(145.50,16.00){\makebox(0,0)[cc]{\tiny{$ts$}}}
\put(142.50,13.00){\makebox(0,0)[cc]{\tiny{$w_0$}}}
\dashline{1}(142.50,21.00)(138.00,16.00)
\dashline{1}(142.50,21.00)(147.50,16.00)
\dashline{1}(142.50,10.50)(138.00,16.00)
\dashline{1}(142.50,10.50)(147.50,16.00)
\put(128.50,22.00){\makebox(0,0)[cc]{\tiny{$st$}}}
\put(128.50,19.00){\makebox(0,0)[cc]{\tiny{$w_0$}}}
\dashline{1}(126.50,23.00)(130.50,23.00)
\dashline{1}(126.50,23.00)(126.50,18.00)
\dashline{1}(126.50,18.00)(130.50,18.00)
\dashline{1}(130.50,23.00)(130.50,18.00)
\put(136.50,22.00){\makebox(0,0)[cc]{\tiny{$ts$}}}
\put(136.50,19.00){\makebox(0,0)[cc]{\tiny{$w_0$}}}
\dashline{1}(134.50,23.00)(138.50,23.00)
\dashline{1}(134.50,23.00)(134.50,18.00)
\dashline{1}(134.50,18.00)(138.50,18.00)
\dashline{1}(138.50,23.00)(138.50,18.00)
\put(132.50,25.00){\makebox(0,0)[cc]{\tiny{$w_0$}}}
\dashline{1}(132.50,28.00)(129.50,25.00)
\dashline{1}(129.50,25.00)(132.50,22.00)
\dashline{1}(132.50,22.00)(135.50,25.00)
\dashline{1}(135.50,25.00)(132.50,28.00)
\end{picture}
\caption{Indecomposable projectives in $\mathcal{O}_0$}\label{figwam}
\end{figure}
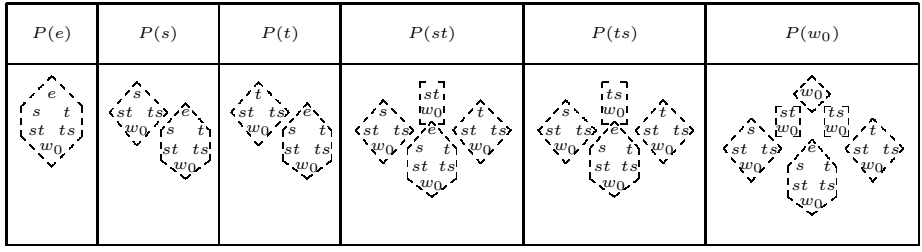

Let $W=\langle s,t\rangle\cong S_3$. Then there are four right
cells and the Hasse diagram of the right order is as follows:
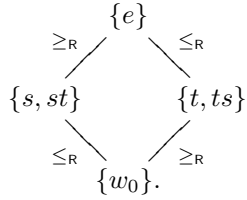
\begin{displaymath}
\xymatrix@!=0.6pc{
 &\{e\}\ar@{-}[dl]_{\geq_\mathsf{R}}\ar@{-}[dr]^{\leq_\mathsf{R}} &\\
\{s,st\}\ar@{-}[dr]_{\leq_\mathsf{R}}&&\{t,ts\}\ar@{-}[dl]^{\geq_\mathsf{R}}\\ &\{w_0\}.&
}
\end{displaymath}
Consider first the case
$\mathbf{R}=\{w_0\}$, where we have $\mathcal{O}_0^{\widehat{\{w_0\}}}=
\mathcal{O}_0$. It contains all simple modules $L(w)$, $w\in S_3$.
The presentation of this category as a module category over a finite dimensional algebra can be found in \cite[5.1.2]{St3}. The
graded filtrations of the indecomposable projective modules
(with indicated Verma subquotients) in this case are shown on
Figure~\ref{figwam}.
The category $\mathcal{P}^{\{w_0\}}$ contains (up to grading shift) a unique indecomposable
module, namely $P^{\{w_0\}}(w_0)$. The algebra $B^{\{w_0\}}=\End_\mg(P^{\{w_0\}}(w_0))$ is the
coinvariant algebra of $W$, see \cite[Endomorphismensatz]{Sperv}.

Below we collect the analogous information for the three other choices for the right cells,
in particular, we present all the algebras which appear there in
terms of quivers and relations.
{\tiny
\begin{displaymath}
\begin{array}{|c||c|c|c|}
\hline
\mathbf{R}   & \{e\} & \{s,st\} & \{t,ts\}  \\
\hline
\text{Simple modules:}& e & e, s, st & e, t, ts\\
\hline
\text{Projective modules:}& e &
\begin{array}{c|c|c}
P(e) & P(s) & P(st)\\
\hline
\begin{array}{c}e\\s\\\text{\hspace{2mm}} \end{array}&
\begin{array}{ccc}&s&\\st&&e\\&s&\end{array}&
\begin{array}{c}st\\s\\st\end{array}
\end{array}&
\begin{array}{c|c|c}
P(e) & P(t) & P(ts)\\
\hline
\begin{array}{c}e\\t\\\text{\hspace{2mm}} \end{array}&
\begin{array}{ccc}&t&\\ts&&e\\&t&\end{array}&
\begin{array}{c}ts\\t\\ts\end{array}
\end{array}
\\
\hline
\text{Quiver of $\mathcal{O}_0^{\mathbf{\hat{R}}}$:}&
e &
\begin{array}{c}\xymatrix{st\ar@/^/[r]^{\alpha}&
s\ar@/^/[r]^{\gamma}\ar@/^/[l]^{\beta}&e\ar@/^/[l]^{\delta}}
\\
\beta\delta=\gamma\alpha=\gamma\delta=0\\
\alpha\beta=\delta\gamma\end{array}
&\begin{array}{c}\xymatrix{ts\ar@/^/[r]^{\alpha}&
t\ar@/^/[r]^{\gamma}\ar@/^/[l]^{\beta}&e\ar@/^/[l]^{\delta}}
\\
\beta\delta=\gamma\alpha=\gamma\delta=0\\
\alpha\beta=\delta\gamma\end{array}
\\
\hline
\text{Quiver of $\mathscr{C}^{\mathbf{R}}$:}&
e &
\xymatrix{st\ar@/^/[r]^{\alpha}&s\ar@/^/[l]^{\beta}}
\begin{array}{l}\alpha\beta\alpha=\\\beta\alpha\beta=0\end{array}&
\xymatrix{ts\ar@/^/[r]^{\alpha}&t\ar@/^/[l]^{\beta}}
\begin{array}{l}\alpha\beta\alpha=\\\beta\alpha\beta=0\end{array}\\
\hline
\end{array}
\end{displaymath}
} In the above example the category $\mathcal{O}^\mathbf{\hat{R}}_0$
always coincides with some parabolic category
$\mathcal{O}^{\mathfrak{p}}_0$. This is not the case in general. The
smallest such example is the right cell $\{s_1s_3,s_1s_3s_2\}$ of
$S_4$.

\subsection{Specht modules}\label{s3.3}

In the special case $W=S_n$ we denote $\mathds{H}=\mathds{H}_n$. The (right) cell modules are exactly the irreducible $\mathds{H}_n$-modules,
\cite[Theorem~1.4]{KLCoxeter}. However, cell modules for different right cells (namely if they are in the same double cell) might be isomorphic. Theorem~\ref{thm5} gives therefore (several) categorifications for each
irreducible $\mathds{H}$-module. If we specialize $v=1$ (i.e. we forget the grading) and work over a field of characteristic zero, the irreducible modules
for the Hecke algebra specialize to irreducible modules for the symmetric group
(for an explicit description see for example \cite{Na}), hence we get
categorifications of Specht modules. In the special situation of
Remark~\ref{rem2} we obtain the categorification of Specht modules constructed
in \cite{KMS}.

Every cell module has a symmetric, non-degenerate, $\mathds{H}_n$-invariant
bilinear form $\langle\cdot,\cdot\rangle$ with values in $\mathbb{Z}[v,v^{-1}]$,
which is unique up to a scalar, see \cite[page~114]{Murphy}. There is a
categorical interpretation of this form as follows: For any $\Z$-graded complex
vector space $M=\oplus_{j\in\Z} M^j$ let $h(M)=\sum_{j\in\Z}(
\operatorname{dim}_\mC M^j)v^{j}\in\Z[v,v^{-1}]$ be the corresponding
Hilbert polynomial. For all $M$, $N\in\mathscr{C}^{\mathbf{R}}$ and all
$i\in\mathbb{Z}$ the vector space $E^i(M,N):=
\operatorname{Ext}^i_{\mathscr{C}^{\mathbf{R}}}(M,N)$ is  $\Z$-graded in the
natural way.  Set $h(E(M,N))=\sum_{i\in\Z}(-1)^i
h(E^i(M,N))$. Let $\operatorname{d}$ denote the graded lift of the standard
duality on $\cO_0$, restricted to the category  $\mathscr{C}^{\mathbf{R}}$.

\begin{proposition} \label{bilinearform}
The form
\begin{displaymath}
\beta(\cdot,\cdot):=h(E(\cdot,\operatorname{d}(\cdot))):
\mathscr{C}^{\mathbf{R}}\times \mathscr{C}^{\mathbf{R}} \rightarrow\Z((v))
\end{displaymath}
descends to a symmetric, non-degenerate, $\mathcal{H}_n^{\Z((v))}$-invariant
bilinear form $\langle\cdot,\cdot\rangle$ on the
$\mathcal{H}_n^{\Z((v))}$-module $[\mathscr{C}^{\mathbf{R}}]^{\Z((v))}$.
The restriction of this form to $\left[\mathcal{P}^{\mathbf{R}}\right]_{\oplus}$
has values in $\mathbb{Z}[v,v^{-1}]$.
\end{proposition}

\begin{proof}
The same as the proof of \cite[Proposition~4]{KMS}.
\end{proof}

\section{Uniqueness of the categorification for type $A$}\label{s4}

In this section we stick to the case where $W=S_n$. In the previous section we
constructed various categorifications for each single Specht module via cell
modules. In this section we will show that all these categorifications are in
fact equivalent. In particular, one can consider the categorification from
\cite{KMS} as a kind of `universal one'.

\subsection{Equivalence of categories}\label{s4.1}

\begin{theorem}[Uniqueness Theorem]\label{thm6}
Let $\mathbf{R}_1$ and $\mathbf{R}_2$ be two right cells of $W=S_n$, which
belong to the same double cell. Then there is an equivalence of categories
\begin{displaymath}
\Phi=\Phi_{\mathbf{R}_1}^{\mathbf{R}_2}:\quad
\mathscr{C}^{\mathbf{R}_1}\overset{\sim}{\rightarrow}
\mathscr{C}^{\mathbf{R}_2},
\end{displaymath}
which (naturally) commutes with projective functors and induces an isomorphism
of $\mathds{H}$-modules $[\mathscr{C}^{\mathbf{R}_1}]\cong
[\mathscr{C}^{\mathbf{R}_2}]$.
\end{theorem}

We will only prove the ungraded version of this theorem. The graded
version follows by standard arguments. For our proof we will need
several new definitions and more notation. For any right cell
$\mathbf{R}$ let $\mathscr{P}(\mathbf{R})$ denote the full additive
subcategories of $\mathcal{O}$, generated by all indecomposable
direct summands of the modules $E\otimes P^{\mathbf{\hat{R}}}(w)$,
$w\in \mathbf{R}$, where $E$ runs through all finite-dimensional
$\mathfrak{g}$-modules. Analogously we define
$\mathscr{P}(\mathbf{\hat{R}})$ using the condition $w\in
\mathbf{\hat{R}}$. Set $\mathcal{O}^{\mathbf{R}}=
\overline{\mathscr{P}(\mathbf{R})}$ and
$\mathcal{O}^{\mathbf{\hat{R}}}=
\overline{\mathscr{P}(\mathbf{\hat{R}})}$.

Denote by $\mathcal{O}_{\op{int}}$ the full subcategory of $\mathcal{O}$, which
consists of all modules with integral support (i.e. those modules $M$ such
that each weight of $M$ is also a weight of some finite-dimensional module).
Further, for $s\in S$ we denote
by $\mathcal{O}^s_{\op{int}}$ the integral part of the $s$-parabolic category,
that is the full subcategory of $\mathcal{O}_{\op{int}}$, which consists of
all modules which have only composition factors of the form  $L(w\cdot\lambda)$,
where $\lambda$ is an integral weight in $\mathfrak{h}_{dom}^*$,  $sw\cdot
\lambda\neq w\cdot \lambda$, and $sw>w$. For these categories we have the
natural inclusion $\mathrm{i}_s:\mathcal{O}^s_{\op{int}}\hookrightarrow
\mathcal{O}_{\op{int}}$ and we denote by $\mathrm{Z}_s$ and
$\hat{\mathrm{Z}}_s$ the left and the right adjoint to this inclusion
respectively. These are the classical {\em Zuckerman functors}.

If  $\mathbf{R}$ is a right cell such that $\mathbf{R}\leq_\mathsf{R} sw_0$,
then we have the natural inclusion $\mathrm{i}_s^{\mathbf{\hat{R}}}:
\mathcal{O}^{\mathbf{\hat{R}}}\hookrightarrow \mathcal{O}^s_{\op{int}}$
and we denote by $\mathrm{Z}_s^{\mathbf{\hat{R}}}$ and
$\hat{\mathrm{Z}}_s^{\mathbf{\hat{R}}}$ respectively the left and the right
adjoint to this inclusion.

Let now $\mathbf{R}_1$ and $\mathbf{R}_2$ be two right cells. Assume
that (see \cite[Proof of Theorem~1.4]{KLCoxeter})
\begin{equation}\label{Lcondition}
\exists\, s,t\in S\text{ and } w\in \mathbf{R}_1\text{ such that }
(st)^3=e, sw\geq w, tw\leq w, tw\in \mathbf{R}_2.
\end{equation}
In this case we have the following picture:
\begin{displaymath}
\xymatrix{
\mathcal{D}^b(\mathcal{O}^s_{\op{int}})
\ar@/^/[rr]^{\mathrm{i}_s}
&&\mathcal{D}^b(\mathcal{O}_{\op{int}})
\ar@/^/[rr]^{\mathcal{L}\mathrm{Z}_t\llbracket-1\rrbracket}
\ar@/^/[ll]^{\mathcal{R}\hat{\mathrm{Z}}_s}
&&
\mathcal{D}^b(\mathcal{O}^t_{\op{int}})
\ar@/^/[ll]^{\mathrm{i}_{t}\llbracket1\rrbracket}
}
\end{displaymath}
For this diagram we denote by $\mathrm{F}$ the composition from the left
to the right and by $\mathrm{G}$ the composition from the right to the left.
Directly from the definitions we have that $(F,G)$ is an adjoint pair of
functors. Furthermore, there are adjoint pairs
$(\mathrm{i}_s^{\mathbf{\hat{R}}_1},\hat{\mathrm{Z}}_s^{\mathbf{\hat{R}}_1})$
and $(\mathrm{Z}_t^{\mathbf{\hat{R}}_2},\mathrm{i}_t^{\mathbf{\hat{R}}_2})$
as follows:
\begin{displaymath}
\xymatrix{
\mathcal{O}^{\mathbf{\hat{R}}_1}\ar@/^/[rr]^{\mathrm{i}_s^{\mathbf{\hat{R}}_1}}
&&
\mathcal{O}^s_{\op{int}}\ar@/^/[ll]^{\hat{\mathrm{Z}}_s^{\mathbf{\hat{R}}_1}}
}\quad\quad\quad\quad
\xymatrix{
\mathcal{O}^t_{\op{int}}\ar@/^/[rr]^{\mathrm{Z}_t^{\mathbf{\hat{R}}_2}}
&&
\mathcal{O}^{\mathbf{\hat{R}}_2}\ar@/^/[ll]^{\mathrm{i}_t^{\mathbf{\hat{R}}_2}}
}
\end{displaymath}

\begin{lemma}\label{l21}
The functors $F$, $G$, $\mathrm{i}_s^{\mathbf{\hat{R}}_1}$,
$\hat{\mathrm{Z}}_s^{\mathbf{\hat{R}}_1}$, $\mathrm{i}_t^{\mathbf{\hat{R}}_2}$
and $\mathrm{Z}_t^{\mathbf{\hat{R}}_2}$ commute with functors of tensoring with
finite-\-di\-men\-si\-o\-nal $\mathfrak{g}$-modules, in particular with
projective functors.
\end{lemma}

\begin{proof}
Since all involved categories are stable under tensoring with
finite-di\-men\-si\-o\-nal $\mathfrak{g}$-modules by definition, all involved
inclusions commute with these functors. We will show how one derives from here
that $\mathrm{Z}_t^{\mathbf{\hat{R}}_2}$ commutes with tensoring with
finite-dimensional $\mathfrak{g}$-modules. For all other functors the arguments
are similar and therefore omitted.

Let $E$ be a finite-dimensional $\mathfrak{g}$-module. For each $M\in
\mathcal{O}^t_{\op{int}}$ from the definition of
$\mathrm{Z}_t^{\mathbf{\hat{R}}_2}$  we have the canonical projection $M\tto
\mathrm{Z}_t^{\mathbf{\hat{R}}_2}M$ with kernel $K_M$. Denote
$\theta:=E\otimes{}_-$, and let $P$ be a projective module in
$\mathcal{O}^t_{\op{int}}$ and $f\in \mathrm{End}_{\mathfrak{g}}(P)$.
Consider the following diagram:
\begin{equation}\label{eqses1}
\xymatrix{
&\theta K_{P}\ar@{^{(}->}[rr]\ar[ld]_{\overline{\theta f}}
\ar[dd]_>>>>>>>{\varphi}&&
\theta P\ar@{->>}[rr]\ar@{=}[dd]\ar[ld]_{\theta f}&&
\theta \mathrm{Z}_t^{\mathbf{\hat{R}}_2}P
\ar[ld]^{\theta\mathrm{Z}_t^{\mathbf{\hat{R}}_2} f}
\ar[dd]_>>>>>>>{\varphi'}\\
\theta K_{P}\ar@{^{(}->}[rr]\ar[dd]_{\varphi} &&
\theta P\ar@{->>}[rr]\ar@{=}[dd]&&\theta\mathrm{Z}_t^{\mathbf{\hat{R}}_2}P
\ar[dd]_>>>>>>>{\varphi'}& \\
&K_{\theta P}\ar@{^{(}->}[rr]\ar[ld]_{\overline{\theta f}}&&
\theta P\ar@{->>}[rr]\ar[ld]_{\theta f}&&
\mathrm{Z}_t^{\mathbf{\hat{R}}_2}\theta P
\ar[ld]^{\mathrm{Z}_t^{\mathbf{\hat{R}}_2} \theta f}
\\
K_{\theta P}\ar@{^{(}->}[rr] &&
\theta P\ar@{->>}[rr]&&\mathrm{Z}_t^{\mathbf{\hat{R}}_2}\theta P&
}
\end{equation}
Both modules, $\theta\mathrm{Z}_t^{\mathbf{\hat{R}}_2} P$ and
$\mathrm{Z}_t^{\mathbf{\hat{R}}_2}\theta P$, are obviously projective in
$\mathcal{O}^{\mathbf{\hat{R}}_2}$. Let $\theta'$ be the adjoint of $\theta$.
Then for any simple module $L\in \mathcal{O}^{\mathbf{\hat{R}}_2}$ we have
\begin{eqnarray*}
\mathrm{Hom}_{\mathfrak{g}}(\theta\mathrm{Z}_t^{\mathbf{\hat{R}}_2} P,L)&=&
\mathrm{Hom}_{\mathfrak{g}}(\mathrm{Z}_t^{\mathbf{\hat{R}}_2} P,\theta' L)\\
&=&\mathrm{Hom}_{\mathfrak{g}}(P,\theta' L)\\
&=&\mathrm{Hom}_{\mathfrak{g}}(\theta P, L)\\
&=&\mathrm{Hom}_{\mathfrak{g}}(\mathrm{Z}_t^{\mathbf{\hat{R}}_2}\theta P, L).
\end{eqnarray*}
Hence $\theta\mathrm{Z}_t^{\mathbf{\hat{R}}_2} P\cong
\mathrm{Z}_t^{\mathbf{\hat{R}}_2}\theta P$. In particular, by definition of
$\mathrm{Z}_t^{\mathbf{\hat{R}}_2}$, we have that $\theta K_{P}$ coincides with
the maximal submodule of $\theta P$, whose head consists only of simple modules
not in $\mathcal{O}^{\mathbf{\hat{R}}_2}$. In particular, the identity map on
$\theta P$ restricts to an isomorphism $\varphi:\theta
K_{P}\overset{\sim}{\longrightarrow} K_{\theta P}$, and induces the isomorphism
$\varphi':\theta\mathrm{Z}_t^{\mathbf{\hat{R}}_2} P
\overset{\sim}{\longrightarrow} \mathrm{Z}_t^{\mathbf{\hat{R}}_2}\theta P$.
It follows that cube on the left and the front, back, top and bottom faces of
\eqref{eqses1} commute. Therefore the face pointing to the right
commutes as well. This implies $\theta\mathrm{Z}_t^{\mathbf{\hat{R}}_2}\cong
\mathrm{Z}_t^{\mathbf{\hat{R}}_2}\theta$ since both functors are right exact.
\end{proof}

\begin{proposition}\label{p22}
Assume that $\mathscr{P}(\mathbf{R}_1)$ has a simple projective module $L$.
Then $\mathscr{P}(\mathbf{R}_2)$ has a simple projective module $L'$ given by
$\mathrm{Z}_t^{\mathbf{\hat{R}}_2}
\mathrm{F}\mathrm{i}_s^{\mathbf{\hat{R}}_1}L$.
\end{proposition}

To prove Proposition~\ref{p22} we will need a series of auxiliary statements.
We start with verifying that the expression $\mathrm{Z}_t^{\mathbf{\hat{R}}_2}
\,\mathrm{F}\,\mathrm{i}_s^{\mathbf{\hat{R}}_1}L$ makes sense, i.e. that it
gives a module:

\begin{lemma}\label{p22-l1}
Let $X=L$ or $X=L(x)$ for some $x\in\mathbf{{R}}_1$. Then
$\mathrm{Z}_t^{\mathbf{\hat{R}}_2}\,\mathrm{F}\,
\mathrm{i}_s^{\mathbf{\hat{R}}_1}X\in \mathcal{O}^{\mathbf{\hat{R}_2}}$.
\end{lemma}

\begin{proof}
The module $X$ does not belong to $\mathcal{O}^t_{\op{int}}$ because of the
conditions \eqref{Lcondition}. Hence by \cite[Proposition~4.2]{EW} we have
$\mathcal{L}_i\mathrm{Z}_t\, X=0$ for $i=0,2$ and $\mathcal{L}_1\mathrm{Z}_t\,
X\in \mathcal{O}^{t}_{\op{int}}$. Thus $\mathrm{F}X\in
\mathcal{O}^{t}_{\op{int}}$ and hence $\mathrm{Z}_t^{\mathbf{\hat{R}}_2}
\,\mathrm{F}\,\mathrm{i}_s^{\mathbf{\hat{R}}_1}X
\in \mathcal{O}^{\mathbf{\hat{R}_2}}$.
\end{proof}

\begin{lemma}\label{p22-l2}
\begin{enumerate}[(i)]
\item\label{p22-l2.1} $L':=\mathrm{Z}_t^{\mathbf{\hat{R}}_2}
\,\mathrm{F}\,\mathrm{i}_s^{\mathbf{\hat{R}}_1}L$ is a simple module.
\item\label{p22-l2.2} For each $L(x)$, $x\in\mathbf{{R}}_1$, the module
$\mathrm{Z}_t^{\mathbf{\hat{R}}_2}\,\mathrm{F}\,
\mathrm{i}_s^{\mathbf{\hat{R}}_1}L(x)$ is simple and has the form
$L(y)$ for some $y\in \mathbf{{R}}_2$. Moreover, the map
$\varphi:x\mapsto y$ is a bijection from $\mathbf{{R}}_1$ to $\mathbf{{R}}_2$.
\end{enumerate}
\end{lemma}

\begin{proof}
Let  $L(z)\in \mathcal{O}_0$ be the (unique) simple module which
translates to $L\in\mathscr{P}(\mathbf{R}_1)$ via translations to
walls (see e.g. \cite[4.12 (3)]{Ja2}). By \cite[Theorem~2]{MS3},
\cite[Theorem~6.3]{AS} and \cite[Theorem~7.8]{AS} we have
\begin{displaymath}
\mathcal{L}_1\mathrm{Z}_t\, L(z)\cong L(tz)\oplus\bigoplus_y L(y)^{a_y},
\end{displaymath}
where $tz\in\mathbf{{R}}_2$, and $a_y\neq 0$ implies that $y\neq tz$ but
both $y$ and $tz$ belong to the same left cell. Since the intersection of a
left and a right cell inside a  common two-sided cell consists of exactly one
element (by the  Robinson-Schensted correspondence, see e.g. \cite[3.1]{Sa}),
the later restrictions give that $a_y\neq 0$ implies $y\not\in \mathbf{R}_2$.
Hence  $\mathrm{Z}_t^{\mathbf{\hat{R}}_2}\mathcal{L}_1\mathrm{Z}_t\, L(z)$ is a
simple module. Translating this onto the walls we obtain that the module $L'$
is simple. This proves \eqref{p22-l2.1} and also \eqref{p22-l2.2} for the
module $L(z)$. For other $x\in \mathbf{{R}}_1$ the proof is just the same as
for $L(z)$. The fact that $\varphi:\mathbf{{R}}_1\to\mathbf{{R}}_2$
is a bijection follows from \cite[Section~4]{KLCoxeter}.
\end{proof}

\begin{lemma}\label{p22-l3}
\begin{enumerate}[(i)]
\item\label{p22-l3.1}
$L=\hat{\mathrm{Z}}_s^{\mathbf{\hat{R}}_1}\,
\mathrm{G}\,\mathrm{i}_t^{\mathbf{\hat{R}}_2}\,L'$.
\item\label{p22-l3.2}
For any $x\in \mathbf{{R}}_1$ we have
$L(x)=\hat{\mathrm{Z}}_s^{\mathbf{\hat{R}}_1}\,
\mathrm{G}\,\mathrm{i}_t^{\mathbf{\hat{R}}_2}\,L(\varphi(x))$.
\end{enumerate}
\end{lemma}

\begin{proof}
Analogous to the proof of Lemma~\ref{p22-l2}.
\end{proof}

As $L\in\mathscr{P}(\mathbf{R}_1)$, the category $\mathscr{P}(\mathbf{R}_1)$
is equivalent to the additive closure of the category with objects  $L\otimes
E$, where $E$ runs through all finite-dimensional $\mathfrak{g}$-modules.  Set $\tilde{\mathrm{F}}=\mathrm{Z}_t^{\mathbf{\hat{R}}_2}
\mathrm{F}\mathrm{i}_s^{\mathbf{\hat{R}}_1}$,
$\tilde{\mathrm{G}}=\hat{\mathrm{Z}}_s^{\mathbf{\hat{R}}_1}
\mathrm{G}\mathrm{i}_t^{\mathbf{\hat{R}}_2}$ and
$\mathscr{Q}=\tilde{\mathrm{F}}\mathscr{P}(\mathbf{R}_1)$.

\begin{lemma}\label{l23}
\begin{enumerate}[(i)]
\item\label{l23.1} The functors $\tilde{\mathrm{F}}$ and  $\tilde{\mathrm{G}}$
define mutually inverse equivalences between $\mathscr{P}(\mathbf{R}_1)$ and
$\mathscr{Q}$.
\item\label{l23.2} $\mathscr{Q}$ is equivalent to the additive  closure
of the category with objects $L'\otimes E$, where $E$ runs through all
finite-dimensional $\mathfrak{g}$-modules.
\end{enumerate}
\end{lemma}

\begin{proof}
We have already seen that $\tilde{\mathrm{F}}L=L'$ and
$\tilde{\mathrm{G}}L'=L$. By Lemma~\ref{l21} we thus have that
\begin{equation}\label{eq12321}
\tilde{\mathrm{G}}\tilde{\mathrm{F}}\,(E\otimes L)\cong E\otimes L\quad
\text{ and }\quad
\tilde{\mathrm{F}}\tilde{\mathrm{G}}\,(E\otimes L')\cong E\otimes L'
\end{equation}
for any finite-dimensional $\mathfrak{g}$-module $E$.  By definition, we
have the adjoint pair $(\tilde{\mathrm{F}},\tilde{\mathrm{G}})$.
Consider the adjunction morphisms $\mathrm{adj}:\tilde{\mathrm{F}}\tilde{\mathrm{G}}\to\mathrm{ID}$
and $\overline{\mathrm{adj}}:\mathrm{ID}\to
\tilde{\mathrm{G}}\tilde{\mathrm{F}}$. Then the adjunction property
says that $\mathrm{adj}_{\tilde{\mathrm{F}}()}\circ
\tilde{\mathrm{F}}(\overline{\mathrm{adj}})=\mathrm{id}$. In particular
$\mathrm{adj}_{E\otimes L'}$ must be surjective, hence an isomorphism
by \eqref{eq12321}. Similarly $\overline{\mathrm{adj}}_{E\otimes L}$
is an isomorphism. This proves statement \eqref{l23.1} and
statement \eqref{l23.2} follows then from \eqref{l23.1} and Lemma~\ref{l21}.
\end{proof}

Let now $\mathscr{Y}_1$ denote the full subcategory of $\mathcal{O}_0$, whose
objects are the $P^{\mathbf{\hat{R}_1}}(x)$ and the $L(x)$, $x\in \mathbf{R}_1$.
Denote further by $\mathscr{Y}_2$ the full subcategory of
$\mathcal{O}_0$ whose objects are $\tilde{\mathrm{F}}\,
P^{\mathbf{\hat{R}_1}}(x)$, $x\in \mathbf{R}_1$, and
$L(y)$, $y\in \mathbf{R}_2$. Lemma~\ref{l23} can be refined as follows:

\begin{lemma}\label{p22-l4}
The functors $\tilde{\mathrm{F}}$ and $\tilde{\mathrm{G}}$ induce mutually
inverse equivalences of categories between $\mathscr{Y}_1$ and
$\mathscr{Y}_2$.
\end{lemma}

\begin{proof}
By definition and Lemma~\ref{l23}, $\tilde{\mathrm{F}}\,
P^{\mathbf{\hat{R}_1}}(x)\in \mathscr{Y}_2$ for all $x\in \mathbf{R}_1$, and
$\tilde{\mathrm{G}}\tilde{\mathrm{F}}\, P^{\mathbf{\hat{R}_1}}(x)\in
\mathscr{Y}_1$ for all $x\in \mathbf{R}_1$. Analogously to the proof of
Lemma~\ref{p22-l2} one shows that for each $x\in \mathbf{R}_1$ we have
$\tilde{\mathrm{F}}\, L(x)\cong L(y)$ for some $y\in \mathbf{R}_2$, and that
for each  $y\in \mathbf{R}_2$ we have $\tilde{\mathrm{G}}\, L(y)\cong L(x)$
for some $x\in \mathbf{R}_1$. Hence $\tilde{\mathrm{F}}:\mathscr{Y}_1\to
\mathscr{Y}_2$ and $\tilde{\mathrm{G}}:\mathscr{Y}_2\to \mathscr{Y}_1$.
That these functors are mutually inverse equivalences is proved in the same
way as in Lemma~\ref{l23}.
\end{proof}

For $x\in \mathbf{R}_1$ set $N_x=\tilde{\mathrm{F}}\,
P^{\mathbf{\hat{R}_1}}(x)$.

\begin{corollary}\label{p22-c5}
For every $x\in \mathbf{R}_1$ we have $P^{\mathbf{\hat{R}_2}}(\varphi(x))
\tto N_x$.
\end{corollary}

\begin{proof}
Using the Lemmas~\ref{p22-l1}-\ref{p22-l4}, for any $x\in \mathbf{R}_1$
and  $y\in \mathbf{R}_2$ we have
\begin{displaymath}
\begin{array}{rcl}
\mathrm{Hom}_{\mathfrak{g}}(N_x,L(y))&=&
\mathrm{Hom}_{\mathfrak{g}}(\tilde{\mathrm{F}}\,
P^{\mathbf{\hat{R}_1}}(x),L(y))\\
&=&\mathrm{Hom}_{\mathfrak{g}}( P^{\mathbf{\hat{R}_1}}(x),
\tilde{\mathrm{G}}\,L(y))\\
&=&\mathrm{Hom}_{\mathfrak{g}}( P^{\mathbf{\hat{R}_1}}(x),
L(\varphi^{-1}(y)))\\
&=&
\begin{cases}
\mathbb{C},&\varphi(x)=y;\\0,& \text{otherwise}.\end{cases}
\end{array}
\end{displaymath}
and the claim follows.
\end{proof}

\begin{lemma}\label{p22-l6}
\begin{enumerate}[(i)]
\item\label{p22-l6.1} For $x,w\in W$ we have $\theta_w^l\theta_x^l\cong
\oplus_{y\geq_{\mathsf{L}}w}(\theta_y^l)^{m_y}$.
\item\label{p22-l6.2} Let $x,w\in W$ be such that $x<_{\mathsf{R}}w$. Then
$\theta_w^l L(x)=0$.
\item\label{p22-l6.3} For each $w\in W$ there exists $x\in W$ such that
$x\sim_{\mathsf{R}}w$ and $\theta^l_w L(x)\neq 0$.
\end{enumerate}
\end{lemma}

\begin{proof}
To prove the first statement we use some ideas from the proof of
\cite[Theorem~11]{Ma3}. Denote by $\sigma$ the unique
anti-automorphism of $\mathds{H}$, which maps $H_w$ to $H_{w^{-1}}$
(and hence $\underline{H}_w$ to $\underline{H}_{w^{-1}}$) for each
$w\in W$. Using \eqref{formula1}  we have:
\begin{multline*}
\underline{H}_w\underline{H}_x=
\sigma(\sigma(\underline{H}_w\underline{H}_x))
=\sigma(\sigma(\underline{H}_x)\sigma(\underline{H}_w))=\\=
\sigma(\underline{H}_{x^{-1}}\underline{H}_{w^{-1}})
=\sum_{y^{-1}\geq_\mathsf{R}w^{-1}}\sigma(a_y\underline{H}_{y^{-1}})
=\sum_{y\geq_\mathsf{L}w}a_y\underline{H}_{y}.
\end{multline*}
Now \eqref{p22-l6.1} follows from Proposition~\ref{prop3}.

Let  $x,w\in W$ be such that $x<_{\mathsf{R}}w$. We have
$P(x)\cong \theta^l_x\Delta(e)\tto L(x)$. Using \eqref{p22-l6.1} we have
$\theta^l_wP(x)\cong \theta^l_w\theta^l_x\Delta(e)\cong
\oplus_{y\geq_{\mathsf{L}}w}P(y)^{m_y}$. At the same time by
Proposition~\ref{prop5}\eqref{prop5.1}, the head of  $\theta^l_w L(x)$ can
contain only $L(y)$ such that  $y\leq_{\mathsf{R}}x$. Hence we have
$y\leq_{\mathsf{R}}x<_{\mathsf{R}}w\leq_{\mathsf{L}}y$, which is not possible.
Therefore  $\theta^l_w L(x)=0$ proving \eqref{p22-l6.2}.

Let $\mathbf{R}$ be the right cell of $w$. Using Lemma~\ref{l21}, we have
\begin{displaymath}
\theta^l_w\mathrm{Z}^{\mathbf{\hat{R}}}\Delta(e)\cong
\mathrm{Z}^{\mathbf{\hat{R}}}\theta^l_w\Delta(e)\cong
\mathrm{Z}^{\mathbf{\hat{R}}}P(w)\cong
P^{\mathbf{\hat{R}}}(w)\neq 0.
\end{displaymath}
Hence $\theta^l_w L(x)\neq 0$ for some simple subquotient $L(x)$ of
$\mathrm{Z}^{\mathbf{\hat{R}}}\Delta(e)$. In particular, $x\leq_{\mathsf{R}}w$
(thanks to the definition of $\mathrm{Z}^{\mathbf{\hat{R}}}$)
and then $x\sim_{\mathsf{R}}w$ follows from \eqref{p22-l6.2}.
\end{proof}

\begin{lemma}\label{p22-l8}
There exists $z\in \mathbf{R}_1$ such that $N_z\in \mathscr{P}(\mathbf{R}_2)$.
\end{lemma}

\begin{proof}
We choose $w,y\in \mathbf{R}_2$ such that $\theta_w^l L(y)\neq 0$
(see Lemma~\ref{p22-l6}\eqref{p22-l6.3}). By Corollary~\ref{p22-c5},
there exists some $x\in \mathbf{R}_1$ such that
$P^{\mathbf{\hat{R}_2}}(y)\tto N_x$. Let $K$ be the kernel of the
latter map. Consider the short exact sequence $K'\hookrightarrow
K\tto K''$, where $K''$ is the maximal quotient of $K$, which
contains only simple subquotients  of the form $L(z)$,
$z<_{\mathsf{R}}w$. By Lemma~\ref{p22-l6}\eqref{p22-l6.2} we have
$\theta_w^l K'\cong \theta_w^l K$. Hence we have the short exact
sequence of the form
\begin{equation}\label{p22-l8-e1}
\theta_w^l K'\hookrightarrow \theta_w^lP^{\mathbf{\hat{R}_2}}(y)
\tto \theta_w^lN_x.
\end{equation}
Note that $\theta_w^lN_x\neq 0$ since $\theta$ is exact,
$\theta_w^lL(y)\neq 0$ and $L(y)$ is the head of $N_x$. If
$\theta_w^l K'=0$, we immediately get that  $0\neq \theta_w^lN_x\in
\mathscr{P}(\mathbf{R}_2)$. But the additive category, generated by
indecomposable modules $N_z$, $z\in\mathbf{R}_1$,  is stable with
respect to projective functors by Lemma~\ref{l21}. This implies that
$N_z\in \mathscr{P}(\mathbf{R}_2)$ for some $z\in \mathbf{R}_1$.

Assume hence that $\theta_w^l K'\neq 0$ and consider an arbitrary short exact
sequence of the form $M'\hookrightarrow \theta_w^l K'\tto M''$ such that $M''$
is simple. Then $M''\cong L(v)$ for some $v\in\mathbf{R}_2$. If we factor $M'$
out  in \eqref{p22-l8-e1} we obtain the short exact sequence
\begin{equation}\label{p22-l8-e2}
L(v)\hookrightarrow X \tto \theta_w^lN_x,
\end{equation}
where $X=\theta_w^lP^{\mathbf{\hat{R}_2}}(y)/M'$. By Corollary~\ref{p22-c5},
The heads of $X$ and  $\theta_w^lN_x$ are isomorphic. Hence the sequence
\eqref{p22-l8-e2} is not split. Apply now the functor $\tilde{\mathrm{G}}$ to
the sequence \eqref{p22-l8-e2}, which basically reduces to the application of
the functor $\mathcal{L}_1\mathrm{Z}_t$ because of the definition of
$\tilde{\mathrm{G}}$. As $\mathcal{L}_2\mathrm{Z}_t\theta_w^lN_x=0$ and
$\mathcal{L}_0\mathrm{Z}_tL(v)=0$ (this follows for example from
Lemma~\ref{p22-l4} and the definition of $\tilde{\mathrm{G}}$), we obtain a
short exact sequence
\begin{equation}\label{p22-l8-e3}
\tilde{\mathrm{G}}L(v)\hookrightarrow \tilde{\mathrm{G}}X \tto
\tilde{\mathrm{G}}\theta_w^lN_x,
\end{equation}
in particular, $\tilde{\mathrm{G}}X\in \mathcal{O}^{\mathbf{\hat{R}_1}}$.
Analogously one shows that $\tilde{\mathrm{F}}\tilde{\mathrm{G}}X\in
\mathcal{O}^{\mathbf{\hat{R}_2}}$, which, together with
Lemma~\ref{p22-l4}, implies that the adjunction morphism induces an
isomorphism $\tilde{\mathrm{F}}\tilde{\mathrm{G}}X\cong X$, and thus
the sequence \eqref{p22-l8-e2} is obtained from the sequence
\eqref{p22-l8-e3} by applying  $\tilde{\mathrm{F}}$. However,
the sequence \eqref{p22-l8-e3} splits as
$\tilde{\mathrm{G}}\theta_w^lN_x$ is projective in
$\mathcal{O}^{\mathbf{\hat{R}_1}}$. Therefore \eqref{p22-l8-e2} must be
split as well, a contradiction. Hence $\theta_w^l K'\neq 0$ is not possible.
This completes the proof.
\end{proof}

\begin{proof}[Proof of Proposition~\ref{p22}.]
To prove Proposition~\ref{p22} it is enough to show that
$\mathscr{Q}=\mathscr{P}(\mathbf{R}_2)$. Let  $\mathscr{Q}_0$ and
$\mathscr{P}(\mathbf{R}_2)_0$ denote the intersections of $\mathcal{O}_0$
with $\mathscr{Q}$ and $\mathscr{P}(\mathbf{R}_2)$  respectively.
The definition of $\mathscr{P}(\mathbf{R}_2)$ and Lemma~\ref{l23}
imply that it is even enough to show that
$\mathscr{Q}_0=\mathscr{P}(\mathbf{R}_2)_0$. From Lemma~\ref{p22-l8}
we know that $\mathscr{Q}_0\cap \mathscr{P}(\mathbf{R}_2)_0$ is not trivial.
As $\mathscr{Q}_0$ is additively closed by Lemma~\ref{l23}\eqref{l23.2}
we have that $\mathscr{Q}_0$ contains some indecomposable projective from
$\mathscr{P}(\mathbf{R}_2)_0$. Applying projective functors and
Theorem~\ref{thm5} we get that $\mathscr{Q}_0$ must contain all indecomposable
projectives from $\mathscr{P}(\mathbf{R}_2)_0$. But by Lemma~\ref{l23}
the categories $\mathscr{Q}_0$ and  $\mathscr{P}(\mathbf{R}_2)_0$ contain
the same number  of pairwise non-isomorphism indecomposable modules. Hence
$\mathscr{Q}_0=\mathscr{P}(\mathbf{R}_2)_0$. This completes the proof.
\end{proof}

Now we are prepared to prove Theorem~\ref{thm6}.

\begin{proof}[Proof of Theorem~\ref{thm6}.]
Assume first that $\mathbf{R}_1$ is of the form described in
Remark~\ref{rem2}. Then $\mathscr{P}(\mathbf{R}_1)$ has a simple
projective module by \cite[Section~3.1]{IS}. Let now
$\mathbf{R}_2$ be any other right cell in the same two-sided cell
as $\mathbf{R}_1$. By  \cite[Proof of Theorem~1.4]{KLCoxeter} there is
a sequence, $\mathbf{R}_1=\mathbf{R}^{(0)}$, $\mathbf{R}^{(2)}$,\dots,
$\mathbf{R}^{(k)}=\mathbf{R}_2$, such that
$(\mathbf{R}^{(i)},\mathbf{R}^{(i+1)})$ satisfies the condition
\eqref{Lcondition} for each $i=0,\dots,k-1$. Inductively applying Lemma~\ref{l23} and Proposition~\ref{p22} provides an equivalence between
$\mathscr{P}(\mathbf{R}_1)$ and $\mathscr{P}(\mathbf{R}_2)$.
This of course induces an equivalence of abelian categories.
\end{proof}

\subsection{Consequences}\label{s4.2}

Let $\mathbf{R}$ be a right cell of $S_n$. From
Theorem~\ref{thm6} and Remark~\ref{rem2} one can deduce the following facts:

\begin{enumerate}[(I)]
\item\label{cons1}
The Koszul grading on the algebra $A$ (\cite{SoKipp}) turns
$\mathrm{End}_{\mathcal{O}_0}(\oplus_{w\in \mathbf{R}}P^{\mathbf{\hat{R}}}(w))$
into a positively graded self-injective symmetric algebra,
\cite[Theorem~5.4]{MS2}.
\item\label{cons2}
The center of $\mathrm{End}_{\mathcal{O}_0}
(\oplus_{w\in \mathbf{R}}P^{\mathbf{\hat{R}}}(w))$ is
isomorphic to the cohomology algebra of the associated
Springer fiber, see \cite[Theorem~2]{Br} and
\cite[Theorem~4.1.1]{St2}.
\item\label{cons3}
For each $w\in \mathbf{R}$ there is a finite-dimensional
$\mathfrak{g}$-module $E$ such that each
$P^{\mathbf{\hat{R}}}(x)$, $x\in \mathbf{R}$, is a direct
summand of $E\otimes L(w)$. This follows from
\cite[Proposition~4.3(ii)]{Irself}.
\item\label{cons4}
The projective modules in $\mathscr{P}(\mathbf{R})$ have all the
same Loewy lengths (\cite[Theorem~5.2]{MS2}).
\end{enumerate}

\subsection{Counter-examples}\label{s4.4}

Perhaps the most remarkable feature of Theorem~\ref{thm6} is that
there is no way to extend this result to the categories
$\mathcal{O}_0^{\mathbf{\hat{R}}}$. For two right cells satisfying
the condition of Remark~\ref{rem2} this was already pointed out in
\cite[Propositions~6]{Kh}. At the same time, in
\cite[Propositions~7]{Kh}, it was shown that the corresponding
$\mathcal{O}_0^{\mathbf{\hat{R}}}$'s are derived equivalent. Even
this weaker statement is not true in the general case. For example,
take $W=S_4$, generated by the simple reflections $s,t,r$ such that
$sr=rs$. Take the two right cells $\mathrm{R}_1=\{sr,srt\}$ and
$\mathrm{R}_2=\{tsr,tsrt\}$. Then we have
$\mathbf{\hat{R}}_1=\{e,s,r,ts,tr,sr,rts,str,srt\}$ whereas
$\mathbf{\hat{R}}_2=\{e,t,ts,tr,tsr,tsrt\}$. In particular, the
categories $\mathcal{O}_0^{\mathbf{\hat{R}}_1}$ and
$\mathcal{O}_0^{\mathbf{\hat{R}}_2}$ have different numbers of
simple modules; hence they cannot be derived equivalent.

For right cells ${\mathbf{\hat{R}}}$ satisfying the condition of
Remark~\ref{rem2}, the categories $\mathcal{O}^{\mathbf{\hat{R}}}$
are special amongst the categories associated with right cells: they
are equivalent to the principal block of some parabolic category
$\mathcal{O}$, in particular are highest weight categories (i.e.
described by quasi-hereditary algebras), see \cite{RC}. This is not
true for arbitrary right cells. The smallest such example is again
the case $W=S_4$ with $\mathrm{R}=\{t,ts,tr\}$. In this case
$\mathbf{\hat{R}}= \{e,t,ts,tr\}$ and we have the following graded
filtrations of projective and standard modules in
$\mathcal{O}_0^{\mathbf{\hat{R}}}$:
\begin{displaymath}
\begin{array}{|c||c|c|c|c|}
\hline
w & e & t & ts & tr\\
\hline\hline
P(w)&
\begin{array}{c}e\\t\\\text{\hspace{2mm}}\end{array}&
\begin{array}{ccc}&t&\\ts&e&tr\\&t&\end{array}&
\begin{array}{c}ts\\t\\ts\end{array}&
\begin{array}{c}tr\\t\\tr\end{array}
\\
\hline
\Delta(w)&
\begin{array}{c}e\\t\end{array}&
\begin{array}{ccc}&t&\\ts&&tr\end{array}&
\begin{array}{c}ts\\\text{\hspace{2mm}}\end{array}&
\begin{array}{c}tr\\\text{\hspace{2mm}}\end{array}\\
\hline
\end{array}
\end{displaymath}
We see that not all projective modules have standard filtrations and
hence $\mathcal{O}_0^{\mathbf{\hat{R}}}$
is not a highest weight category.

\section{Tensor products and parabolic induction}\label{s5}

In this section we show how one can categorify some standard
representation theoretical operations like tensor products and
parabolic induction. As application we categorify induced cell
modules. Up to equivalence, the resulting categories depend only on
the isomorphism class of the cell module, not on the actual cell
module itself.

\subsection{Inner and outer tensor products}\label{s5.1}

Let  $W$ and $W'$ be arbitrary finite Weyl groups with sets of simple
reflections $S$ and $S'$. Let $\mathds{H}$, $\mathds{H}'$ be the corresponding
Hecke algebras. If $M$ is a right $\mathds{H}$-module and $M'$ is a right
$\mathds{H}'$-module then the {\it outer tensor product} $M\boxtimes M'$ is the
right $\mathds{H}\otimes\mathds{H}'$-module whose underlying space is
$M\otimes M'$ and the module structure is given by $m\otimes m'(h\otimes
h')=mh\otimes m'h'$ for $m\in M$, $m'\in M'$, $h\in\mathds{H}$ and
$h'\in\mathds{H}'$.

Let now $X$ and $Y$ be right $\mathds{H}$-modules. Then the {\it inner tensor
product} $X\otimes Y$ is the right $\mathds{H}$-module whose underlying space
is $X\otimes Y$ and the module structure is given by
$(x\otimes y).h=x.h\otimes y.h$ for $x\in X$, $y\in Y$ and $h\in\mathds{H}$.

Given two categories $\ccC_1$ and $\ccC_2$ let $\ccC_1\oplus\ccC_2$ be the
category with objects being pairs $(C_1,C_2)$, where $C_i$ is an object in
$\ccC_i$, and the morphisms from an object $(A_1,A_2)$ to an object $(B_1,B_2)$
being pairs of morphisms $(f_1,f_2)$, where $f_i:A_i\rightarrow B_i$ for
$i=1,2$. We assume that each of these categories is either equivalent to a
module categories over some finite dimensional algebra $A$ or at least
equivalent to its (bounded) derived category. Then $\op{Gr}(\ccC_1\otimes\ccC_2)
\cong \op{Gr}(\ccC_1)\otimes_\mZ \op{Gr}(\ccC_2)$ and hence also
$[\cC_1\oplus\cC_2]\cong [\ccC_1]\otimes [\ccC_2]$. Given two functors
$F_i:\ccC_i\rightarrow\ccC_i$, $i=1,2$, then we denote by $F_1\boxtimes F_2$
the endofunctor of $\ccC_1\oplus\ccC_2$ which maps $(A_1,A_2)$ to
$(F_1(A_1),F_2(A_2))$ and $(f_1,f_2)$ to $(F_1(f_1),F_2(f_2))$. The following
result gives a categorification of the outer and inner tensor products:

\begin{proposition}[Tensor products]\label{tensorproducts}
Assume we are given a right $\mathds{H}$-module $M$ and a  right
$\mathds{H}'$-module $M'$ together with a categorification
$(\ccS,\mathcal{E},\{\mathrm{F}_s\}_{s\in S})$
of $M$ with respect to  the generators $\underline{H}_s$, $s\in S$, of
$\mathds{H}$; and a categorification  $(\ccS',\mathcal{E}',
\{\mathrm{F}'_{s'}\}_{s'\in S'})$ of $M'$ with respect to the generators
$\underline{H}_{s'}$, where $s'\in S'$, of $\mathds{H}'$. Then we have:
\begin{enumerate}[(i)]
\item\label{tensorproducts.1}
The tuple
\begin{displaymath}
(\ccS\oplus\ccS', \mathcal{E}\otimes\mathcal{E}',
\{\mathrm{F}_s\boxtimes\mathrm{F}'_{s'}\}_{s\in S,s'\in S'})
\end{displaymath}
is a categorification of $M\boxtimes M'$ with respect to the
generators $\underline{H}_{s}\otimes \underline{H}_{s'}$, $s\in S$,
$s'\in S'$.
\item\label{tensorproducts.2}
If both $M$, and $M'$ are right $\mathds{H}$-modules then
\begin{displaymath}
(\ccS\oplus\ccS',\mathcal{E}\otimes\mathcal{E}',
\{\mathrm{F}_{s}\boxtimes\mathrm{F}'_{s}\}_{s\in S})
\end{displaymath}
is a categorification of $M\otimes M'$ with respect to the generators
$\underline{H}_s$, where $s\in S$, of $\mathds{H}$.
\end{enumerate}
\end{proposition}

\begin{proof}
This follows directly from the definitions.
\end{proof}

\subsection{Examples of parabolic induction}\label{s5.2}

Let now $W'$ be a parabolic subgroup of $W$ which corresponds to a subset
$S'\subset S$. Let $\mathds{H}'=\mathds{H}(W',S')$ be the corresponding
subalgebra of $\mathds{H}$, and let $M$ be a (right) $\mathds{H}'$-module.
The purpose of this section is to give a categorification of the induced
module $\mathrm{Ind}_{\mathds{H}'}^{\mathds{H}} M=
M\otimes_{\mathds{H}'}\mathds{H}$, where $M$ is a cell module over $\mathds{H}'$. We start by recalling examples from the literature.

\subsubsection{Sign parabolic module}

The assignment $H_s\mapsto -v$ for all $s\in S'$ defines a surjection
$\mathds{H}'\tto \mathbb{Z}[v,v^{-1}]$ and hence defines on
$\mathbb{Z}[v,v^{-1}]$ the structure of an $\mathds{H}'$--bimodule. Consider
the {\em sign parabolic} $\mathds{H}$-module $\mathcal{N}=
\mathbb{Z}[v,v^{-1}]\otimes_{\mathds{H}'}\mathds{H}$. The set $\{N_x=1\otimes
H_x\}$, where $x$ runs through the set $(W'\backslash W)_{short}$ of shortest
coset representatives in $W'\backslash W$, forms a basis  of $\mathcal{N}$.
The action of $\underline{H}_s$, $s\in S$, in this basis is given by
(see \cite[Section~3]{SoKipp}):
\begin{displaymath}
N_x\underline{H}_s=
\begin{cases}
N_{xs}+vN_x, & \text{ if } xs\in (W'\backslash W)_{short}, xs>x;\\
N_{xs}+v^{-1}N_x, & \text{ if } xs\in (W'\backslash W)_{short}, xs<x;\\
0, & \text{ if } xs\not\in (W'\backslash W)_{short}.
\end{cases}
\end{displaymath}
It is easy to see that the specialization $v=1$ gives the $W$-module
$\mathrm{Ind}_{W'}^W M$, where $M$ is the {\em sign} $W'$-module, that is
$M=\mathbb{Z}$ with the alternating action $1\, s=-1$ for all $s\in S'$.

\subsubsection*{Its categorification}
Let $\mathfrak{p}\supseteq \mathfrak{b}$ be the parabolic subalgebra of
$\mathfrak{g}$ corresponding to $S'$. Let further
$\mathcal{O}_0^{\mathfrak{p}}$ be the locally $\mathfrak{p}$-finite part of
$\mathcal{O}_0$ (in the sense of \cite{RC}). This is the full extension closed
subcategory of $\mathcal{O}_0$, generated by the simple modules $L(w)$, $w\in
(W'\backslash W)_{short}$. Finally, let
$\mathcal{O}_0^{\mathfrak{p},\mathbb{Z}}$ be the graded version of
$\mathcal{O}_0^{\mathfrak{p}}$ (as defined in \cite{BGS}). Let
$\Delta^{\mathfrak{p}}(w)$ denote the corresponding standard graded lift of the
generalized Verma module, i.e. the corresponding standard module in
$\mathcal{O}_0^{\mathfrak{p},\mathbb{Z}}$ with head concentrated in degree $0$.
The category $\mathcal{O}_0^{\mathfrak{p}}$ has finite homological dimension,
and hence we have a unique isomorphism $\mathcal{E}^{\mathfrak{p}}$ of
$\mathbb{Z}[v,v^{-1}]$-modules as follows:
\begin{eqnarray*}
\mathcal{E}^{\mathfrak{p}}: \quad \mathcal{N} & \overset{\sim}{\longrightarrow}
& \left[\mathcal{O}_0^{\mathfrak{p},\mathbb{Z}}\right]
\\  N_w & \mapsto & \left[\Delta^{\mathfrak{p}}(w)\right].
\end{eqnarray*}
The following result is well-known (see for example
\cite[Proposition~1.5]{StDuke}):

\begin{proposition}\label{p31}
$(\mathcal{O}_0^{\mathfrak{p},\mathbb{Z}},\mathcal{E}^{\mathfrak{p}},
\{\theta_s^l\}_{s\in S})$ is a categorification of $\mathcal{N}$ with
respect to the generators $\underline{H}_s$, $s\in S$.
\end{proposition}

\subsubsection{Permutation parabolic module}
\label{permmodule}

The assignment $H_s\mapsto v^{-1}$ for all $s\in S'$
defines a surjection $\mathds{H}'\tto \mathbb{Z}[v,v^{-1}]$
and hence determines on $\mathbb{Z}[v,v^{-1}]$ the structure of an
$\mathds{H}'$--bimodule. The {\em permutation parabolic}
$\mathds{H}$-module is defined as follows: $\mathcal{M}=
\mathbb{Z}[v,v^{-1}]\otimes_{\mathds{H}'}\mathds{H}$.
There is the standard basis of
$\mathcal{M}$ given by $\{M_x=1\otimes H_x\}$, where $x$ runs through $(W'\backslash W)_{short}$. The action of $\underline{H}_s$,
$s\in S$, in this basis is given as follows (see \cite[Section~3]{SoKipp}):
\begin{displaymath}
M_x\underline{H}_s=
\begin{cases}
M_{xs}+vM_x, & \text{ if } xs\in (W'\backslash W)_{short}, xs>x;\\
M_{xs}+v^{-1}M_x, & \text{ if } xs\in (W'\backslash W)_{short}, xs<x;\\
(v+v^{-1})M_x, & \text{ if } xs\not\in (W'\backslash W)_{short}.
\end{cases}
\end{displaymath}
It is easy to see that the specialization $v=1$ gives the
$W$-module $\mathrm{Ind}_{W'}^W M$, where $M$ is the
{\em trivial} $W'$-module, that is $M=\mathbb{Z}$ with the trivial action
$1\, s=1$ for all $s\in S'$. The module $\mathrm{Ind}_{W'}^W M$
is usually called the {\em permutation module}, see \cite[2.1]{Sa}.

\subsubsection*{Its categorification} Let $\mathfrak{p}$ be as in the previous example. Let
$\mathscr{P}(\mathfrak{p})$ be the additive category, closed with
respect to the shift of grading, and generated by the indecomposable
projective modules $\mathtt{P}(w)\in\cO_0^\mZ$, where $w$ runs
through the set $(W'\backslash W)_{long}$ of longest coset
representatives in $W'\backslash W$. The category
$\mathcal{O}_0^{\mathfrak{p}-pres,\mathbb{Z}}=
\overline{\mathscr{P}(\mathfrak{p})}$ is the graded version of the
category $\mathcal{O}_0^{\mathfrak{p}-pres}$ from \cite{MS} (see
also Subsection~\ref{s3.244}). The simple objects of
$\mathcal{O}_0^{\mathfrak{p}-pres}$ are in a natural bijection with
$w\in (W'\backslash W)_{long}$. For $w\in (W'\backslash W)_{long}$
denote by $\Delta^{\mathfrak{p}-pres}(w)$ the standard object of
$\mathcal{O}_0^{\mathfrak{p}-pres,\mathbb{Z}}$ corresponding to $w$
and with the head concentrated in degree $0$ (\cite[Theorem 2.16,
Lemma 7.2]{MS}). Let $w'_0$ be the longest element of $W'$. All
this defines a unique homomorphism $\mathcal{E}^{\mathfrak{p}-pres}$
of $\mathbb{Z}[v,v^{-1}]$-modules as follows:
\begin{eqnarray*}
\mathcal{E}^{\mathfrak{p}-pres}: \quad
\mathcal{M} & \overset{\sim}{\longrightarrow} & \left[\mathcal{O}_0^{\mathfrak{p}-pres,\mathbb{Z}}\right]
\\ M_{w'_0w} & \mapsto & \left[\Delta^{\mathfrak{p}-pres}(w)\right].
\end{eqnarray*}
The category $\mathcal{O}_0^{\mathfrak{p}-pres,\mathbb{Z}}$ does not have
finite homological dimension in general, however, projective dimension of
all standard modules is finite. Hence we have
$\left[\Delta^{\mathfrak{p}-pres}(w)\right]\in
[\mathscr{P}(\mathfrak{p})]_{\oplus}$ for all $w\in (W'\backslash W)_{long}$.
This can be extended to the following (see for example \cite[Theorem~7.7]{MS}):

\begin{proposition}\label{p32}
\begin{enumerate}[(i)]
\item\label{p32.1}
$(\mathcal{O}_0^{\mathfrak{p}-pres,\mathbb{Z}},
\mathcal{E}^{\mathfrak{p}-pres}, \{\theta_s^l\}_{s\in S})$ is a
precategorification where\-as
$(\mathscr{P}(\mathfrak{p}),\mathcal{E}^{\mathfrak{p}-pres},
\{\theta_s^l\}_{s\in S})$ is a categorification  of $\mathcal{M}$ with
respect to the generators $\underline{H}_s$, $s\in S$.
\item\label{p32.2}
The homomorphism $\mathcal{E}^{\mathfrak{p}-pres}$ extends uniquely to the
${\mathbb{Z}((v))}$-categorification
$(\mathcal{O}_0^{\mathfrak{p}-pres,\mathbb{Z}},
\mathcal{E}^{\mathfrak{p}-pres},\{\theta_s^l\}_{s\in S})$
of  $\mathcal{M}^{\mathbb{Z}((v))}$ with respect to the
generators $\underline{H}_s$, $s\in S$.
\end{enumerate}
\end{proposition}

\subsection{The `unification': The category $\mathcal{O}\{\mathfrak{p},\mathscr{A}\}$}\label{s5.3}

In Subsection~\ref{s5.2} we used certain parabolic categories to
categorify the sign module, but also used categories of presentable
modules to categorify the permutation parabolic modules. Both depend
on a fixed parabolic $\p\subset\mg$. In this section we actually
want to put these two approaches under one roof using a series
$\mathcal{O}\{\mathfrak{p},\mathscr{A}\}$ of categories, depending
on the (fixed) $\p$, $\g$ and a (varying) category $\mathscr{A}$.
The categorifications from Subsection~\ref{s5.2} will then emerge
for a special choice of $\mathscr{A}$.

The categories $\mathcal{O}\{\mathfrak{p},\mathscr{A}\}$ were first
introduced in \cite{FKM} - as certain parabolic generalizations of
the category $\mathcal{O}$ which led to properly stratified
algebras. The setup was afterwards extended in \cite[6.2]{Ma} to
include general stratified algebras (in the sense of \cite{CPS}).
Here we present a slight variation of the original definition. This
variation seems to be more natural for us, and is better adapted to
the examples we work with.

Let $\tilde{\mathfrak{a}}$ be a reductive complex finite dimensional Lie
algebra with semisimple part $\mathfrak{a}$ and center
$\mathfrak{z}(\tilde{\mathfrak{a}})$. Let $\mathscr{A}$ be a full
subcategory of the category of finitely generated $\tilde{\mathfrak{a}}
$-modules. Then $\mathscr{A}$ is called an {\it admissible} category (of $\tilde{\mathfrak{a}}$-modules) if the following holds:
\begin{enumerate}[(L1)]
\item\label{lll2} $\mathscr{A}$ is stable under $E\otimes _-$ for each simple
finite dimensional $\tilde{\mathfrak{a}}$-module $E$;
\item\label{lll3} the action of $Z(\tilde{\mathfrak{a}})$ gives
a decomposition of  $\mathscr{A}$ into a direct sum of full subcategories,
each of which is equivalent to a module category over a
finite-dimensional self-injective associative algebra;
\item\label{lll4} the action of $\mathfrak{z}(\tilde{\mathfrak{a}})$ on
any object $M$ from $\mathscr{A}$ is diagonalizable.
\end{enumerate}
Since the functors $E\otimes _-$ and $E^*\otimes _-$ are both left and right
adjoint to each other on the category of all $\tilde{\mathfrak{a}}$-modules,
(L\ref{lll2}) implies that $E\otimes _-$ is in fact exact (as endo-functor of
$\mathscr{A}$). (L\ref{lll3}) guarantees that $\mathscr{A}$ is abelian, has
enough projectives (which are also injective) and that each object of
$\mathscr{A}$ has finite length (with respect to the abelian structure of
$\mathscr{A}$, but not as a $\tilde{\mathfrak{a}}$-module in general).

Given an admissible $\mathscr{A}$, we can construct a series of categories
$\mathcal{O}\{\mathfrak{p},\mathscr{A}\}$ as follows: We take a semisimple
(or reductive) Lie algebra $\g$ with a chosen Borel subalgebra $\mathfrak{b}$,
and require that $\mathfrak{p}\supset \mathfrak{b}$ is a parabolic subalgebra
of $\g$ such that $\mathfrak{p}=\tilde{\mathfrak{a}}\oplus
\mathfrak{n}_\mathfrak{p}$ is the Levi decomposition of $\mathfrak{p}$. Given
these data it makes sense to make the following definition:

\begin{definition}\label{defopl}
{\rm
The category
$\mathcal{O}\{\mathfrak{p},\mathscr{A}\}$ is
the full subcategory of the category of
$\mathfrak{g}$-modules given by all objects which are

\begin{enumerate}[(PL1)]
\item\label{OPLone} finitely generated,
\item\label{OPLtwo} locally $\mathfrak{n}_\mathfrak{p}$-finite,
\item\label{OPLthree} direct sums of objects from $\mathscr{A}$
when viewed as $\tilde{\mathfrak{a}}$-modules.
\end{enumerate}
}
\end{definition}

\subsection{Special cases of $\mathcal{O}\{\mathfrak{p},\mathscr{A}\}$}
\label{special}
\subsubsection*{Category $\mathcal{O}$}
If $\mathfrak{p}=\mathfrak{b}$ then $\tilde{\mathfrak{a}}=\mathfrak{h}$ is
abelian. Let $\mathscr{A}$ be the category of all finite dimensional semisimple
$\mathfrak{h}$-modules. This category is obviously admissible. The category
$\mathcal{O}\{\mathfrak{p},\mathscr{A}\}$ in this case is nothing else than the
usual category $\mathcal{O}=\mathcal{O}(\mathfrak{g},\mathfrak{b})$. Note that
the property (PL\ref{OPLthree}) in this case just means that the modules from
$\mathcal{O}\{\mathfrak{p},\mathscr{A}\}$ have a weight space decomposition.
The category $\mathcal{O}$ is a highest weight category with standard modules
given by the Verma modules.

\subsubsection*{The parabolic category $\cO^\p$}
If $\mathfrak{p}$ is any parabolic and $\mathscr{A}$ is the category of finite
dimensional semi-simple $\tilde{\mathfrak{a}}$-modules, then
$\mathcal{O}\{\mathfrak{p},\mathscr{A}\}$ is the parabolic category
$\mathcal{O}^\mathfrak{p}$. The category $\mathcal{O}^\p$ is a highest weight
category with standard modules given by the parabolic Verma modules.

\subsubsection*{The category $\mathcal{O}_0^{\mathfrak{p}-pres}$}\label{3}
Let $\mathfrak{p}$ be any parabolic subalgebra with Weyl group $W'$ and longest
element $w'_0$. Consider the corresponding indecomposable projective module
$P^{\mathfrak{a}}(w'_0\cdot0)$  in the category $\mathcal{O}(\mathfrak{a},
\mathfrak{a}\cap\mathfrak{b})$ corresponding to $\mathfrak{a}$. Let
$\mathscr{A}$ be the smallest abelian category which contains this
$P^\mathfrak{a}(w'_0\cdot0)$ and is closed under tensoring with
finite-dimensional simple $\mathfrak{a}$-modules and taking quotients. Extend
$\mathscr{A}$ to a category of $\tilde{\mathfrak{a}}$-modules by allowing
diagonalizable action of $\mathfrak{z}(\tilde{\mathfrak{a}})$.
Then the category $\mathscr{A}$ is
admissible and $\mathcal{O}\{\mathfrak{p},\mathscr{A}\}$ is the category of
modules which are presentable by the $P(w\cdot \lambda)\in\cO$, where $w$
runs through $(W'\backslash W)_{long}$ and $\lambda$ is an integral weight in
$\h^*_{dom}$ (for details see \cite{MS}). This category is also equivalent to
the category of Harish-Chandra bimodules with generalized trivial integral
central character  from the left hand side and the singular central character
given by $W'$ from the right hand side (for details see e.g. \cite{BG},
\cite[Kapitel 6]{Ja2}). This category $\mathcal{O}\{\mathfrak{p},\mathscr{A}\}$
is not a highest weight category in general, but still equivalent to a module
category over a so-called properly stratified algebra, see \cite{MS}.

\subsection{From highest weight categories to stratified algebras}
As usual, the category  $\mathcal{O}\{\mathfrak{p},\mathscr{A}\}$ decomposes
into a direct sum of full subcategories, each of which is equivalent to a
module category over a finite-dimensional associative algebra. Any block of
the (parabolic) category $\cO$ is a highest weight category, hence the
associated algebra is quasi-hereditary. In general, this is not true for a
block of $\mathcal{O}\{\mathfrak{p},\mathscr{A}\}$ (see for example \cite{MS}).
The algebras which appear from blocks of $\mathcal{O}\{\mathfrak{p},
\mathscr{A}\}$ are however always {\em weakly properly stratified} in the sense
of \cite[Section~2]{Fr}. The proof of this fact is completely analogous to the
properly stratified case, and we refer to \cite[Section 3]{FKM} for details.

A weakly properly stratified structure of an algebra means the following: the
isomorphism classes of simple modules are indexed by a partially pre-ordered
set $I$ and we have so-called {\em standard} and {\em proper standard} modules
(both indexed by $I$ again) such that projective modules have
standard filtrations, i.e. filtrations with subquotients isomorphic to standard
modules, and standard modules have proper standard filtrations. Which
subquotients are allowed to occur in the above filtrations and in the
Jordan-H{\"o}lder filtrations of proper standard modules is given by the
partial pre-order (for a precise definition we refer to \cite{Fr}).

The modules defining the stratified structure are given in terms of
parabolically induced modules as follows: If $V$  is any
$\tilde{\mathfrak{a}}$-module, we consider $V$ as a $\mathfrak{p}$-module
by letting $\mathfrak{n}_\mathfrak{p}$ act trivially and define the {\it
parabolically induced module}
\begin{displaymath}
\Delta(\mathfrak{p}, V):=U(\g)\otimes_{U(\mathfrak{p})}V.
\end{displaymath}
If $V$ is a simple object of $\mathscr{A}$ then
$\Delta(\mathfrak{p},V)$ is a {\em proper standard} module; if $V$
is projective then $\Delta(\mathfrak{p},V)$ is a {\em standard}
module. The dual construction (using conduction) gives rise to {\it
(proper) costandard module}. If $V$ is a simple
$\tilde{\mathfrak{a}}$-module, then $\Delta(\mathfrak{p},V)$ is
usually called a {\em generalized Verma module}, or simply a {\em
GVM}.

Let $\mathscr{F}(\Delta)$ denote the full subcategories of
$\mathcal{O}\{\mathfrak{p},\mathscr{A}\}$, given by all modules, which admit a
standard filtration, that is a filtration, whose subquotients are standard
modules. Analogously one defines $\mathscr{F}(\overline{\Delta})$ for modules
with proper standard filtration, $\mathscr{F}(\nabla)$ for modules with
costandard filtration, and $\mathscr{F}(\overline{\nabla})$ for modules with
proper costandard filtration. In this notation the property to be weakly
properly stratified is equivalent to the claim that all projective modules in
$\mathcal{O}\{\mathfrak{p},\mathscr{A}\}$ belong to both
$\mathscr{F}(\Delta)$ and $\mathscr{F}(\overline{\Delta})$.

We would like to point out that weakly properly stratified algebras form a
strictly bigger class than properly stratified algebras as the classes of
simple modules might be only partially pre-ordered, not partially ordered.
As a consequence there could be non-isomorphic standard modules $\Delta_1$
and $\Delta_2$ such that
$\Hom(\Delta_1,\Delta_2)\not=0\not=\Hom(\Delta_2,\Delta_1)$ (which
will be in fact the case in almost all the examples occurring from
now on in this paper).

\subsection{Parabolic induction via $\mathcal{O}\{\mathfrak{p},\mathscr{A}\}$}\label{s5.4}

Let us return to the case $W=S_n$ with some fixed parabolic subgroup
$W'=S_{i_1}\times S_{i_2}\times\cdots\times S_{i_r}$, where
$i_1+i_2+\dots +i_r=n$. Let $\mathds{H}$ and $\mathds{H}'$ be the
corresponding Hecke algebras. Assume we are given a right cell $\mathbf{R}'$ of
$W'$, then $\mathbf{R}'=\mathbf{R}'_{i_1}\times
\mathbf{R}'_{i_2}\times\cdots\times \mathbf{R}'_{i_r}$  for some right
cells $\mathbf{R}'_{i_j}$ in $S_{i_j}$.

Recall from Theorem~\ref{thm6} the categorification
$\mathscr{C}^{\mathbf{R}'_{i_j}}$ of the right cell module associated with
$\mathbf{R}'_{i_j}$. From Subsection~\ref{s5.1} we deduce that the outer
product,  call it $\mathscr{C}^{\mathbf{R}'}$, of these categories categorifies
the cell  module corresponding to $\mathbf{R}'$. The objects of
$\mathscr{C}^{\mathbf{R}'}$ are certain
$\tilde{\mathfrak{a}}:=\mathfrak{gl}_{i_1}\oplus
\mathfrak{gl}_{i_2}\oplus\cdots\oplus \mathfrak{gl}_{i_r}$-modules. Let
$\mathscr{P}$ denote the additive closure of the category of all modules,
which have the form $E\otimes P$, where $P\in \mathscr{C}^{\mathbf{R}'}$
is projective and $E$ is a simple finite-dimensional
$\tilde{\mathfrak{a}}$-module. Set $\mathscr{A}^{\mathbf{R}'}=
\overline{\mathscr{P}}$.

\begin{lemma}\label{lem51}
$\mathscr{A}^{\mathbf{R}'}$ is admissible.
\end{lemma}

\begin{proof}
As translations  are exact, condition (L\ref{lll2}) is satisfied by
definition. Condition (L\ref{lll4}) follows again from the
definitions as  the action of $\mathfrak{z}(\tilde{\mathfrak{a}})$
on any simple finite-dimensional $\tilde{\mathfrak{a}}$-module is
diagonalizable. It is left to check (L\ref{lll3}). By definition,
$\mathscr{A}^{\mathbf{R}'}$ is a subcategory of
$\mathcal{O}\{\tilde{\mathfrak{a}}, \tilde{\mathfrak{a}}
\cap\mathfrak{b})$. The block decomposition of the latter (with
respect to  the action of the center of $Z(\tilde{\mathfrak{a}})$)
induces a block decomposition of the former. Since  translations are
exact and send projectives to projectives,
$\mathscr{A}^{\mathbf{R}'}$ has enough projectives. These projective
modules are also injective by \eqref{cons1} from
Subsection~\ref{s4.2}. Therefore the condition (L\ref{lll3}) follows
from the definitions and \cite[Section~5]{Au}.
\end{proof}

By Lemma~\ref{lem51}, the category
$\mathcal{O}\{\mathfrak{p},\mathscr{A}^{\mathbf{R}'}\}$ is defined.
By construction, it is a subcategory of $\mathcal{O}$ and hence inherits a
decomposition from the block decomposition of $\mathcal{O}$ which we call
the block decomposition of
$\mathcal{O}\{\mathfrak{p},\mathscr{A}^{\mathbf{R}'}\}$.
Denote by $\mathcal{O}\{\mathfrak{p},\mathscr{A}^{\mathbf{R}'}\}_0$
the block of $\mathcal{O}\{\mathfrak{p},\mathscr{A}^{\mathbf{R}'}\}$
corresponding to the trivial central character. We note that one
can show that $\mathcal{O}\{\mathfrak{p},\mathscr{A}^{\mathbf{R}'}\}_0$ is
indecomposable by invoking Theorem~\ref{thm6}. We omit a proof,
since the result will not be relevant for the following.

From the definition of $\mathcal{O}\{\mathfrak{p},\mathscr{A}^{\mathbf{R}'}\}$
we have that simple objects in
$\mathcal{O}\{\mathfrak{p},\mathscr{A}^{\mathbf{R}'}\}_0$ are in a natural
bijection with the elements from $W$ of the form $xw$, where $w\in
(W'\backslash W)_{short}$ and $x\in \mathbf{R}'$. Denote by
$\Delta({\mathfrak{p}},xw)$ and $\overline{\Delta}(\mathfrak{p},xw)$
the standard, respectively proper standard, module in
$\mathcal{O}\{\mathfrak{p},\mathscr{A}^{\mathbf{R}'}\}_0$ corresponding to
$xw$. Dually we also have the (proper) costandard module
$\nabla({\mathfrak{p}},xw)$ ($\overline{\nabla}(\mathfrak{p},xw)$).
For details see \cite[Section~2]{Fr}. Finally, let $P(\mathfrak{p},xw)$ be
the projective cover of $\Delta(\mathfrak{p},xw)$ in
$\mathcal{O}\{\mathfrak{p},\mathscr{A}^{\mathbf{R}'}\}$.

Set
\begin{gather*}
\mathds{I}(\mathbf{R}')=\{(x,w)\,:\,x\in \mathbf{R}',
w\in (W'\backslash W)_{short}\},\\
\mathds{J}(\mathbf{R}')=\{y\in W\,:\, y\geq_{\mathsf{R}}\mathbf{R}',
y\neq xw \text{ for any }(x,w)\in \mathcal{I}(\mathbf{R}')\}.
\end{gather*}
In particular, by above the set $\mathds{I}(\mathbf{R}')$ indexes bijectively
the isomorphism classes of indecomposable projective modules in $\mathcal{O}\{\mathfrak{p},\mathscr{A}^{\mathbf{R}'}\}_0$.
From Subsection~\ref{s3.2} and the definitions it follows that for any
$(x,w)\in \mathds{I}(\mathbf{R}')$ the module $P(\mathfrak{p},xw)$ is the
quotient of $P(xw)$ modulo the trace of all $P(y)$,
$y\in \mathds{J}(\mathbf{R}')$. In particular, all
projectives in $\mathcal{O}\{\mathfrak{p},\mathscr{A}^{\mathbf{R}'}\}_0$
are gradable and hence the endomorphism ring $B$ of a minimal projective
generator of $\mathcal{O}\{\mathfrak{p},\mathscr{A}^{\mathbf{R}'}\}_0$
inherits a $\mZ$-grading from the ring $A$ from Subsection~\ref{s25.5}.
We denote by $\mathcal{O}\{\mathfrak{p},\mathscr{A}^{\mathbf{R}'}\}_0^\mZ$
the category of finite-dimensional graded $B$-modules.

Let $S(\mathbf{R}')\otimes_{\mathds H'}\mathds{H}$ be the induced cell module.
By definition, it has a $\mZ[v,v^{-1}]$-basis given by $\Delta_{x,w}:=
\underline{H}_x\otimes H_w$, where $(x,w)\in \mathds{I}(\mathbf{R}')$.
Hence we can define a homomorphism of $\mZ[v,v^{-1}]$-modules  as follows:
\begin{eqnarray}\label{Psi}
\Psi_{\mathbf{R}'}:S(\mathbf{R}')\otimes_{\mathds H'}\mathds{H}
&\overset{\sim}{\longrightarrow}&\left[\mathcal{O}\{\mathfrak{p},
\mathscr{A}^{\mathbf{R}'} \}_0^{\mathbb{Z}}\right]\\
\Delta_{x,w} &\mapsto &\left[\Delta({\mathfrak{p}},xw)\right].
\nonumber
\end{eqnarray}
Let $\mathscr{P}(\mathfrak{p},\mathscr{A}^{\mathbf{R}'})$
denote the additive category of all (graded) projective modules in
$\mathcal{O}\{\mathfrak{p},\mathscr{A}^{\mathbf{R}'}\}_0$.
We obtain the following main result:

\begin{theorem}[Categorification of induced cell modules]\label{thm53}
{\tiny .}

\begin{enumerate}[(i)]
\item\label{thm53.0}
The map $\Psi_{\mathbf{R}'}$ is a homomorphism of $\mathds{H}$-modules.
\item\label{thm53.1}
$\big(\mathcal{O}\{\mathfrak{p},\mathscr{A}^{\mathbf{R}'}\}_0^{\mathbb{Z}}\;,
\Psi_{\mathbf{R}'},\{\theta_s^l\}_{s\in S}\big)$ is
precategorification of the induced (right) cell
$\mathds{H}$-module $S(\mathbf{R}')\otimes_{\mathds H'}\mathds{H}$
whereas  $\big(\mathscr{P}(\mathfrak{p},\mathscr{A}^{\mathbf{R}'}),
\Psi_{\mathbf{R}'},\{\theta_s^l\}_{s\in S}\big)$ is
a categorification of this module with
respect to the  generators $\underline{H}_s$, $s\in S$.
\item\label{thm53.2}
The map $\Psi_{\mathbf{R}'}$ defines a  $\mZ((v))$-categorification
$\big(\mathcal{O}\{\mathfrak{p},\mathscr{A}^{\mathbf{R}'}\}_0^{\mathbb{Z}}\;,
\Psi_{\mathbf{R}'},\{\theta_s^l\}_{s\in S}\big)$ of the  induced cell $\mathds{H}^{\mZ((v))}$-module
$S(\mathbf{R}')^{\mZ((v))}\otimes_{(\mathds{H}')^{\mZ((v))}}
\mathds{H}^{\mZ((v))}$ with respect to
the generators $\underline{H}_s$, $s\in S$.
\end{enumerate}
\end{theorem}

\begin{proof}
In order to prove our theorem we only have to prove that
$\Psi_{\mathbf{R}'}$ is a homomorphism of $\mathds{H}$-modules. In other words,
we have to compare the combinatorics of the action of $\mathds{H}$ on
$S(\mathbf{R}')\otimes_{\mathds H'}\mathds{H}$ with the combinatorics of the
action of $\{\theta_s^l\}_{s\in S}$ on
$\mathcal{O}\{\mathfrak{p},\mathscr{A}^{\mathbf{R}'}\}_0^{\mathbb{Z}}$.
Fix $(x,w)\in \mathds{I}(\mathbf{R}')$  and $s\in S$.

If $ws\in (W'\backslash W)_{short}$ then the definition of
$\mathds{H}$ (see Subsection~\ref{s25.2}) gives
\begin{displaymath}
\Delta_{x,w}\underline{H}_s=
\begin{cases}
\Delta_{x,ws} +v \Delta_{x,w}, & \text{ if } ws>w;\\
\Delta_{x,ws} +v^{-1} \Delta_{x,w}, & \text{ if } ws<w.
\end{cases}
\end{displaymath}
If $ws\not\in (W'\backslash W)_{short}$ we have that $ws=s'w$ for $s'\in S\cap
W'$. In particular $ws>w$ and the definition of $\mathds{H}$ gives
$H_w\underline{H}_s=\underline{H}_{s'}H_w$. Therefore
\begin{equation}\label{eqthm53.1}
\Delta_{x,w}\underline{H}_s=
(\underline{H}_x\otimes H_w)\underline{H}_s=
\underline{H}_x\otimes\underline{H}_{s'}H_w=
\underline{H}_x\underline{H}_{s'}\otimes H_w,
\end{equation}
and $\underline{H}_x\underline{H}_{s'}$ can be computed using the
definition of $S(\mathbf{R}')$, i.e. \eqref{formula1}.

Now let us compare this with the combinatorics of the translation functors. Assume first that $ws=s'w\not\in (W'\backslash W)_{short}$. If $M$ is a
$\tilde{\mathfrak{a}}$-module and $E$ is a finite-dimensional
$\mathfrak{g}$-module (which we can also view as a finite-dimensional
$\tilde{\mathfrak{a}}$-module), then the Poincar{\'e}-Birkoff-Witt Theorem
implies the so-called tensor identity
$U(\mathfrak{g})\otimes_{U(\mathfrak{p})}(E\otimes
M)\cong E\otimes U(\mathfrak{g})\otimes_{U(\mathfrak{p})}M$ as $\tilde{\mathfrak{a}}$-modules (in both cases the $\tilde{\mathfrak{a}}$-module
structure is given by restriction). This implies that the computation of
$[\theta_s^l\Delta({\mathfrak{p}},xw)]$ reduces to the computation of
$[\theta_{s'}^l V]$, where $V$ is the indecomposable projective module in
$\mathscr{A}^{\mathbf{R}'}$ such that $\Delta(\mathfrak{p},xw)=
U(\mathfrak{g})\otimes_{U(\mathfrak{p})}V$. From Theorem~\ref{thm5}
it follows that the result is given by \eqref{formula1}. Hence it perfectly
fits with the computation of \eqref{eqthm53.1}. Finally, let us assume that
$ws\in (W'\backslash W)_{short}$. We have to show that in this case
\begin{equation}\label{eqthm53.2}
[\theta_s^l\Delta({\mathfrak{p}},xw)]=
\begin{cases}
[\Delta({\mathfrak{p}},xws)]+v[\Delta({\mathfrak{p}},xw)],& ws>w\\
[\Delta({\mathfrak{p}},xws)]+v^{-1}[\Delta({\mathfrak{p}},xw)],& ws<w.
\end{cases}
\end{equation}
Since all (proper) standard modules  are parabolically induced, from
our observation about the parabolic induction and projective
functors above it follows that projective functors preserve
both $\mathscr{F}(\Delta)$ and $\mathscr{F}(\overline{\Delta})$.
By duality, projective  functors also preserve both
$\mathscr{F}(\nabla)$ and $\mathscr{F}(\overline{\nabla})$. Hence
$\theta_s^l\Delta({\mathfrak{p}},xw)\in \mathscr{F}(\Delta)$ and we only
have to compute which standard modules occur in the standard filtration of
$\Delta({\mathfrak{p}},xw)$ and with which multiplicity.
From \cite[4.1]{Fr} it follows that the multiplicity of
$\Delta({\mathfrak{p}},y)\langle k\rangle$ in the standard filtration of
$\theta_s^l\Delta({\mathfrak{p}},xw)$ equals the dimension of
\begin{displaymath}
\mathrm{Hom}_{\mathcal{O}\{\mathfrak{p},\mathscr{A}^{\mathbf{R}'}
\}_0^{\mathbb{Z}}}(\theta_s^l\Delta({\mathfrak{p}},xw),
\overline{\nabla}({\mathfrak{p}},y)\langle k\rangle).
\end{displaymath}
Write $\theta_s^l=\theta_s^{out}\theta_s^{on}$, where
$\theta_s^{on}$ and $\theta_s^{out}$ are the graded translations onto and
out of the $s$-wall (see \cite[Corollary 8.3]{Stgrad}). Adjunction properties \cite[Theorem 8.4]{Stgrad} give
\begin{eqnarray*}
&&\mathrm{Hom}_{\mathcal{O}\{\mathfrak{p},\mathscr{A}^{\mathbf{R}'}
\}_0^{\mathbb{Z}}} (\theta_s^{out}\theta_s^{on}\Delta^{\mathfrak{p}}(xw),
\overline{\nabla}({\mathfrak{p}},y)\langle k\rangle)\\
&=&
\mathrm{Hom}_{\mathcal{O}\{\mathfrak{p},\mathscr{A}^{\mathbf{R}'}
\}_0^{\mathbb{Z}}} (\theta_s^{on}\Delta^{\mathfrak{p}}(xw),
\theta_s^{on}\overline{\nabla}({\mathfrak{p}},y)\langle k+1\rangle).
\end{eqnarray*}
A character argument shows that $\theta_s^{on}\Delta({\mathfrak{p}},xw)$ is a graded lift of a standard module
and $\theta_s^{on}\overline{\nabla}({\mathfrak{p}},y)$ is a graded lift of a proper standard module on the wall, and a direct calculation (using \cite[Theorem 8.1]{Stgrad}) gives
\begin{displaymath}
\mathrm{Hom}_{\mathcal{O}\{\mathfrak{p},\mathscr{A}^{\mathbf{R}'}
\}_0^{\mathbb{Z}}} (\theta_s^{on}\Delta({\mathfrak{p}},xw),
\theta_s^{on}\overline{\nabla}({\mathfrak{p}},y)\langle k+1\rangle)=
\begin{cases}
\mathbb{C},& y=xws, k=0;\\
\mathbb{C},& y=xw, k=1, ws>w;\\
\mathbb{C},& y=xw, k=-1, ws<w;\\
0,&\text{ otherwise. }
\end{cases}
\end{displaymath}
Formula \eqref{eqthm53.2} follows and the proof is complete.
\end{proof}

\subsection{Uniqueness of categorification}

Assume that we are still in the situation of Subsection~\ref{s5.4}.

\begin{proposition}\label{cunique}
Let $\mathbf{R}'_1$ and $\mathbf{R}'_2$ by two right cells of $W'$
inside the same two-sided cell. Then the categories
$\mathcal{O}\{\mathfrak{p},\mathscr{A}^{\mathbf{R}_1'}\}$ and
$\mathcal{O}\{\mathfrak{p},\mathscr{A}^{\mathbf{R}'_2}\}$ are
equivalent.
\end{proposition}

\begin{proof}
The equivalence between $\mathscr{A}^{\mathbf{R}'}$ and
$\mathscr{A}^{\mathbf{R}'_2}$, constructed in Theorem~\ref{thm6}
extends in a straightforward way to an equivalence between the
categories $\mathcal{O}\{\mathfrak{p},\mathscr{A}^{\mathbf{R}'}\}$
and $\mathcal{O}\{\mathfrak{p},\mathscr{A}^{\mathbf{R}'_2}\}$.
\end{proof}

\section{Combinatorics and filtrations of induced cell modules}

In this section we first introduce a non-degenerate bilinear form on
the induced cell modules and establish a categorical version of it.
As a result we get four different distinguished bases in any induced
cell modules which we then will interpret via four distinguished
classes of objects in the corresponding categorification. Afterwards
we describe the resulting refined Kazhdan-Lusztig combinatorics and
also introduce a natural filtration on induced cell modules which
are induced from a natural counterpart on their categorifications.

\subsection{Different bases and the combinatorics of induced cell modules}

Assume that we are still in the situation of Subsection~\ref{s5.4}. Any
cell module $S(\mathbf{R}')$ has a unique up to a scalar non-degenerate
symmetric bilinear $\mathds{H}'$-invariant form $\langle\cdot,\cdot\rangle$.
We normalize this form such that its categorification is given by
Proposition~\ref{bilinearform}. We first state the following easy lemma:

\begin{lemma}
The induced module $S(\mathbf{R}')\otimes_{\mathds H'}\mathds{H}$
has a non-degenerate symmetric  $\mathds{H}$-invariant bilinear form
$(\cdot, \cdot)$ with values in $\mZ[v,v^{-1}]$ given by
\begin{eqnarray*}
(m\otimes H_x, n\otimes H_y)=\delta_{x,y}\langle m,n\rangle,
\end{eqnarray*}
for any $x$, $y\in (W'\backslash W)_{short}$ and $m$, $n\in S(\mathbf{R}')$.
\end{lemma}

\begin{proof}
The form is obviously symmetric and non-degenerate, since so is $\langle\cdot, \cdot\rangle$. It is left to show the $\mathds{H}$-invariance. Let $s\in
S\subset W$ and $m, n,x,y$ as above.  For the rest of the proof it would
be convenient to set $X=(m\otimes H_xH_s, n\otimes H_y)$ and
$Y=(m\otimes H_x, n\otimes H_yH_s)$.

Assume first that $xs, ys\in (W'\backslash W)_{short}$. If $xs>x$ and $ys>y$
then $xs\neq y$, $x\neq ys$ and hence  $X=(m\otimes H_{xs}, n\otimes
H_y)=\delta_{xs,y}\langle m,n\rangle=0$, and $Y=(m\otimes H_x, n\otimes
H_{ys})=\delta_{x,ys}\langle m,n\rangle=0$.  If $xs>x$ and $ys<y$ then
$X=(m\otimes H_{xs}, n\otimes H_y)=\delta_{xs,y}\langle m,n\rangle$, and
$Y=(m\otimes H_x, n\otimes H_{ys}+(v^{-1}-v)H_{y})=\delta_{x,ys}\langle
m,n\rangle=\delta_{xs,y}\langle m,n\rangle$ (as $x\neq y$, and $xs=y$ if
and only if $x=ys$). If $xs<x$ and $ys>y$ then the argument is analogous (by
symmetry). If $xs<x$ and $ys<y$ then $X=(m\otimes H_{xs}+(v^{-1}-v)H_{x},
n\otimes H_y)=(v^{-1}-v)\delta_{x,y}\langle m,n\rangle$, and $Y=(m\otimes H_x,
n\otimes H_{ys}+(v^{-1}-v)H_{y})=(v^{-1}-v)\delta_{x,y}\langle m,n\rangle$.

Now let us assume $xs\notin (W'\backslash W)_{short}$, $ys\in (W'\backslash
W)_{short}$. We write $xs=s'x$ and get $X=(m H_{s'}\otimes H_{x}, n\otimes
H_y)=\delta_{x,y}\langle m,n\rangle=0$. On the other hand, $Y=(m\otimes H_x,
n\otimes H_{y}H_s)=0$ as $x\not\in{y,ys}$.

Finally let us assume $xs\notin (W'\backslash W)_{short}$,  $ys\notin
(W'\backslash W)_{short}$. We write $xs=s'x$ and $ys=ty$, where $s', t\in W'$
are simple reflections. Then $X=(mH_{s'}\otimes H_x, n\otimes H_y)\not=0$
implies $x=y$, and then also  $s'=t$. The same holds if $Y=(m\otimes H_x,
nH_t\otimes H_y)\not=0$; and both terms have the same value, namely $\langle
mH_{t},n\rangle=\langle m,nH_{t}\rangle$, since $\langle\cdot,\cdot\rangle$ is
$\mathds{H}'$-invariant.  The statement of the lemma follows.
\end{proof}

The involution $h'\mapsto \overline{h'}$ on $\mathds{H}'$ restricts
to an involution on any right cell module and is on the other hand
itself the restriction of the involution $h\mapsto\overline{h}$ on
$\mathds{H}$. Therefore, we get an involution
\begin{eqnarray*}
S(\mathbf{R}')\otimes_{\mathds H'}\mathds{H}&\rightarrow&
S(\mathbf{R}')\otimes_{\mathds H'}\mathds{H} \\
m\otimes h&\mapsto&\overline{m\otimes
h}:=\overline{m}\otimes\overline{h}.
\end{eqnarray*}
For $(x,w)\in\mathds{I}(\mathbf{R}')$ we define the
{\it Kazhdan-Lusztig element} $\underline{H}_x\boxdot{H}_w$ as the unique self-dual element in $S(\mathbf{R}')\otimes_{\mathds H'}\mathds{H}$,
satisfying
\begin{displaymath}
\underline{H}_x\boxdot{H}_w\in \underline{H}_x\otimes {H}_w +
\sum_{(x',w')\in\mathds{I}(\mathbf{R}')}v\mathbb{Z}[v]
\underline{H}_{x'}\otimes {H}_{w'}.
\end{displaymath}
The existence and uniqueness of such elements is obtained by
standard arguments (see e.g. \cite[Theorem~2.1]{SoKipp}).

The induced module $S(\mathbf{R}')\otimes_{\mathds H'}\mathds{H}$ has then
four distinguished bases:
\begin{itemize}
\item the {\it Kazhdan-Lusztig basis} (or short {\it KL basis}) given by the
Kazhdan-Lusztig elements ${\underline{H}_x\boxdot{H}_w}$, where $x\in
\mathbf{R}'$ and $w\in (W'\backslash W)_{short}$.
\item the {\it Kazhdan-Lusztig-standard basis} (or short {\it KL-s basis}) given by the elements $\Delta_{x,w}=\underline{H}_x\otimes H_w$, where $x\in
\mathbf{R}'$ and $w\in (W'\backslash W)_{short}$.
\item the {\it dual Kazhdan-Lusztig basis} (or short {\it dual KL basis}) which
is the dual of the KL-basis with respect to the form $(\cdot, \cdot)$.
\item the {\it dual Kazhdan-Lusztig-standard basis} (or short
{\it dual KL-s basis}).
\end{itemize}
These bases have the following categorical interpretation:

\begin{theorem}[Combinatorics]\label{combinatorics}
The isomorphism $\Psi_{\mathbf{R}'}$ from \eqref{Psi} defines the
following correspondences:

\begin{tabular}[thb]{ccc}
\\ \text{KL-s  basis}&$\leftrightarrow$&
\text{isoclasses of standard lifts of standard modules}\\
\\ \text{KL basis}&$\leftrightarrow$&
\text{isoclasses of standard lifts of indecomposable projectives}\\
\\ \text{dual KL-s basis}&$\leftrightarrow$&
\text{isoclasses of standard lifts of proper standard  modules}\\
\\ \text{dual KL basis}&$\leftrightarrow$&
\text{isoclasses of standard lifts of simple modules}. \\
\end{tabular}
\end{theorem}

\begin{proof}
Let $(x,w)\in\mathds{I}(\mathbf{R}')$. The isomorphism class
$\left[\Delta({\mathfrak{p}},xw)\right]$ is mapped to $\Delta_{x,w}$
by definition, hence the first statement of the theorem is obvious.
Note that for $w=e$, the module $\Delta({\mathfrak{p}},x)$ is always
projective and $\Delta_{x,e}=\underline{H}_x\boxdot{H}_e=
\underline{H}_x\otimes {H}_e$. This provides the starting point
for an induction argument which proves the remaining part of the theorem.

Before we do the induction argument we have to recall a few facts. First
recall that for $s\in S$ the functor $\theta_s^l$ sends projectives to
projectives, since it is left adjoint  to an exact functor. The usual weight
argument also shows that if $ws\in (W'\backslash W)_{short}$ and $ws>w$ then
\begin{eqnarray}
\label{eq:dec}
\theta_s^lP(\mathfrak{p},xw)\cong
P(\mathfrak{p},xws)\bigoplus_{(y,z)\not=(x,w)} a_{y,z}P(\mathfrak{p},yz),
\end{eqnarray}
at least if we forget the grading. Since the category
$\mathcal{O}\{\mathfrak{p},\mathscr{A}^{\mathbf{R}'}\}$ is by definition a
subcategory of $\mathcal{O}$, we could take the projective cover
$P(xw)\in\mathcal{O}_0^{\mathbb{Z}}$ of $P(\mathfrak{p},xw)\in
\mathcal{O}\{\mathfrak{p},\mathscr{A}^{\mathbf{R}'}\}_0^{\mathbb{Z}}$,
and the decomposition \eqref{eq:dec} is controlled by that of
$\theta_s^lP(xw)$. In particular, it is of the form as stated in
\eqref{eq:dec}(even as graded modules).

Assume now that the statement of the theorem is true for some $(x,w)\in
\mathds{I}(\mathbf{R}')$ and let $s\in W$ such that $ws>w$ and  $ws\in
(W'\backslash W)_{short}$. From Theorem~\ref{thm53} we know that
$\theta_s^lP(\mathfrak{p},xw)$ corresponds to  $H:=(\underline{H}_x\boxdot{H}_w)
\underline{H}_s$ under $\Psi_{\mathbf{R}'}^n$.  In particular
$H=\underline H_{x}\otimes {H}_{ws}+\sum_{(x',w')\not=
(x,ws)}\beta_{(x',w')}(v)\underline{H}_{x'}\otimes {H}_{w'}$, where
$\beta_{x',w'}(v)\in\mZ[v]$. From \eqref{eq:dec} and the explanation
afterwards we get then that  $P(\mathfrak{p},xw)$ corresponds to
\begin{displaymath}
H':=(\underline H_x\boxdot H_w)\underline{H}_s-\beta_{x',w'}(0)
\underline{H}_{x'}\otimes H _{w'}.
\end{displaymath}

Note that $H'$ can be characterized as the unique self-dual element with the
property that $H'\in \underline{H}_x\otimes {H}_{ws}+\sum_{(x',w')}v\mZ[v]
\underline{H}_{x'}\otimes{H}_{w'}$. The same characterization holds
for the  element $\underline H_x\boxdot{H}_{ws}$. Hence
$\underline{H}_x\boxdot{H}_{ws}$ is mapped to
$\left[P(\mathfrak{p},xws)\right]$ under the isomorphism
$\Psi_{\mathbf{R}'}$ and the second part of the theorem follows.

It is not difficult to verify that the bilinear form $(\cdot,\cdot)$
has again a categorical version as in
Proposition~\ref{bilinearform}. In particular, the isomorphism
classes of simple modules are dual to the ones of indecomposable
projective modules. Finally, the proper standard modules form a dual
basis to the costandard modules thanks to the duality on
$\mathcal{O}\{\mathfrak{p}, \mathscr{A}^{\mathbf{R}'}\}$ and the
usual Ext-orthogonality between standard and proper costandard
module (\cite[Theorem 3]{Fr}). The theorem follows.
\end{proof}

\begin{example}
{\rm
Let $W=S_3=<s,t>$ and $W'=<s>\cong S_2\times S_1$ and choose the right cell
$\mathbf{R}'$ of $W'$ corresponding to the (longest) element $s$. Then
$(W'\backslash W)_{short}=\{e,t,ts\}$. The categorification
$\mathscr{C}^{\mathbf{R}'}$ is then equivalent to the category of graded
$R=\mC[x]/(x^{2})=\End_{\mathfrak{gl}(2)}(\mathtt{P}(s\cdot0))$ as in
Example~\ref{ex6}, and $\mathcal{O}\{\mathfrak{p},\mathscr{A}^{\mathbf{R}'}\}
\cong\cO_0^{\p-pres}$ from Subsection~\ref{special}. The module
$\Delta(\mathfrak{p},se)$ is projective, hence $\Delta(\mathfrak{p},se)=
P(\mathfrak{p},se)$. A direct calculation shows that the projective module
$P(\mathfrak{p},st)$ has a standard-filtration of length two,  with
$\Delta(\mathfrak{p},st)$ as a quotient, and $\Delta(\mathfrak{p},se)$ as
a submodule; whereas $P(\mathfrak{p},sts)$ has a standard filtration with
$\Delta(\mathfrak{p},sts)$ occurring as a quotient, $\Delta(\mathfrak{p},st)$
as a subquotient, and $\Delta(\mathfrak{p},se)$ as a submodule (see the
detailed example in \cite[Section 9]{MS}). On the combinatorial side, the
standard basis element $\underline{H}_s\otimes H_e$ is a self-dual KL-basis
element. The element $\underline{H}_s\otimes H_t+v\underline{H}_s\otimes H_e$
is a KL-basis element. Now, $\underline{H}_s\otimes \underline{H}_{ts}=
\underline{H}_s\otimes(H_{ts}+v(H_t+H_s)+v^2H_e)$ is self-dual and equal to
$\underline{H}_s\otimes H_{ts}+v \underline{H}_s\otimes H_t+ \underline{H}_s
\otimes H_e+v^2\underline{H}_s\otimes H_e.$ Hence subtracting
$\underline{H}_s\otimes H_e$ gives $\underline{H}_s\boxdot
H_{ts}=\underline{H}_s\otimes H_{ts}+v\underline{H}_s\otimes
H_{t}+v^2\underline{H}_s\otimes H_{e}.$
}
\end{example}

\subsection{Stratifications of induced modules}\label{s5.5}

Let us come back to the examples in Subsection~\ref{s5.2} and assume $W=S_n$
with parabolic subgroup $W'$. Let $\mu$ be the composition of $n$ which defines
$W'$ and let $\la$ be the corresponding partition. Consider again the
permutation module  $\mathcal{M}=\mathcal{M}^\la$ and the irreducible cell
module $S(\la)$, which specializes to the irreducible Specht module $S^\la$
corresponding to $\la$. This is naturally a submodule of $\mathcal{M}^\la$.
Over the complex numbers, however, $\mathcal{M}^\la$ is completely reducible
and contains $S(\la)^{\mathbb{C}}$ as a unique direct summand.
Furthermore, over the  complex numbers, any  finite dimensional (right)
$\mathds{H}^{\mathbb{C}}$-module $M$ has a decomposition into isotypic
components. This special feature is however not independent of the
ground field (as Specht module are only indecomposable but not irreducible in
general), in particular it is not an integral phenomenon. However, there
is a natural filtration  of $\mathcal{M}$ by Specht modules which always exists
(see e.g. \cite[4.10 Corollary]{Mathas}).  The purpose of this subsection  is
to give a very natural categorical construction of a somewhat rougher
filtration on all induced cell modules. The idea is to use the notion of
Gelfand-Kirillov-dimension.

Consider the category
$\mathcal{O}\{\mathfrak{p},\mathscr{A}^{\mathbf{R}'}\}_0^{\mathbb{Z}}$.
The objects of this category are certain $\mathfrak{gl}_n$-modules. Any such
module $M$ has a well-defined {\it Gelfand-Kirillov-dimension}
$\operatorname{GKdim}(M)$. Recall the following easy facts:

\begin{lemma}
\label{GK}
\begin{enumerate}
\item For any $s\in S\subset W$ we have
$\operatorname{GKdim}(\theta_s^lM)\leq\operatorname{GKdim}(M)$.
\item
$\operatorname{GKdim}(M)=\operatorname{max}\{\operatorname{GKdim}(L_j)\}$,
where $L_j$ runs through the composition factors of $M$.
\end{enumerate}
\end{lemma}

\begin{proof}
See for example \cite[Lemmas 8.6, 8.8 and 8.7(1)]{Ja2}.
\end{proof}

For any positive integer $j$ we define
$(\mathcal{O}\{\mathfrak{p},\mathscr{A}^{\mathbf{R}'}\}_0^{\mathbb{Z}})_{\leq
j}$  to be the full subcategory of
$\mathcal{O}\{\mathfrak{p},\mathscr{A}^{\mathbf{R}'}\}_0^{\mathbb{Z}}$
consisting of all modules which have Gelfand-Kirillov dimension at
most $j$. From the Lemma above it follows that this subcategory is
closed under taking submodules, quotients and extensions, and also
stable under translations through walls. Therefore, we have a
filtration of the $\mathds{H}$-module $\left[\mathcal{O}\{\mathfrak{p},
\mathscr{A}^{\mathbf{R}'}\}_0^{\mathbb{Z}}\right]$.
For simplicity we relabel the filtration such that we have:
\begin{multline*}\label{gkdime1}
\{0\}\subsetneq \left[(\mathcal{O}\{\mathfrak{p},\mathscr{A}^{\mathbf{R}'}
\}_0^{\mathbb{Z}})_1\right] \subsetneq \left[(\mathcal{O}\{\mathfrak{p},
\mathscr{A}^{\mathbf{R}'}\}_0^{\mathbb{Z}})_2\right] \subsetneq \cdots\\
\cdots\subsetneq\left[(\mathcal{O}\{\mathfrak{p},\mathscr{A}^{\mathbf{R}'}
\}_0^{\mathbb{Z}})_r\right]= \left[\mathcal{O}\{\mathfrak{p},
\mathscr{A}^{\mathbf{R}'} \}_0^{\mathbb{Z}}\right]
\end{multline*}

The set of partitions of $n$ is ordered via the so-called dominance
ordering which we denote by $\unrhd$. Given two partitions
$\nu=\nu_1\geq \nu_2\ldots$ and $\mu=\mu_1\geq \mu_2\ldots$ we have
$\nu\unrhd\mu$ if and only if
$\sum_{j=1}^i\nu_j\geq\sum_{j=1}^i\mu_j$ for any $i\geq 1$. The
simple composition factors of the module $\mathcal{M}^\la$ are all
of the form $S(\mu)$, where $\la\unlhd\mu$ (see e.g. \cite[4.10,
Exercise 1]{Mathas} or \cite[Corollary~2.4.7]{Sa}).

The following result is the technical formulation of a fact which is
quite easy to describe: For every induced cell module
$\mathds{H}$-module $S(\mathbf{R}')\otimes_{\mathds H'}\mathds{H}$
we have a corresponding categorification, hence an attached category
$\cC$, of modules over some Lie algebra. The Gelfand-Kirillov
dimension induces a filtration on $\cC$ that corresponds to a
filtration of $S(\mathbf{R}')\otimes_{\mathds H'}\mathds{H}$ which
is an analogue of the Specht filtration of the induced cell module
given by the dominance ordering. More precisely we have the
following:

\begin{theorem}\label{GKdim}
Assume that we are in the setup of Subsection~\ref{s5.4}. For $i\geq 0$  set
\begin{displaymath}
Q_i=\{v\in S(\mathbf{R}')\otimes_{\mathds{H}'}\mathds{H}\,:\,
\Psi_{\mathbf{R}'}(v)\in \left[(\mathcal{O}\{\mathfrak{p},
\mathscr{A}^{\mathbf{R}'}\}_0^{\mathbb{Z}})_i\right]\}.
\end{displaymath}
Then we have:
\begin{enumerate}[(i)]
\item \label{GKdim.1}
$Q_0\subsetneq Q_1\subsetneq\dots\subsetneq Q_r=
S(\mathbf{R}')\otimes_{\mathds{H}'}\mathds{H}$ is a filtration of the induced
cell module $S(\mathbf{R}')\otimes_{\mathds{H}'}\mathds{H}$.
\item \label{GKdim.2}
Assume, $S(\la)$ occurs in the $i$-th and $S(\mu)$ in the $j$-th
filtration in \eqref{GKdim.1} respectively. Then $\la\lhd\mu$
implies $i<j$ (in other words: if $\la\lhd\mu$ then $S(\la)$ occurs
earlier than $S(\mu)$ as a subquotient of \eqref{GKdim.1}).
\item \label{GKdim.3}
All subquotients of the filtration \eqref{GKdim.1} are direct sums of Specht
modules.
\item \label{GKdim.4}
In the permutation module $\mathcal{M}^{\la}$ the Specht submodule
$S(\la)$ coincides with $Q_1$ (i.e. it is given by the subcategory
of modules of the minimal possible Gelfand-Kirillov-dimension from
$\mathcal{O}\{\mathfrak{p},\mathscr{A}^{\mathbf{R}'}\}_0^{\mathbb{Z}}$).
\end{enumerate}
\end{theorem}

\begin{proof}
The statement \eqref{GKdim.1} follows from Theorem~\ref{thm53} and
definitions.
To prove \eqref{GKdim.2} recall that there is Joseph's explicit formula
(see e.g. \cite[10.11 (2)]{Ja2})
\begin{equation}\label{gkdimf1}
2\operatorname{GKdim} (L(w\cdot0))=n(n-1)-\sum_{i}\mu_i(\mu_i-1),
\end{equation}
where $\mu$ is the shape of the tableaux associated with $w\in S_n$ via the Robinson-Schensted correspondence (in particular, simple modules in the same right cell have the same Gelfand-Kirillov dimension). Hence the statement
\eqref{GKdim.2} follows from Lemma~\ref{question} below, since for two
partitions $\mu$ and $\nu$ of $n$ we have $\sum_{i}\mu_i(\mu_i-1)<
\sum_{i}\nu_i(\nu_i-1)$ if and only  $\sum_{i}\mu_i^2<\sum_{i}\nu_i^2$.

\begin{lemma}\label{question}
Let $\mu$ and $\nu$ be partitions of $n$ and $l$ be the maximum of the
lengths of the partition. Then $\mu\lhd\nu$ implies
$\sum_{i=1}^n\mu_i^2<\sum_{i=1}^n\nu_i^2$.
\end{lemma}

\begin{proof}
If $l=2$ then $2(\mu_1^2+\mu_2^2)=(\mu_1+\mu_2)^2+(\mu_1-\mu_2)^2<
(\nu_1+\nu_2)^2+(\nu_1-\nu_2)=2(\nu_1^2+\nu_2^2)$. We will do induction on $l$.
Without loss of generality assume $\mu_i\not=\nu_i$ for $1\leq i\leq l$. Choose
now $i$ minimal such that $\mu_i<\nu_i$, but  $\mu_{i+1}>\nu_{i+1}$ and set
$m:=\min\{\mu_i-\nu_i, \nu_{i+1}-\mu_{i+1}\}$. It is easy to check that we get
a new partition $\sigma$, where $\sigma_i=\mu_i-m$, $\sigma_{i+1}=\mu_{i+i}+m$
and $\sigma_j=\mu_j$ for all other $j$. Note that $\sigma_k=\nu_k$ for some
$k\in\{i,i+1\}$. So we may apply the induction hypothesis to the partitions
$\sigma$ and $\nu$ with the common part removed. On the other hand
$(\mu_i,\mu_{i+1})\lhd(\sigma_{i}, \sigma_{i+1})$ satisfies the assumption of
the lemma, hence $\mu_i^2+\mu_{i+1}^2<\sigma_i^2+\sigma_{i+1}^2$ and so
$\sum_{j=1}^l\mu_j^2<\sum_{j=1}^l\sigma_j^2<\sum_{j=1}^l\nu_j^2$.
\end{proof}

From \eqref{GKdim.2} and \cite[Theorem~5.1]{Ge} it follows that the
indexing partitions of the Specht modules occurring in a fixed
subquotient of the filtration from \eqref{GKdim.1} are not
comparable in the right order. This implies \eqref{GKdim.3}.  The
claim \eqref{GKdim.4} follows immediately from \eqref{GKdim.2} and
\cite[Corollary~2.4.7]{Sa} (see the remark before the formulation of
the theorem).
\end{proof}

\begin{remark}\label{gkremark}
{\rm Using \cite[Corollary~2.4.7]{Sa} and \cite[4.10
Corollary]{Mathas} one can construct the following natural integral
filtration of the permutation module $\mathcal{M}^{\lambda}$: For
the first step of the filtration we take the submodule $S(\la)$
(note again that $\la$ is the minimal partition, with respect to the
dominance order, amongst the partitions indexing the subquotients of
$\mathcal{M}^{\lambda}$). To construct the second step in the
quotient we take the direct sum of all Specht modules, whose
partitions are minimal elements in the dominance order among all
other partitions which occur; and so on. For $n\leq 6$ the
constructed {\em dominance order filtration} will coincide with the
one given by Theorem~\ref{GKdim}\eqref{GKdim.1}. However, already
for $n=7$ one gets that the filtration given by
Theorem~\ref{GKdim}\eqref{GKdim.1} is a proper refinement of the
dominance order filtration. For example if we take $n=7$ and the
permutation module corresponding to the sign representation (this
permutation module is isomorphic  to the regular representation of
the Hecke algebra), it turns out that  the `dominance order
filtration' contains a step in which the subquotients  are Specht
modules corresponding to partitions $(5,1,1)$ and $(4,3)$. However,
$5^2+1^2+1^2=27\neq 25=4^2+3^2$ and hence these Specht  modules
occur in different layers of the filtration given by
Theorem~\ref{GKdim}\eqref{GKdim.1} because of \eqref{gkdimf1}. }
\end{remark}

\section{An alternative categorification of the permutation module}\label{s6}

In this section we propose an alternative categorification of the
permutation parabolic modules. The connection to the
categorification from Subsection~\ref{permmodule} is not completely
obvious (but can be made precise using \cite[6.4-6.5]{MOS}).

Let $W'$ be a subgroup of $W$ and let $\la\in\h^*_{dom}$ be an
integral weight with stabilizer $W'$ with respect to the dot-action.
The isomorphism classes of the Verma modules in $\mathcal{O}_\la$
are exactly given by the $M(x\cdot\la)$, where $x\in
(W/W')_{short}$.

For any simple reflection $s\in S$, the {\it twisting functor}
$T_s:\cO\rightarrow \cO$ (see Subsection~\ref{s25.5}) preserves
blocks, in particular induces $T_s:\cO_\la\rightarrow{}\cO_\la$. The
most convenient description (for our purposes) of these functors is
given in \cite{KM} in terms of partial coapproximation: Let
$M\in{}\cO_\la$ be projective. Let $M'\subset M$ be the smallest
submodule such that $M/M'$ has only composition factors of the form
$L(x\cdot\la)$, where $sx>x$. Then $M\mapsto M'$ defines a functor
$T_s$ from the additive category of projective modules in
${}\cO_\la$ to ${}\cO_\la$. This functor extends in a unique way to
a right exact functor $T_s:\cO_\la\rightarrow{}\cO_\la$ (for details
see \cite{KM}). From this definition of $T_s$ it is immediately
clear that this functor is gradable. More precisely, we have the
following:

\begin{lemma}\label{twist}{\rm (\cite[Proposition 5.1]{FKS})}
For any simple reflection $s\in W$ and integral weight $\la\in\h^*_{dom}$,
the twisting functor $T_s:{}\cO_\la\rightarrow{}\cO_\la$ is gradable. A graded
lift is unique up to isomorphism and shift in the grading.
\end{lemma}

\begin{proposition}\label{SchurWeyl}
Let $s\in S$.
\begin{enumerate}
\item The twisting functor $T_s$ is right exact, and exact when restricted the subcategory $\mathcal{V}_\la$ of $\cO_\la$ of modules having a filtration with subquotients isomorphic to Verma modules.
\item One can choose a graded lift ${\bf T}_s$ satisfying the following
properties:
\begin{align} \label{Rrel}
\begin{split}
&\big[{\bf T}_s M(x\cdot\la)\big]\\
&=\begin{cases}
[( M(sx\cdot\la)]+(v^{-1}-v)[( M(x\cdot\la)]&
\text{ if $sx<x$, $sx\in W/W'_{short}$},\\
[( M(sx\cdot\la))]&\text{ if $sx>x$, $sx\in (W/W')_{short}$},\\
v^{-1}[( M(x\cdot\la))]&\text{ if $sx\notin (W/W')_{short}$.}
\end{cases}
\end{split}
\end{align}
\item There is an isomorphism of (left) $\mZ[W]$-modules
\begin{eqnarray*}
\Psi_\la:\mZ[W]\otimes_{\mZ[W']}\mZ&\longrightarrow& [\mathcal{D}^b(\cO_\la)]\\
x\otimes 1&\longmapsto&\left[ M(x\cdot\la)\right],
\end{eqnarray*}
where the $\mZ[W']$-structure on $\mZ$ is trivial, and the
$\mZ[W]$-structure on the right hand side is induced by the action
of the left derived twisting functors $\mathcal{L} T_s$.
\end{enumerate}
\end{proposition}

\begin{proof}
The first statement follows directly from  \cite[Lemma 2.1]{AS}.
If we forget the grading (and put $v=1$) then the second statement follows
directly from \cite[Theorem 6.2, Definition 5.1 (ii)]{AL} and implies the
last statement. For the graded setup we refer to the proof of
\cite[Proposition 5.2]{FKS}.
\end{proof}

\section{Remarks on Schur-Weyl dualities}\label{s7}

For completeness we would like to formulate here a categorical version of the Schur-Weyl duality generalizing the approach of \cite{FKS}. Complete proofs and also a geometric interpretation in terms of generalized Steinberg varieties will appear in \cite{SS}.

\subsection{For permutation parabolic modules}\label{s7.1}

We assume again the setup of Subsection~\ref{s5.2}. Let $\la, \mu\in\h^*_{dom}$
be integral. If $F:\cO_\la\rightarrow\cO_\mu$ is a projective functor then it
induces a homomorphism $F^G:[\cO_\la]\rightarrow [\cO_\mu]$. Since finite
direct sums of projective functors are again projective functors, they form a
monoid. On the other hand, the composition of two projective functors (if
defined) is again a projective functor. The same holds if we work in the graded
setup with graded translation functors between the graded versions
$\cO_\la^{\mathbb{Z}}$ and $\cO_\mu^{\mathbb{Z}}$ of $\cO_\la$ and $\cO_\mu$
(see \cite{BGS}). This means we have the additive category of (graded)
projective functors from $ \cO_\la^{\mathbb{Z}}$ to $\cO_\mu^{\mathbb{Z}}$ with
its complexified split Grothendieck group $\left[\text{ projective functors:
}\cO_\la\rightarrow\cO_\mu\right]_{\oplus}^{\mC}$.

\begin{theorem}
With the notation from Subsection~\ref{s5.2} we have the following:
There is an isomorphism of $\mC[v,v^{-1}]$-modules
\begin{eqnarray*}
\Psi_{\la,\mu}:\quad\left[\text{ projective functors: }\cO_\la\rightarrow\cO_\mu\right]_{\oplus}^{\mC}&\cong&
\mathrm{Hom}_{\mathds{H}^{\mC}}(\mathcal{M}^\la,\mathcal{M}^\mu)\\
F&\mapsto&\Psi_\mu^{-1} F^G\Psi_\la.
\end{eqnarray*}
\end{theorem}

The latter result is true for any reductive complex Lie algebra $\mathfrak{g}$.
In the following we assume however $\mathfrak{g}=\mathfrak{sl}_n$. For any
Young subgroup $S_\la$ of $S_n$ we {\bf pick} some integral weight
$\la\in\h^*_{dom}$ where $W_\la\cong S_\la$.  Let $\Lambda$ be the set of all
these $\lambda$'s. For any positive integer $d$ let $\Lambda(d)$ denote the
subset of $\Lambda$ whose elements correspond to partitions with at most $d$
rows. The complexified Grothendieck group of all projective functors from
$\oplus_{\la\in\Lambda(d)}\cO_\lambda^{\mathbb{Z}}$ to
$\oplus_{\la\in\Lambda(d)}\cO_\lambda^{\mathbb{Z}}$ has also a multiplication
induced from the composition of projective functors which induces a ring
structure. Let $\op{Func}(d)$ denote the complexification of this Grothendieck
ring of all projective functors from
$\oplus_{\la\in\Lambda(d)}\cO_\lambda^{\mathbb{Z}}$ to
$\oplus_{\la\in\Lambda(d)}\cO_\lambda^{\mathbb{Z}}$. Finally let
$\mathbf{S}^{\mC}_{\mZ,v}(d,n)=\End_{\mathds{H}}(\oplus_{\la\in
\Lambda(d)}M^\la)$ be the (generic) Schur algebra attached to the
numbers $d$, $n$. Then the following holds

\begin{theorem}
\label{Schuralgebra}
There is an isomorphism of $\mC[v,v^{-1}]$-algebras
\begin{eqnarray*}
\operatorname{Func}(d)\;\cong \mathbf{S}_{\mZ,v}^{\mC}(d,n).
\end{eqnarray*}
\end{theorem}

The double centralizer property (see \cite[Theorem 4.19]{Mathas}) of the Hecke algebra $\mathds{H}^{\mC}$ for the symmetric group $S_n$ is an isomorphism
\begin{eqnarray*}
\mathds{H}^{\mC}\cong\End_{\mathbf{S}^{\mC}_{\mZ,v}(d,n)}
(\oplus_{\la\in\Lambda}\mathcal{M}^\la).
\end{eqnarray*}

It is well-known (see \cite[Theorem 3.2]{AS}) that twisting functors commute
naturally (in the sense of \cite{Khom}) with translation functors. From
Proposition~\ref{SchurWeyl} we know that the permutation parabolic modules can
be categorified via certain singular blocks of category $\cO$ together with the
action of the twisting functors. Together with the remarks of this section one
can deduce the following categorical version of the double centralizer property:
{\it The left derived functors of the graded versions of twisting functors
categorify the above action of the Schur algebra and commute naturally with
projective functors.}

\subsection{For sign parabolic modules}\label{s7.2}

Here we get the analogous result using Koszul duality (\cite{BGS}). Translation
functors should be replaced by the so-called Zuckerman functors and twisting
functors should be replaced by Irving's shuffling functors. For the Koszul
duality of these functors see \cite[Section 6]{MOS} and \cite{Steen}.

\section{Properties of $\mathscr{X}:=\mathcal{O}\{\mathfrak{p},
\mathscr{A}^{\mathbf{R}'}\}$ in case of type $A$}\label{s8}

This section describes in more detail the categories
$\mathcal{O}\{\mathfrak{p},\mathscr{A}^{\mathbf{R}'}\}$, which were
used to categorify induced cell modules, in the special case where
$\mg=\mathfrak{sl}_n$. We will describe projective-injective
modules, the associated Serre functor and show that the categories
are always Ringel self-dual.

From now on we assume that we are in the situation of
Subsection~\ref{s5.4} and will use the notation introduced there.
Additionally we assume that the Lie algebra $\mg$ is of type $A$.

We fix a right cell $\mathbf{R}'$ of $W'$ and for simplicity put
\begin{displaymath}
\mathscr{X}:=\mathcal{O}\{\mathfrak{p},\mathscr{A}^{\mathbf{R}''}\}_0.
\end{displaymath}
For $(x,w)\in \mathds{I}(\mathbf{R}')$ we denote by $L^{\mathscr{X}}(xw)$ the
simple object of $\mathscr{X}$, which corresponds to $x$ and $w$. We also have
the corresponding standard module $\Delta^{\mathscr{X}}(xw)$, proper standard
module $\overline{\Delta}^{\mathscr{X}}(xw)$, indecomposable projective module
$P^{\mathscr{X}}(xw)$, and indecomposable injective module
$I^{\mathscr{X}}(xw)$. We further denote by $T^{\mathscr{X}}(xw)$ the
indecomposable tilting module in $\mathscr{X}$ whose standard filtration starts
with a submodule $\Delta^{\mathscr{X}}(xw)$ (see \cite[4.2]{Fr} for its
existence and properties). We denote by $w_0'$ the longest element in
$W'\subset W$ and $\overline{w}=w'_0w_0$ the longest element in $(W'\backslash
W)_{short}$.

\subsection{Irving-type properties}\label{s8.1}

The following theorem is a generalization of both \cite[Main result]{Irself}
and  \cite[Proposition~4.3]{Irself}.

\begin{theorem}\label{irving}
Let $(x,w)\in \mathds{I}(\mathbf{R}')$. Then
the following conditions are equivalent:
\begin{enumerate}[(a)]
\item \label{irving.1} $L^{\mathscr{X}}(x,w)$ occurs in the socle
of some standard module from $\mathscr{X}$.
\item \label{irving.2} $L^{\mathscr{X}}(x,w)$ occurs in the socle
of some proper standard module from $\mathscr{X}$.
\item \label{irving.25} $L^{\mathscr{X}}(x,w)$ occurs in the socle
of some tilting module from $\mathscr{X}$.
\item \label{irving.3} $P^{\mathscr{X}}(x,w)$ is injective.
\item \label{irving.4} $P^{\mathscr{X}}(x,w)$ is tilting.
\item \label{irving.5} $xw\in W$ belongs to the same right cell
$\tilde{\mathbf{R}}$ of $W$ as $x\overline{w}$.
\end{enumerate}
\end{theorem}

\begin{remark}\label{rem070106-1}
{\rm As $\mathbf{R}'$ is a right cell of $W'$, we have that with our fixed
$\overline{w}$ all the $y\overline{w}$, where $y$ runs through $\mathbf{R}'$,
are in the same right cell of $W$. We denote this right cell by
$\tilde{\mathbf{R}}$, see the condition \eqref{irving.5} above. }
\end{remark}

\begin{proof}
Since the parabolic induction is exact, Consequence \eqref{cons1} from
Subsection~\ref{s4.2}, and the definition of proper standard modules as induced
simple modules, implies that $\overline{\Delta}^{\mathscr{X}}(xw)$ is a
submodule of $\Delta^{\mathscr{X}}(xw)$, hence
\eqref{irving.2}$\Rightarrow$\eqref{irving.1}. Since any standard module has a
proper standard filtration we also have
\eqref{irving.1}$\Rightarrow$\eqref{irving.2}.

Analogously, as $\Delta^{\mathscr{X}}(xw)\subset T^{\mathscr{X}}(xw)$ and
tilting modules have standard filtrations, the equivalence
\eqref{irving.1}$\Leftrightarrow$ \eqref{irving.25} is clear.

Consequence \eqref{cons1}  from Subsection~\ref{s4.2} implies that for each
$x\in \mathbf{R}'\subset W'\subset W$ the module
$\Delta^{\mathscr{X}}(x\overline{w})$ is both standard and costandard, hence
tilting, and that the socle of $\Delta^{\mathscr{X}}(x\overline{w})$ is
isomorphic to $L^{\mathscr{X}}(x\overline{w})$. Let $\theta$ be a projective
functor and $\theta'$ be its adjoint. For any
$(y,w)\in  \mathds{I}(\mathbf{R}')$ we have
\begin{displaymath}
\mathrm{Hom}_{\mathscr{X}}(L^{\mathscr{X}}(yw),\theta
\Delta^{\mathscr{X}}(x\overline{w}))=
\mathrm{Hom}_{\mathscr{X}}(\theta' L^{\mathscr{X}}(yw),
\Delta^{\mathscr{X}}(x\overline{w})).
\end{displaymath}
Since projective functors respect the right order (Proposition~\ref{prop5}),
the latter space can be non-zero only if $yw\geq_{\mathsf{R}}x\overline{w}$ in
the right order. Since $\overline{w}$ is the longest element in $(W'\backslash
W)_{short}$ and $\mathbf{R}'$ is a right cell, it follows that $yw$
is in the same right cell than $x\overline{w}$. From the proof of
Theorem~\ref{thm53}  (namely from the formula \eqref{eqthm53.2}) it follows
that, translating the tilting module $\Delta^{\mathscr{X}}(x\overline{w})$
inductively through the walls, we obtain, as direct summands, all
indecomposable tilting modules in $\mathscr{X}$. The equivalence
\eqref{irving.25}$\Leftrightarrow$\eqref{irving.5} follows.

A module which is projective and injective, is in particular tilting (since it
has a standard and a proper costandard filtration). On the other hand, a
tilting module has by definition a  standard filtration and a proper costandard
filtration, but by the construction  described above even a costandard
filtration and a proper standard filtration? Hence the dual module of a tilting
module is again tilting. By weight arguments, it is isomorphic to the original
tilting module. Hence a projective tilting module is also injective and so
\eqref{irving.3}$\Leftrightarrow$\eqref{irving.4}.

By Proposition~\ref{p22}, the category $\mathscr{A}^{\mathbf{R}'}$ has a
simple projective module. Using this and \cite[Theorem~1]{DOF} one shows that
$\mathcal{O}\{\mathfrak{p},\mathscr{A}^{\mathbf{R}'}\}$ has a simple projective
module (this statement also follows from Consequence~\eqref{cons3} and
\cite[3.1]{IS}). Translating this module out of the wall one gets that there
is at least one indecomposable projective module in $\mathscr{X}$ which is also
injective.  As we have seen already, this module must be then of the form
$P^{\mathscr{X}}(xw)$ for some  $(x,w)\in\mathds{I}(\mathbf{R}')$ such that
$xw\in\tilde{\mathbf{R}}$. Applying to $P^{\mathscr{X}}(xw)$ projective functors
we get that $P^{\mathscr{X}}(yu)$ is both projective and injective for all
$(y,u)\in\mathds{I}(\mathbf{R}')$ such that $yu\in\tilde{\mathbf{R}}$. Hence,
finally, \eqref{irving.3}$\Leftrightarrow$\eqref{irving.5}.
\end{proof}

\begin{remark}
{\rm The following statements from Theorem~\ref{irving} do not require the
additional assumption that $\mg$ is of type $A$: $\eqref{irving.1}
\Leftrightarrow\eqref{irving.2}\Leftrightarrow\eqref{irving.25}\Leftrightarrow
\eqref{irving.5}$, $\eqref{irving.3}\Leftrightarrow\eqref{irving.4}
\Rightarrow\eqref{irving.25}$. We use that $\mg$ is of type $A$ when we refer
to \cite[3.1]{IS} in the last paragraph of the proof (in particular, using
\cite{IS} the complete statement of Theorem~\ref{irving} extends to some other
special cases treated in \cite{IS}, but not to the general case because of the
counterexample from \cite[5.1]{IS}). We believe, however, that the whole Theorem~\ref{irving}  holds for arbitrary type, but do not have a complete
argument. Basically, to complete the proof for arbitrary type one has to show
that $\mathscr{X}$ always  contains a projective-injective module. }
\end{remark}

\subsection{Double centralizer and the center}\label{s8.2}

Recall that an algebra $R$ has the {\it double centralizer property} with respect to an $R$-module $M$, if there is an algebra isomorphism
\begin{displaymath}
R\cong\mathrm{End}_{\mathrm{End}_{R}(M)}(M).
\end{displaymath}
If now $R$ has a double centralizer property with respect to a module $M$ and
$\ccC\cong \mathrm{mod-}R$, then we also say that $\ccC$ has the double
centralizer property (with  respect to the image of $M$ under the equivalence).
We call a module $M\in \mathrm{mod-}R$ projective-injective, if it is both
projective and injective. If it is a direct sum of non-isomorphic
indecomposable projective-injective modules, exactly one from each isomorphism
class, then we call the module a {\it full basic projective-injective module}.

The following statement is a generalization of \cite[Theorem~5.2(ii)]{MS2}:

\begin{proposition}\label{pr991}
The category $\mathscr{X}$ satisfies the double centralizer property
with respect to any full basic projective-injective module.
\end{proposition}

To prove the statement we first need a generalization of \cite[Lemma~4.7]{MS2}:

\begin{lemma}\label{l992}
Let $\mathbf{R'}$ be a right cell in $W'\subset W$ and $x\in \mathbf{R'}$.
Then the socle of $\Delta^{\mathscr{X}}(xe)\in \mathscr{X}$ is simple.
\end{lemma}

\begin{proof}
Theorem~\ref{irving} ensures the existence of projective-injective tilting
modules in $\mathscr{X}$. Translation functors preserve the category of
projective-injective tilting modules. Any translation of a module with standard
filtration has a standard filtration. Further, from the combinatorics in the
proof of Theorem~\ref{thm53} it follows that any standard module can be
translated to a module, whose standard filtration contains
$\Delta^{\mathscr{X}}(xe)$ as a subquotient. Therefore, any
projective-injective tilting module can be translated to some
projective-injective tilting module $T$, whose standard filtration contains
$\Delta^{\mathscr{X}}(xe)$ as a subquotient. The module  $T$ contains
$\Delta^{\mathscr{X}}(xe)$ as a submodule since $\Delta^{\mathscr{X}}(xe)$ is
projective, and hence $T^{\mathscr{X}}(xe)$ is a direct summand of $T$. Thus
$T^{\mathscr{X}}(xe)$ is projective-injective, in particular, has
simple socle. As $\Delta^{\mathscr{X}}(xe)\hookrightarrow
T^{\mathscr{X}}(xe)$, the claim follows.
\end{proof}

\begin{proof}[Proof of Proposition~\ref{pr991}]
Let $x\in\mathbf{R'}$. Then the inclusion
$\Delta^{\mathscr{X}}(xe)\hookrightarrow T^{\mathscr{X}}(xe)$
extends to a short exact sequence of the following form:
\begin{equation}\label{eq991-1}
0\to \Delta^{\mathscr{X}}(xe)\to T^{\mathscr{X}}(xe)\to K\to 0,
\end{equation}
where $K\in \mathscr{F}(\Delta^{\mathscr{X}})$. The module
$T^{\mathscr{X}}(xw)$ is projective-injective by Lemma~\ref{l992}.
Projective functors are exact and preserve $\mathscr{F}(\Delta^{\mathscr{X}})$.
Hence, applying to \eqref{eq991-1} appropriate projective functors
and taking the direct sum over all $(y,w)\in \mathds{I}(\mathbf{R}')$, we
get an exact sequence
\begin{displaymath}
0\to P^{\mathscr{X}}\to M_1\to M_2\to 0,
\end{displaymath}
where $P^{\mathscr{X}}$ is a projective generator for $\mathscr{X}$, while
$M_1$ is projective-injective and $M_2\in \mathscr{F}(\Delta^{\mathscr{X}})$.
By Theorem~\ref{irving}, the injective envelope of $M_2$ is projective. The
statement now follows from \cite[Theorem~2.8]{KSX}.
\end{proof}

\begin{corollary}\label{c993}
Let $Q^{\mathscr{X}}$ denote a full basic  projective-injective module of
$\mathscr{X}$. Then the centres of $\mathscr{X}$ and $\mathrm{End}_{\mathscr{X}}(Q^{\mathscr{X}})$ are isomorphic.
\end{corollary}

\begin{proof}
The center of $\mathscr{X}$ is isomorphic to the center of
$\mathrm{End}_{\mathscr{X}}(P^{\mathscr{X}})$, where $P^{\mathscr{X}}$ is a
projective generator. Thanks to Proposition~\ref{pr991}, the centres of
$\mathrm{End}_{\mathscr{X}}(P^{\mathscr{X}})$ and
$\mathrm{End}_{\mathscr{X}}(Q^{\mathscr{X}})$ are
isomorphic (see \cite[Theorem~5.2(ii)]{MS2} for details).
\end{proof}

Because of Proposition~\ref{cunique} we can now assume that  there is a
parabolic subgroup, $W''$ of $W'$ such that $\mathbf{R}'$ contains the element
$w''_0w'_0$, where $w'_0$ and $w''_0$ denote the longest elements in $W'$ and
$W''$ respectively. Set $S''=W''\cap S'$. Let  $\mathfrak{q}$ denote the
parabolic subalgebra of $\mathfrak{g}$, which contains $\mathfrak{b}$ and such
that the Weyl group of its Levi factor is $W''$. Both  $\mathscr{X}$
and $\mathcal{O}_0^{\mathfrak{q}}$ are subcategories of the category
$\cO_0$ for $\mathfrak{g}$ and we have the following result:

\begin{lemma}\label{Tscoincide}
\begin{enumerate}[(i)]
\item\label{Tscoincide.1} $\mathscr{X}$ is a subcategory of
$\mathcal{O}_0^{\mathfrak{q}}$.
\item\label{Tscoincide.2} The projective-injective modules in $\mathscr{X}$
and $\mathcal{O}_0^{\mathfrak{q}}$, considered as objects in $\cO_0$, coincide.
\end{enumerate}
\end{lemma}

\begin{proof}
For $s\in S''$ we obviously have $sw''_0w'_0>w''_0w'_0$. As $w''_0w'_0\in
\mathbf{R}'$ and $\mathbf{R}'$ is a right  cell, it follows that $sxw>xw$ for
all $s\in S''$ and  $(x,w)\in\mathds{I}(\mathbf{R}')$. In particular,
$L^{\mathscr{X}}(xw)\in \mathcal{O}_0^{\mathfrak{q}}$
$(x,w)\in\mathds{I}(\mathbf{R}')$, which implies \eqref{Tscoincide.1}.

Consider now the element $w''_0w'_0\overline{w}=w''_0w_0'w_0'w_0=w_0''w_0$.
Then the module $P^\mathfrak{q}(w_0''w_0)$ is
projective-injective in $\mathcal{O}_0^{\mathfrak{q}}$
(see \cite{KMS} for details). As a $\mathfrak{g}$-module, the module
$P^\mathscr{X}(w''_0w'_0\overline{w})$ has simple top $L(w_0''w_0\cdot 0)$.
Hence, as $P^\mathscr{X}(w''_0w'_0\overline{w})\in
\mathcal{O}_0^{\mathfrak{q}}$ by \eqref{Tscoincide.1}, we get that
$P^\mathscr{X}(w''_0w'_0\overline{w})$ is a quotient of
$P^\mathfrak{q}(w_0''w_0)$.

On the other hand, from the existence of a simple projective module in
$\mathcal{O}^{\mathfrak{q}}$ (see \cite[3.1]{IS}) it follows that
$P^\mathfrak{q}(w_0''w_0)$ is a direct summand of some translation of
$L(w_0''w_0)$ (see Conseqeuence~\ref{cons3} in Subsection~\ref{s4.2}),
which, in turn, is the simple quotient of
$P^\mathscr{X}(w''_0w'_0\overline{w})$. Hence
$P^\mathfrak{q}(w_0''w_0)$ is a quotient of some translation of
$P^\mathscr{X}(w''_0w'_0\overline{w})$. As $P^\mathfrak{q}(w_0''w_0)$
has simple top, it follows that the only possibility is that
$P^\mathfrak{q}(w_0''w_0)$ is a quotient of
$P^\mathscr{X}(w''_0w'_0\overline{w})$.

The above implies that the $\mathfrak{g}$-modules $P^\mathfrak{q}(w_0''w_0)$
and $P^\mathscr{X}(w_0''w_0'\overline{w})$ are isomorphic, and
the claim \eqref{Tscoincide.2} follows by applying projective functors.
\end{proof}

Lemma~\ref{Tscoincide} implies the following result:

\begin{proposition}\label{pcentre}
The algebra $\mathrm{End}_{\mathscr{X}}(Q^{\mathscr{X}})$ is symmetric.
The center of $\mathscr{X}$ is isomorphic to the center of $\mathcal{O}_0^{\mathfrak{q}}$.
\end{proposition}

\begin{proof}
By Lemma~\ref{Tscoincide}, the first statement is nothing else than
\cite[Theorem 4.6]{MS2}. The second statement is given by Corollary~\ref{c993}.
\end{proof}

\begin{remark}\label{rem994}
{\rm Recall our assumption that $\mg$ is of type $A$. In this case
the center of $\mathcal{O}^{\mathfrak{q}}_0$ has a nice geometric
description: it is isomorphic to the cohomology algebra of a certain
Springer fiber. This is described in \cite{Br} and \cite{St3}. }
\end{remark}

\subsection{The Serre functor for $\mathcal{D}^p(\mathscr{X})$}\label{s8.3}

Let $\mathcal{D}^p(\mathscr{X})$ denote the full subcategory of
$\mathcal{D}^b(\mathscr{X})$ given by perfect complexes, that is complexes
which are quasi-isomorphic to finite complexes of projective objects from
$\mathscr{X}$.

Recall that if $\mathscr{C}$ is a $\Bbbk$-linear additive category
with finite-dimensional homomorphism spaces, then a {\em Serre functor}
on $\mathscr{C}$ is an auto-equivalence $\mathrm{F}$ of $\mathscr{C}$
such that the bifunctors $(X,Y)\mapsto\mathscr{C}(X,\mathrm{F}Y)$ and
$(X,Y)\mapsto\mathscr{C}(Y,X)^*$ are isomorphic (here, $*$ denotes the ordinary
duality of vector spaces).

Denote by $\mathrm{Coapp}_{\mathbf{R}'}:\mathscr{X}\to\mathscr{X}$ the functor
of partial coapproximation with respect to a fixed full basic
projective-injective module $Q^{\mathscr{X}}$. It is constructed as follows
(see \cite[2.5]{KM} for details): If $M\in \mathscr{X}$ then
$\mathrm{Coapp}_{\mathbf{R}'}(M)$ is obtained from $M$ by first maximally
extending $M$ using simple modules, which do not occur in the top of
$Q^{\mathscr{X}}$, and afterwards deleting all occurrences of such
modules in the top part.

\begin{proposition}\label{prserre}
The functor $\mathcal{R}\,\mathrm{Coapp}_{\mathbf{R}'}^2$ is a Serre
functor for $\mathcal{D}^{p}(\mathscr{X})$.
\end{proposition}

\begin{proof}
Thanks to Proposition~\ref{pr991} and Proposition~\ref{pcentre}, we are in the
situation of \cite[Theorem~3.7]{MS2}, except that the category $\mathscr{X}$
usually does not have finite global dimension. Using
\cite[Proposition~20.5.5(i)]{Gi} (see \cite[4.3]{MS2} for details), one can
get rid of the assumption of finite global dimension by  working with the
category of perfect complexes instead of the bounded  derived category.
\end{proof}

\subsection{Ringel self-duality of $\mathscr{X}$}\label{s8.4}

Consider the module
\begin{displaymath}
T^{\mathscr{X}}=\bigoplus_{(x,w)\in\mathds{I}(\mathbf{R}')}T^{\mathscr{X}}(xw).
\end{displaymath}
Based on \cite{Ri}, the algebra $\mathrm{End}_{\mathscr{X}}(T^{\mathscr{X}})$
is called the {\em Ringel dual} of the algebra
$\mathrm{End}_{\mathscr{X}}(P^{\mathscr{X}})$, see \cite{Fr}. If
$\mathbf{R}'=\{e\}$, the category $\mathscr{X}$ is {\it Ringel self-dual}
(that is $\mathrm{End}_{\mathscr{X}}(P^{\mathscr{X}}) \cong
\mathrm{End}_{\mathscr{X}}(T^{\mathscr{X}})$) by \cite[Section~7]{So4}. If
$\mathbf{R}'=\{w'_0\}$, the category $\mathscr{X}$ is Ringel self-dual
by \cite[Theorem~3]{FKM2}, see also \cite[Proposition~4.9]{MS2}. The following
theorem generalizes both these results:

\begin{theorem}\label{trsd}
The category $\mathscr{X}$ is Ringel self-dual for each
$\mathbf{R}'$.
\end{theorem}

\begin{proof}
We retain all assumptions and notation from Subsection~\ref{s8.2}
(especially the ones before Lemma~\ref{Tscoincide}). To prove this
statement we will construct an endofunctor
$\mathrm{F}=\mathrm{F}_2\mathrm{F}_1$ on $\mathcal{O}$ which maps
$P^{\mathscr{X}}$ to $T^{\mathscr{X}}$ preserving the endomorphism
ring. The functor $\mathrm{F}_2$ is an auto-equivalence of
$\mathcal{O}$ which is easy to describe: Since $\mg$ is assumed to
be of type $A$, the Dynkin diagram has an involution which is on any
$A_n$-component just the flip mapping the $i$-th vertex to the
$(n+1-i)$-th vertex. This involution induces an automorphism $\phi$ of
$\mg$, and $F_2$ maps a module $M$ to $M^\phi$, the same vector
space with the $\mg$-action twisted by $\phi$. The functor
$\mathrm{F}_1$ is more complicated. Let $w_0=s_{i_1}s_{i_2}\cdots
s_{i_{l(w_0)}}$ be a reduced expression. Consider the twisting functor
\begin{displaymath}
\mathrm{T}:=\mathrm{T}_{i_1}\mathrm{T}_{i_2}\cdots \mathrm{T}_{i_{l(w_0)}}
:\mathcal{O}\to \mathcal{O}
\end{displaymath}
(see Subsection~\ref{s25.5} and then \cite{AS}, \cite{So4} for
details). This functor is right exact, commutes with projective
functors, and $\mathcal{L}\mathrm{T}$ is a self-equivalence of
$\mathcal{D}^b(\mathcal{O})$,  see \cite{AS}. We define
$\mathrm{F}_1= \mathcal{L}_{l(w''_0)}\mathrm{T}$ and claim that
$\mathrm{F}=\mathrm{F}_2\mathrm{F}_1$ does the required job. The
arguments to deduce this are very much along the lines of
\cite[Proposition~4.4]{MS2}. Here we just outline the arguments
leaving to work out the details (following
\cite[Proposition~4.4]{MS2}) to the reader.

Denote by $\mathfrak{k}$ the semi-simple part of the Levi factor of
$\mathfrak{q}$. Then each finite-dimensional simple
$\mathfrak{k}$-module $M$ comes along with its so-called
BGG-resolution (see \cite{BGG}), that is a resolution by (direct
sums of) Verma modules, The (exact) parabolic induction functor from
$\mathfrak{k}$ to $\mathfrak{a}$ can be applied to the
BGG-resolution, and we obtain a resolution of
$M'=\cU(\ma)\otimes_{\cU(\mathfrak{k}+(\mathfrak{b}\cap\mathfrak{a}))}M$
by (direct sums of) Verma $\mathfrak{a}$-modules.

From the uniqueness result, Remark~\ref{rem2} and \cite[3.1]{IS}, there
is a simple, projective object in $\mathscr{A}^{\mathbf{R}'}$ which is
parabolically induced from a simple finite-dimensional
$\mathfrak{k}$-module (see also Subsection~\ref{s9.3} for details).
This is the $M'$ we want to consider. Its resolution by (direct sums
of) Verma $\mathfrak{a}$-modules gives rise to a resolution of the projective
module $\Delta(\mathfrak{p},M'):=U(\mg)\otimes_{\cU(\p)}M'$ by (direct sums
of) Verma $\mathfrak{g}$-modules. Then $\mathcal{L}\mathrm{T}
\Delta(\mathfrak{p},M')=
\mathcal{L}_{l(w''_0)}\mathrm{T}\Delta(\mathfrak{p},M')$ (following the
arguments in \cite{MS2}), and the latter becomes a dual parabolic Verma module.

From the construction of $\mathscr{X}$ we know that each projective
in $\mathscr{X}$ can be obtained as a direct summand of some
translation of $\Delta(\mathfrak{p},M')$. The previous paragraph
says that $\Delta(\mathfrak{p},M')$ has an (explicitly given)
resolution by (direct sums of) Verma $\mathfrak{g}$-modules.

From \cite{AL} and \cite{AS} (see also Proposition~\ref{SchurWeyl}), we
have explicit formulas for the action of the functor $\mathrm{T}$ on Verma
modules. Using these formulas one shows by a direct computation that
$\mathrm{F}$ maps $\Delta(\mathfrak{p},M')$ to a tilting module from
$\mathcal{O}\{\mathfrak{p},\mathscr{A}^{\mathbf{R}'}\}$. As
$\mathrm{F}$ commutes (up to the automorphism defining
$\mathrm{F}_2$) with projective functors, it follows that
$\mathrm{F}$ sends projective modules from $\mathscr{X}$ to tilting
modules from $\mathscr{X}$. Finally, as both $\mathrm{F}_2$ and
$\mathrm{T}$ are equivalences, it also follows that $\mathrm{F}$
preserves the endomorphism ring. This completes the proof.
\end{proof}

\section{The rough structure of generalized Verma modules}\label{s9}

In this section we want to apply the results of the paper to
determine the `rough structure' of generalized Verma modules. We
will start by giving some background information.

\subsection{Basic questions}\label{s9.1}

Let $\mathfrak{g}$ be a Lie algebra with the triangular decomposition
$\mathfrak{g}=\mathfrak{n}_-\oplus\mathfrak{h}\oplus \mathfrak{n}_+$. Let
$\mathfrak{p}\supset \mathfrak{h}\oplus \mathfrak{n}_+$  be a parabolic
subalgebra of $\mathfrak{g}$, and $V$ a simple  $\mathfrak{p}$-module,
annihilated by the nilpotent radical of $\mathfrak{p}$. The module
\begin{displaymath}
\Delta(\mathfrak{p},V)=U(\mathfrak{g})\otimes_{U(\mathfrak{p})}V
\end{displaymath}
is usually called the {\em generalized Verma module} (or simply GVM)
associated with $\mathfrak{p}$ and $V$.

If $\mathfrak{p}=\mathfrak{b}$ and $V$ is one-dimensional then we get an ordinary Verma module. If $\mathfrak{p}$ is arbitrary, but $V$ still finite dimensional, then the resulting module is a parabolic generalized Verma module as studied for example in \cite{Ja}.

The most basic questions about GVMs are:
\begin{itemize}
\item In which case is $\Delta(\mathfrak{p},V)$ irreducible?
\item If $\Delta(\mathfrak{p},V)$ is  not irreducible: which simple
$\mathfrak{g}$-modules occur as subquotients of $\Delta(\mathfrak{p},V)$, and what are their multiplicity (in case this makes sense at all)?
\end{itemize}
These questions were studied in special cases by many authors,
we refer the reader to \cite[Introduction]{KM2} for a more detailed
survey. The answer to the questions above is also of interest in theoretical physics, since the structure of generalized Verma modules determines the structure of Verma modules for (super)algebras appearing in conformal field theory (see for example \cite{Se} for an affine setup).

The most general known facts in the theory of generalized Verma modules are the main results of \cite{KM2} (based on \cite{MiSo}) under the assumption that the module $V$ has minimal possible annihilator: \cite[Theorem~22]{KM2} gives an
explicit criterion for the irreducibility of $\Delta(\mathfrak{p},V)$; and \cite[Theorem~23]{KM2} describes what is called the {\it rough}
structure of $\Delta(\mathfrak{p},V)$, defined as follows:
each $\Delta(\mathfrak{p},V)$ has a unique simple quotient, denoted
by $L(\mathfrak{p},V)$. If $V'$ is another simple $\mathfrak{p}$-module
with minimal annihilator then \cite[Theorem~23]{KM2} says that
the multiplicity $[\Delta(\mathfrak{p},V):L(\mathfrak{p},V')]$ is
well-defined (in particular, it is always finite); an explicit formula for its computation in terms of Kazhdan-Lusztig polynomials is also provided.
In general, this does not describe the structure of $\Delta(\mathfrak{p},V)$ completely: $\Delta(\mathfrak{p},V)$ may have many other subquotients, it even might be of infinite length (because of the example due to Stafford, see
\cite[Theorem 4.1]{Stafford}). No reasonable information about this so-called
{\it fine structure} of $\Delta(\mathfrak{p},V)$ is known so far.

In what follows we want to explain how one can drop the restriction on the minimality of the annihilator of $V$ by applying the techniques we have developed so far in this paper. Following the approach proposed in \cite{MiSo} and developed further in \cite{KM2}, an essential part of the argument is an improved answer to the so-called `Kostant's problem' for certain simple and induced modules.

\subsection{Kostant's problem}\label{s9.2}

Let $\mathfrak{g}$ be a complex reductive finite-dimensional Lie algebra. For
every $\mathfrak{g}$-module $M$ we have the bimodule $\mathscr{L}(M,M)$ of all
$\mathbb{C}$-linear endomorphisms of $M$, on which the adjoint action of
$U(\mathfrak{g})$ is locally finite (that means any vector
$f\in\mathscr{L}(M,M)$ lies inside a finite dimensional subspace which is
stable under the adjoint action defined as $x.f(m)=x(f(m))-f(xm)$ for $x\in\mg$, $m\in M$). Initiated by \cite{Jo}, {\it Kostant's problem} became the standard terminology  for the following question concerning an arbitrary $\g$-module M:
\bigskip

{\em
Is the natural injection $U(\mathfrak{g})/\mathrm{Ann}(M)\hookrightarrow
\mathscr{L}(M,M)$
surjective?
}
\bigskip

Although there are several classes of
modules for which the answer is known to be positive (see
\cite{Jo}, \cite{Ma2} and references therein), a complete answer to this problem seems to be far away - not even for simple
highest weight modules the problem is solved.  There is even an instance of
a simple highest weight module for which the answer is negative. The details of such an example (which was first mentioned in \cite[9.5]{Jo}) will be discussed in Subsection~\ref{B2}.

In the following we will show that for certain
simple and induced modules which appeared  already in the present
paper, the answer to Kostant's problem is positive. We conjecture that this is always the case for simple highest weight modules for Lie algebras of type $A$. Theorem~\ref{cor2-061107} shows that the answer only depends on the left cell associated with a simple highest weight module.

\subsection{General assumptions}\label{ass}

For the rest of the paper let $\mathfrak{g}$ be an arbitrary complex reductive
finite-dimensional Lie algebra with a fixed triangular decomposition
$\mathfrak{g}=\mathfrak{n}_-\oplus\mathfrak{h}\oplus\mathfrak{n}_+$.
Let $\mathfrak{p}\supset \mathfrak{h}\oplus\mathfrak{n}_+$ be a parabolic
subalgebra of $\mathfrak{g}$ with Levi factor $\mathfrak{a}'$ and nilpotent
radical $\mathfrak{n}$. Finally, denote by $\mathfrak{a}$ the semi-simple part
of $\mathfrak{a}'$. Then $\mathfrak{a}$ is a semi-simple finite-dimensional
Lie algebra with induced triangular decomposition. We {\bf assume} that
$\mathfrak{a}$ is of type $A$, that means:
\bigskip

{\em
We assume that
$\mathfrak{a}\cong\bigoplus_{i\in I}\mathfrak{sl}_{k_i}$,
where $I$ is some finite set and $k_i\in\{2,3,\dots\}$.
}
\bigskip

\subsection{Kostant's problem for the IS-module}\label{s9.3}

We have to start with some technical statements which involve explicit definitions of certain weights. Assume  $\mathfrak{a}\cong\mathfrak{sl}_{k}$ for some $k\geq 2$. We consider $\mathfrak{a}$ as a subalgebra of
$\mathfrak{gl}_{k}$ in the canonical way. In particular, all simple highest
weight $\mathfrak{gl}_{k}$-modules are simple highest weight
$\mathfrak{a}$-modules via restriction. Let $\alpha_i$ ($i=1,\dots,k-1$) be
the list of simple roots of $\mathfrak{a}$ in the usual ordering. As before
we denote the Weyl group of $\mathfrak{a}$ by $W'$ and note that $W'\cong S_k$.
Let $r$ be a  partition of $k$ of length $s$, that is
$r=(r_1,\dots,r_s)\in\mathbb{N}^s$, $r_1+\dots+r_s=k$ and $r_1\geq
r_2\geq\dots$. Set $r_0=0$. Depending on $r$,  we define
$\pi=\{\alpha_i\,:i\in I\}$, where
{\small
\begin{displaymath}
I=\{1,2,\dots,r_1-1\}\cup \{r_1+1,\dots,
r_1+r_2-1\}\cup\dots\cup\{r_1+\dots+r_{s-1}+1,\dots, r_1+\dots+r_{s}-1\}
\end{displaymath}
}
and then the $\mathfrak{gl}_{k}$-weight $\nu$ as
\begin{displaymath}
\nu=(b_1,\dots,b_k),\quad
b_{r_{j-1}+m}=r_j-m, \text{ for } m\in\{1,\dots,r_{j}\}.
\end{displaymath}
In \cite[Proposition~3.1]{IS}, it is shown that
$\nu$ is the only $\pi$-dominant weight in
$W'\nu$ and hence the corresponding simple highest
weight module $L(\nu-\rho)$ is a projective simple module in
the parabolic category $\mathcal{O}$ associated with $\pi$. This is what we call the {\it simple projective IS-module}.

Denote by $\mu$ the weight such that $\mu+\rho$ is the dominant
weight in $W'\nu$. To proceed we have to construct Weyl group
elements $x_\nu$, $x_\mu$ such that  $L(\nu-\rho)$ is the
translation of $L(x_\nu\cdot0)$ to $\cO_\mu$, and $L(\mu)$ is the
translation of $L(x_\mu\cdot0)$ to $\cO_\mu$.

Let $\xi=(\xi_1,\dots, \xi_k)$ be a $k$-tuple of non-negative
integers. We convert the coordinates of $\xi$ into the sequence
$(\eta_1,\dots,\eta_k)$ without repetitions, which differs from
$(k-1,k-2,\dots,0)$ only by a permutation, and satisfies
$\eta_j<\eta_k$ if $j<k$ or $\xi_j<\xi_k$ (in practice we first
replace all occurring zeros from the left to the right by $0$,
$1,\dots, m_0$, where $m_0+1$ is the total number of zeros in $\xi$,
then all occurring ones by $m_0+1$, $m_0+2$ etc.). Applying this
procedure to $\nu$ and $\mu+\rho$ we obtain weights $\nu'+\rho$ and
$\mu'+\rho$ from the orbit $W'(k-1,k-2,\dots,0)$. Then
$\nu'+\rho=x_\nu(k-1,k-2,\dots,0)$ and
$\mu'+\rho=x_\mu(k-1,k-2,\dots,0)$ for some $x_\nu$, $x_\mu\in W'$.
By construction, $L(x_\nu\cdot0)$ and $L(x_\mu\cdot0)$ are simple
highest weight modules with the desired properties described above.

\begin{example}{\rm
Consider the case where $\ma=\mathfrak{sl}_4$ with the three simple reflections
$s_1$, $s_2$, $s_3$, where $s_1$ and $s_3$ commute. The partition $r=(2,2)$
gives $\pi=\{\alpha_1,\alpha_3\}$ and $\nu=(1,0,1,0)$. Then
$\mu+\rho=(1,1,0,0)$, $\nu'+\rho=(2,0,3,1)$ and  $\mu'+\rho=(2,3,0,1)$. Hence
$x_\nu=s_2s_1s_3$, $x_\mu=s_1s_3$.
}
\end{example}

Our crucial technical observation is the following

\begin{lemma}\label{lem1-061107}
$x_\nu$ and $x_\mu$ belong to the same left cell.
\end{lemma}

\begin{proof}
We prove this by induction on $k$. If $k=2$, there is nothing to  prove. If
$r_1>r_2$, then in both, $\nu'+\rho$ and $\mu'+\rho$, the element
$k-1$ stays at the leftmost place, and the induction hypothesis
applies to the remaining parts of $\nu'+\rho$ and $\mu'+\rho$.
The only tricky part is therefore the case $r_1=r_2$, which may in fact
easily be reduced inductively to the case $r_1=r_2=\dots=r_s$.
Consider first the case $s=2$. Then $x_\nu$ is the following permutation
on $\{0,\dots,k-1\}$, which we consider as an element of $W'$:
{\small
\begin{displaymath}
x_\nu=
\left(
\begin{array}{ccccccccc}
0&1&\dots&r_1-1&r_1& r_1+1& r_1+2&\dots&m\\
m-1& m-3&\dots&2& 0&m &m-2 &\dots&1
\end{array}
\right).
\end{displaymath}
}
where $m=r_1+r_2-1$. Since $0<2<m$, we can apply Knuth transformation
(see \cite[Definition~3.6.8]{Sa}) to interchange
$0$ and $m$ in the second row of the above permutation.
This can be continued until $m$ appears at the second
left position, where the procedure stops.
Since the Knuth transformations preserve left cells
(\cite[Lemma~3.6.9]{Sa}), the new permutation $\sigma$ will be in
the same left cell as $x_\nu$. Now in $\sigma$ and $x_\mu$
the first two elements coincide. So, applying the induction hypothesis
to the remaining parts, we get that $x_\mu
$ and $x_\nu$ are in the same left cell. The case $s>2$ follows now
inductively. We omit the details.
\end{proof}

The following result is crucial an its proof is based on the
categorification results from Section~\ref{s6} and Section~\ref{s7}:

\begin{theorem}\label{cor2-061107}
\begin{enumerate}[(i)]
\item\label{cor2-061107.1} The modules
$L(\nu-\rho)$ and $L(\mu)$ have the same annihilator.
\item\label{cor2-061107.2} For any projective functor $\theta$
we have
\begin{displaymath}
\dim\mathrm{Hom}_{\mathfrak{a}}(L(\nu-\rho),\theta L(\nu-\rho))=
\dim\mathrm{Hom}_{\mathfrak{a}}(L(\mu),\theta L(\mu)).
\end{displaymath}
\item\label{cor2-061107.3}
Kostant's problem has a positive answer for $L(\mu)$ and $L(\nu-\rho)$.
\item \label{cor2-061107.2a} For any projective functor $\theta$
we have
\begin{displaymath}
\dim\mathrm{Hom}_{\mathfrak{a}}(L(x\cdot0),\theta L(x\cdot0))=
\dim\mathrm{Hom}_{\mathfrak{a}}(L(y\cdot0),\theta L(y\cdot0))
\end{displaymath}
whenever $x$ and $y$ are in the same left cell of $W'$. In particular,
Kostant's problem has a positive answer for $L(x\cdot0)$ if and only if it
has a positive answer for $L(y\cdot0)$.
\end{enumerate}
\end{theorem}

\begin{proof}
The annihilators of the modules  $L(\nu')$ and $L(\mu')$ coincide since
$x_\nu$ and $x_\mu$ belong to the same left cell by Lemma~\ref{lem1-061107}. The statement \eqref{cor2-061107.1} is now obtained by translating
to the wall and applying \cite[5.4 (3)]{Ja2}.

We will see later that the statement \eqref{cor2-061107.2} follows from \eqref{cor2-061107.2a}. To prove \eqref{cor2-061107.2a} we have to work
much harder. The principal idea is the following: Given two simple modules in
the same block, and indexed by elements in the same left cell, then
Proposition~\ref{SchurWeyl} tells us that they are connected via twisting
functors. These twisting functors commute with projective functors and
therefore they could be used to obtain estimates for the dimensions of
homomorphism spaces, which would result in \eqref{cor2-061107.2a}. Let us make
this idea precise. Assume $x\in W$ and $s$ is a simple reflection such
that $sx<x$, and the elements $sx$ an $x$ belong to the same left cell.
For the twisting functor  $\mathrm{T}_s:\cO\rightarrow\cO$ we have:
\begin{eqnarray}
\label{Ts}
\mathrm{Hom}_{\mathfrak{a}}(\theta L(x\cdot0),L(x\cdot0))&\cong&
\mathrm{Hom}_{\mathfrak{a}}(\mathcal{L}\mathrm{T}_s\theta L(x\cdot0),
\mathcal{L}\mathrm{T}_s L(x\cdot0))\\
&\cong&
\mathrm{Hom}_{\mathfrak{a}}(\theta\mathcal{L}\mathrm{T}_s L(x\cdot0),
\mathcal{L}\mathrm{T}_s L(x\cdot0))\nonumber\\
&\cong&
\mathrm{Hom}_{\mathfrak{a}}(\theta\mathrm{T}_s L(x\cdot0),
\mathrm{T}_s L(x\cdot0))\nonumber\\
&\cong&\mathrm{Hom}_{\mathfrak{a}}(\mathrm{T}_s\theta L(x\cdot0),
\mathrm{T}_s L(x\cdot0))\nonumber.
\end{eqnarray}
by \cite[Corollary 4.2, Theorem 2.2, Theorem 6.1, Theorem 3.2]{AS}. Moreover,
we also have $\mathrm{T}_sL(x\cdot0)\not=0$ and a short exact sequence
\begin{eqnarray}
0\to U\longrightarrow\mathrm{T}_s L(x\cdot0)
\overset{\mathrm{nat}}{\longrightarrow} L(x\cdot0)\to 0,
\end{eqnarray}
where $\mathrm{nat}$ is the evaluation at $L(x\cdot0)$ of the
natural transformation from $\mathrm{T}_s$ to the identity functor,
given by \cite[Theorem~4]{KM}, and $U$ is the kernel of
$\mathrm{nat}$. Further, the module $U$ is semi-simple, and has
$L(sx\cdot0)$ as a simple subquotient with multiplicity one by
\cite[Section~7]{AS}. As  $L(x\cdot0)$ is simple and $s$-infinite,
the module $U$ coincides with the maximal $s$-finite submodule of
$\mathrm{T}_sL(x\cdot0)$, see \cite[Proposition~5.4]{AS} and
\cite[2.5 and Theorem~10]{KM}.

Analogously we have a short exact sequence
\begin{eqnarray}
0\to U'\longrightarrow\mathrm{T}_s \theta L(x\cdot0)
\overset{\mathrm{nat}}{\longrightarrow}
\theta L(x\cdot0)\to 0,
\end{eqnarray}
where $U'$ is just the kernel of $\mathrm{nat}$. As all simple
submodules of the socle of $\theta L(x\cdot0)$ are $s$-infinite, the
module $U'$ again coincides with the maximal $s$-finite submodule of
$\mathrm{T}_s\theta L(x\cdot0)$. This implies $U'\cong\theta U$.
Now, any non-zero homomorphism $f\in \mathrm{Hom}_{\mathfrak{a}}(\theta
L(x\cdot0),L(x\cdot0))$ is automatically surjective and gives rise to a
diagram as follows (in which  the square of solid arrows commutes):
\begin{eqnarray*}
\xymatrix{
\theta U\ar@{^{(}->}[r]
\ar@{.>}[d]^{f'}&\mathrm{T}_s\theta L(x\cdot0)
\ar@{->>}[r]^{\mathrm{nat}}
\ar@{->>}[d]^{\mathrm{T}_s f}&\theta L(x\cdot0)\ar@{->>}[d]^{f}\\
U\ar@{^{(}->}[r]&\mathrm{T}_s L(x\cdot0)
\ar@{->>}[r]^{\mathrm{nat}}&L(x\cdot0),
}
\end{eqnarray*}
inducing the map $f'$. We claim that  the map $f'$ restricts  to a
non-zero map
\begin{eqnarray*}
\overline{f}'\in\mathrm{Hom}_{\mathfrak{a}}(\theta L(sx\cdot0),L(sx\cdot0)).
\end{eqnarray*}
We first claim that the cokernel of $f'$ does not contain any simple
module $L(z)$, where $z$ is in the same left cell as $x$. By the
Snake Lemma, the cokernel of $f'$ embeds into $\theta L(x\cdot0)$.
From Proposition~\ref{prop5}\eqref{prop5.1} it follows that $\theta
L(x\cdot0)$ only has composition factors indexed by $z$'s either in
the same right cell as $x$ or in smaller right cells. From the
Robinson-Schensted algorithm it is directly clear that smaller
rights cells intersect the left cell of $x$ trivially.
Robinson-Schensted also implies that any given left and right cell
inside the same two-sided cell intersect in exactly one point; so
the only possible $z$ is $z=x$.
Since $L(x\cdot0)$ is not a composition factor of $U$ by
\cite[Theorem 6.3 (ii)]{AS} the claim follows. In particular,
$L(sx\cdot0)$ occurs in the image of $f'$. Let now $L(z\cdot0)$ be a
simple subquotient of $U$. If $z$ belongs to a smaller two-sided
cell than $sx$, then the arguments of Theorem~\ref{GKdim} imply that
the GK-dimension of $L(z\cdot0)$ is strictly smaller than that of
$L(sx\cdot0)$. Hence $\mathrm{Hom}_{\mathfrak{a}}(\theta
L(z\cdot0),L(sx\cdot0))=0$. If $z$ is in the same left cell as $sx$ then by
Proposition~\ref{prop5}\eqref{prop5.1} and Robinson-Schensted the inequality
$\mathrm{Hom}_{\mathfrak{a}}(\theta L(z\cdot0),L(sx\cdot0))\not=0$ is
possible only for $z=sx$.  This implies $\overline{f}'\not=0$ and
the claim follows. Hence we get the inequality
\begin{equation}\label{abs123}
\dim\mathrm{Hom}_{\mathfrak{a}}(\theta L(x\cdot0),L(x\cdot0))\leq
\dim\mathrm{Hom}_{\mathfrak{a}}(\theta L(sx\cdot0),L(sx\cdot0)).
\end{equation}

Analogously one obtains the inequality
\begin{equation}\label{abs1234}
\dim\mathrm{Hom}_{\mathfrak{a}}(\theta L(x\cdot0),L(x\cdot0))\leq
\dim\mathrm{Hom}_{\mathfrak{a}}(\theta L(tx\cdot0),L(tx\cdot0))
\end{equation}
in the case when $sx$ belongs to the smaller left cell than $x$ and
$t$ is a simple reflection such that $(st)^3=e$ and the element $tx$
belongs to the same left cell as $x$.

Since left cell modules are irreducible, using \eqref{abs123} and
\eqref{abs1234} inductively one also obtains the opposite inequalities, which implies that
\begin{equation}\label{abs123456}
\mathrm{Hom}_{\mathfrak{a}}(\theta L(x\cdot0),L(x\cdot0))=
\mathrm{Hom}_{\mathfrak{a}}(\theta L(y\cdot0),L(y\cdot0))
\end{equation}
if $x$ and $y$ are in the same left cell.
By \cite[6.8(3)]{Ja2} for any $\mathfrak{a}$-module $M$
and any simple  finite-dimensional $\mathfrak{a}$-module $F$
we have
\begin{equation}\label{abs12345}
[\mathscr{L}(M,M):F]=
\mathrm{Hom}_{\mathfrak{a}}(F\otimes M,M).
\end{equation}
Since the modules $L(x\cdot0)$ and $L(y\cdot0)$ have the same annihilator by
\cite[5.25]{Ja2}, the formulas \eqref{abs12345} and \eqref{abs123456} imply
that Kostant's problem has a positive answer either for both  $L(x\cdot 0)$
and $L(y\cdot 0)$ or for none of them. This proves \eqref{cor2-061107.2a}.

By Lemma~\ref{lem1-061107}, both $L(\nu-\rho)$ and $L(\mu)$ are obtained by translating two simple modules, indexed by elements from the same left cell, from $\cO_0$ to a fixed singular block. Hence, using
Proposition~\ref{SchurWeyl}, the statement \eqref{cor2-061107.2} is proved
just in the same way as the statement \eqref{cor2-061107.2a} is proved above.

As in \eqref{cor2-061107.2a}, the statement \eqref{cor2-061107.2} implies that
Kostant's problem has a positive answer either for both $L(\nu-\rho)$ and
$L(\mu)$ or for none of them. Since $\mu$ is dominant, Kostant's problem has
a positive  answer for $L(\mu)$ by \cite[6.9]{Ja2}. Hence the answer to
Kostant's problem for $L(\nu-\rho)$ is positive as well. This completes the
proof.
\end{proof}

\subsection{A negative answer to Kostant's problem:  type $B_2$}\label{B2}

The answer to Kostant's problem is not positive in general. The following
negative example was constructed first in \cite[9.5]{Jo}: Let
\begin{displaymath}
W'=\{e,s,t,st,ts,sts,tst,stst=tsts\}
\end{displaymath}
be the Weyl group of type $B_2$ with the two simple reflections $s$ and $t$ as generators. Then the requirement of Kostant's problem fails for $L(ts)$ and $L(st)$. What goes wrong in our arguments? The right cells of $W'$ are $\{e\}$, $\{s,st,sts\}$, $\{t,ts,tst\}$, $\{stst\}$. In particular, the elements $s$ and $sts$ are both in the same left and in the same right cell. In our arguments we used several times that the intersection of a given left with a given right cell contains at most one element. An easy direct calculation also shows that $\mathrm{Hom}_{\mathcal{O}}(\theta L(s),L(s))=\mC=
\mathrm{Hom}_{\mathcal{O}}(\theta L(sts),L(sts))$, whereas $\mathrm{Hom}_{\mathcal{O}}(\theta L(ts),L(ts))=\mC^2$ for $\theta=\theta_s\theta_t\theta_s$. Hence Theorem~\ref{cor2-061107} \eqref{cor2-061107.2a} fails in this case.

\subsection{Coker-categories and their equivalence}

Recall the setup from Sub\-sec\-ti\-on~\ref{ass}, in particular that
$\mathfrak{a}$ is assumed to be of type $A$. Let $V$ be an arbitrary
simple $\mathfrak{a}'$-module. Let $L$ be the $\mathfrak{a}$-module obtained
by restriction. For simplicity
\begin{eqnarray}\label{easylife}
\begin{array}{c}
\text{\it we assume for the moment that $L$ has}\\
\text{\it integral and regular central character}
\end{array}
\end{eqnarray}
and refer to Remark~\ref{simplicity} for the general case.

Now, $V$ is determined uniquely by the underlying simple
$\mathfrak{a}$-module $L$ and some functional $\eta$ on the center
of $\mathfrak{a}'$. We first construct an admissible category
attached to this data:

Let $L(x\cdot 0)$ be a simple highest weight module with the same annihilator
as $L$. Without loss of generality we assume that $x$ is contained in a right
cell associated with a parabolic subalgebra $\p$ as in Remark~\ref{rem2}
(the latter is possible as $x$ can be chosen arbitrarily in its left cell by
\cite[5.25]{Ja2}). By \cite[3.1]{IS}, there is a block $\cO_\mu^\p$  (for some
integral weight $\mu\in\mathfrak{h}^*_{dom}$) which contains exactly one
simple (highest weight) module $L(y\cdot\mu)$, and this module is also
projective. We assume that $ys<y$ for any simple reflection $s$ such that
$s\cdot \mu=\mu$. The module $L(y\cdot\mu)$ is a tensor product of simple
highest weight modules over all  simple components of $\mathfrak{a}$. Each of
the factors has the form $L(\nu)$,  where $\nu$ is as in Subsection~\ref{s9.3}.
Because of our assumptions we also have that $L(y\cdot\mu)$ is the
translation of $L(y\cdot 0)$ to the $\mu$-wall, and that $x$ and $y$ belong
to the same right cell (Consequence~\ref{cons3} in Subsection~\ref{s4.2}).

\begin{proposition}\label{F}
There is some projective functor
$F:\cO(\mathfrak{a},\mathfrak{a}\cap\mathfrak{b})\rightarrow
\cO(\mathfrak{a},\mathfrak{a}\cap\mathfrak{b})$ such that
$FL(x\cdot 0)\cong \bigoplus_{i=1}^k L(y\cdot\mu)$ for some finite
number $k>0$.
\end{proposition}

\begin{proof}
First we claim that there is a projective functor $\theta$ such that
$L(y\cdot0)$ occurs as a composition factor in $\theta L(x\cdot0)$. Indeed,
recall that the elements $x$ and $y$ are in the same  right cell  of $W'$.
Consider the basis of simple modules for the categorification of the cell
module (corresponding to $x$ and $y$) given by Theorem~\ref{thm5}. As cell
modules are irreducible in type $A$, there is a projective functor $\theta$
such that $[\theta L(x\cdot0)]$ has a non-zero coefficient at $[L(y\cdot0)]$,
when expressed with respect to the basis of simple modules. This means exactly
that $L(y\cdot0)$ occurs as a composition factor in $\theta L(x\cdot0)$.

Let $\theta'$ be the translation to the $\mu$-wall.  Then the functor
$F=\theta'\theta$ satisfies  the requirement of the lemma as the module
$L(y\cdot\mu)$ is simple projective and is a unique simple modules in its
parabolic block (see the definition of $L(y\cdot\mu)$ and \cite[3.1]{IS}).
\end{proof}

We fix $F$ as in Proposition~\ref{F} and set $N:=FL$.

\begin{lemma}\label{Nnotzero}
\begin{enumerate}[(i)]
\item\label{Nnotzero.1} $N=FL\not=0$.
\item\label{Nnotzero.2}
$\Ann_{U(\mathfrak{a})}N=\Ann_{U(\mathfrak{a})}L(y\cdot\mu)=
\Ann_{U(\mathfrak{a})}L(\mu)$.
\end{enumerate}
\end{lemma}

\begin{proof}
Since $\Ann L=\Ann L(x\cdot0)$, we have
\begin{displaymath}
\Ann_{U(\mathfrak{a})}N=\Ann_{U(\mathfrak{a})} F L=
\Ann_{U(\mathfrak{a})} F L(x\cdot0)
\end{displaymath}
(see \cite[5.4]{Ja2}). The second statement follows then directly from
Proposition~\ref{F}, Theorem~\ref{cor2-061107} and the definition of
$L(y\cdot\mu)$. Since $FL(x\cdot0)\not=0$, we also have $FL\not=0$.
\end{proof}

If $\mg$ is any complex Lie algebra and $Q$ a $\mg$-module, then we denote by $\op{Coker}(Q\otimes E)$ the full subcategory of $\mathfrak{g}\mathrm{-mod}$, which consists of all modules $M$ having
a presentation $X\to Y\tto M$, where both $X$ and $Y$ are direct summands of
$Q\otimes E$ for some finite dimensional module $E$. In particular, if we choose the Lie algebra to be $\ma$, then we have the two categories $\op{Coker}(L(y\cdot\mu)\otimes E)$ and $\op{Coker}(N\otimes E)$.

\begin{lemma}\label{projLy}
$L(y\cdot\mu)$ is projective in $\op{Coker}(L(y\cdot\mu)\otimes E)$.
\end{lemma}

\begin{proof}
Of course $L(y\cdot\mu)$ is contained in $\op{Coker}(L(y\cdot\mu)\otimes E)$.
On the other hand, $L(y\cdot\mu)\in\cO_\mu^\p$ is projective. It is even
projective in $\cO^\p$. The latter is stable under tensoring with finite
dimensional modules, hence contains $\op{Coker}(L(y\cdot\mu)\otimes E)$ as a
full subcategory. The statement follows.
\end{proof}

Unfortunately, we do not know how to prove directly that the module $N$ is
semi-simple.  To get around this problem we have to make sure that there is a
`nice' simple subquotient $\overline{N}$ of $N$. Let $G$ be the adjoint
functor of $F$. Then the adjunction morphism $a: L\rightarrow GFL$ is
injective, since $L$ is simple. Since $G$ is exact, there must therefore
be a simple subquotient $\overline{N}$ of $N$ such that $G\overline{N}$
contains $L$ as a quotient.

Let $\chi_\mu=\Ann_{Z(\ma)}M(\mu)$ be the central character of the Verma module
with highest weight $\mu$. Then we denote by $\mathcal{M}(\mu)$ the category of
all $\ma$-modules $M$ such that $(\chi_\mu)^nM=0$ for some $n$ (depending on
$M$), i.e. $M$ has generalized central character $\chi_\mu$. With this notation
the following holds:

\begin{lemma}\label{projfuncN}
Let $\theta:\mathcal{M}(\mu)\rightarrow\mathcal{M}(\mu)$ be an indecomposable projective functor which is not the identity. Then $\theta  L(\mu)=0$,
$\theta N=0$ and $\theta\overline{N}=0$.
\end{lemma}

\begin{proof}
We first show the statement for $L(\mu)$: if $\theta:\cO_\mu\rightarrow\cO_\mu$
is an indecomposable projective functor then $\theta L(\mu)\not=0$ means that
$\theta$ is the identity functor. To see this, take the projective cover
$P(\mu)$ of $L(\mu)$. Then
\begin{eqnarray}\label{hom}
\Hom_{\ma}(P(\mu),\theta L(\mu))&=&\Hom_\ma(\theta' P(\mu),L(\mu)),
\end{eqnarray}
where $\theta'$ is the adjoint functor of $\theta$. Note that $\theta'$ is an
indecomposable projective functor if so is $\theta$. The classification theorem
of projective functors gives $\theta' M(\mu)=P(\zeta)$ for some $\zeta$. If we
assume the space $\eqref{hom}$ to be non-trivial then we have $\zeta=\mu$, which
forces (by the classification theorem again) $\theta'$ to be the identity
functor, and then $\theta$ is the identity functor as well.

Assume therefore $\eqref{hom}$ is trivial, but $\theta L(\mu)\not=0$. Recall
the categorification result of Proposition~\ref{SchurWeyl} and extend the scalars to $\mathbb{C}$. Together with Theorem~\ref{Schuralgebra}  we get
that $\theta$ induces an endomorphism of the complexified
Grothendieck group of $\cO_\mu$. The module $L(\mu)$ has minimal
Gelfand-Kirillov dimension and is contained in the categorification of the
irreducible (Specht) submodule of $\mathcal{M}^\mu$ corresponding to the
partition given by $\mu$. The endomorphism of the parabolic permutation module
given by $\theta$ is a homomorphism of the symmetric group which underlies the
Hecke algebra $\mathbb{H}$ and restricts to an endomorphism of the irreducible
submodule which has to be a multiple $c\in \mathbb{C}$ of the identity. But
since $0\not=\theta L(\mu)$ has at most the same Gelfand-Kirillov dimension as
$L(\mu)$ (by Lemma~\ref{GK}), we deduce that $c\not=0$. On the other hand
the fact that both sides of the equality \eqref{hom} are equal to
$0$ is equivalent to the statement that $L(\mu)$ does not occur as a
composition factor in $\theta L(\mu)$, a contradiction. In particular,
$c\not=0$ forces $\theta$ to be the identity functor. Hence the claim is true
for $L(\mu)$.

Assume again that $\theta$ is not the identity functor. To see that $\theta
N=\theta FL=0$ we consider the annihilator $\Ann \theta
N=\Ann_{\cU(\ma)}(\theta N)$. By \cite[6.35 (1)]{Ja2} we have
\begin{eqnarray}\label{eq:Ann}
\cU(\ma)/\Ann \theta N&=&\theta^{l}\theta^{r}(\cU(\ma)/\Ann N),
\end{eqnarray}
where $\cU(\ma)/\Ann N$ is considered as a $\cU(\ma)$-bimodule, and
$\theta^{l}$ is the projective functor $\theta$ when considering left
$\cU(\ma)$-modules, whereas $\theta^{r}$ is the projective functor $\theta$
when considering right $\cU(\ma)$-modules (see also Subsection \ref{s25.4}).
On the other hand we have an equality of bimodules
\begin{eqnarray}\label{eq:anns}
\cU(\ma)/\Ann N&=&\cU(\ma)/\Ann L(\mu)=\mathscr{L}(L(\mu),L(\mu)).
\end{eqnarray}
Here the first equality is \cite[5.4]{Ja2}. The second equality is given by the
natural map, since Kostant's problem has a positive answer in this case
(Theorem~\ref{cor2-061107}\eqref{cor2-061107.3}). Putting everything together
we get
\begin{eqnarray*}
\cU(\ma)/\Ann \theta N&=&\theta^{l}\theta^{r}(\cU(\ma)/\Ann
N)=\theta^{l}\theta^{r}\mathscr{L}(L(\mu),L(\mu))\\
&\cong&\mathscr{L}(\theta L(\mu),\theta L(\mu))=0.
\end{eqnarray*}
For the penultimate isomorphism we refer to \cite[6.33(6)]{Ja2}. It
follows that $\theta N=0$. As $\theta$ is exact and $\overline{N}$ is a
subquotient of $N$ we also get that $\theta\overline{N}=0$.
\end{proof}

\begin{proposition}\label{smallN}
The following holds:
\begin{enumerate}[(i)]
\item\label{smallN1} $\Ann_{U(\mathfrak{a})}\overline{N}=
\Ann_{U(\mathfrak{a})}L(y\cdot\mu)= \Ann_{U(\mathfrak{a})}L(\mu)$.
\item\label{smallN2} $\overline{N}$ is projective in
$\op{Coker}(\overline{N}\otimes E)$.
\item\label{smallN3} Kostant's problem has a positive
solution for the module $\overline{N}$.
\end{enumerate}
\end{proposition}

\begin{proof}
Of course, $\Ann_{U(\mathfrak{a})}\overline{N}\supseteq
\Ann_{U(\mathfrak{a})}{N}$. Let us assume $\Ann_{U(\mathfrak{a})}\overline{N}$
is strictly bigger than $\Ann_{U(\mathfrak{a})}{N}$. Choose $z$, such that
$\Ann_{U(\mathfrak{a})}\overline{N}=\Ann_{U(\mathfrak{a})}L(z\cdot\mu)$. Since
$\Ann_{U(\mathfrak{a})}N=\Ann_{U(\mathfrak{a})}L(y\cdot\mu)$
(Lemma~\ref{Nnotzero}), it follows that $z$ is strictly smaller than $y$ in the
left order. Hence, also strictly smaller than $x$ in the left order. On the
other hand (by definition of the modules) $L(y\cdot0)$ can be obtained as a
subquotient in a translation of $L(z\cdot0)$.
Proposition~\ref{prop5}~\eqref{prop5.1} tells then that  $y\leq_{\mathsf{R}}z$.
In particular, $y$ is smaller than or equal to $x$ in the twosided order. This
contradicts the fact that $z$ should be strictly smaller than $y$ in the left
order. The  first statement follows.

By definition $\overline{N}\in \op{Coker}(\overline{N}\otimes E)$. Moreover, $\overline{N}\in\mathcal{M}(\mu)$ as $\overline{N}$ is simple
(\cite[Proposition~2.6.8]{Di}). Hence we can apply Lemma~\ref{projfuncN}
and obtain that the intersection of $\mathcal{M}(\mu)$ with the additive
closure of $\overline{N}\otimes E$ consists just of direct sums of copies of
$\overline{N}$. Since $\overline{N}$ is simple, the cokernel of any
homomorphism between direct sums of copies of $\overline{N}$ is
isomorphic to a direct sum of copies of $\overline{N}$ as well
(\cite[Proposition~2.6.5(iii)]{Di}). Hence the intersection of
$\mathcal{M}(\mu)$ with  $\op{Coker}(\overline{N}\otimes E)$ also consists just
of direct sums of copies of $\overline{N}$. This implies  that $\overline{N}$
is projective in  $\op{Coker}(\overline{N}\otimes E)$.

By Theorem~\ref{cor2-061107} we know that
Kostant's problem has a positive answer for the module $L(y\cdot\mu)$. So, by part \eqref{smallN1}, it would suffice to show that
\begin{displaymath}
\dim\mathrm{Hom}_{\mathfrak{a}}(\overline{N},\theta \overline{N})=
\dim\mathrm{Hom}_{\mathfrak{a}}(L(y\cdot\mu),\theta L(y\cdot\mu))
\end{displaymath}
for all indecomposable projective functors $\theta$. This is true if $\theta$ is not the identity functor (Lemma~\ref{projfuncN}), otherwise Schur's Lemma (\cite[2.6.5]{Di}) does the job.
\end{proof}

Finally we get the following main result:

\begin{theorem}\label{thm-eqv1}
With the notation from above there is an equivalence of categories
\begin{eqnarray*}
\op{Coker}(L(y\cdot\mu)\otimes E)\cong \op{Coker}(\overline{N}\otimes E).
\end{eqnarray*}
\end{theorem}

\begin{proof}
By Lemma~\ref{projLy} and Proposition~\ref{smallN}, $L(y\cdot\mu)$
is projective in $\op{Coker}(L(y\cdot\mu)\otimes E)$ and
$\overline{N}$ is projective in $\op{Coker}(\overline{N}\otimes E)$.
By Theorem~\ref{cor2-061107}\eqref{cor2-061107.3}, Kostant's problem
has a positive solution for $L(y\cdot\mu)$. By Proposition~\ref{smallN},
Kostant's problem has a positive solution for $\overline{N}$. Hence
the claim follows from Proposition~\ref{smallN} and \cite[Theorem~5]{KM2}.
\end{proof}

\subsection{Categories of induced modules and
their equivalence}\label{s9.4}

We extend the categories $\op{Coker}(L(y\cdot\mu)\otimes E)$ and
$\op{Coker}(\overline{N}\otimes E)$ of $\mathfrak{a}$-modules to
categories of $\mathfrak{a}'$-modules by allowing arbitrary scalar
actions of the center of $\mathfrak{a}'$. Abusing notation we denote
the resulting categories by the same symbols.

\begin{lemma}\label{s9.2-lem1}
$\op{Coker}(L(y\cdot\mu)\otimes E)$ and
$\op{Coker}(\overline{N}\otimes E)$ are both admissible.
\end{lemma}

\begin{proof}
The conditions (L\ref{lll2}) and (L\ref{lll4}) are clear by definition,
so we have only to check the condition  (L\ref{lll3}). By Lemma~\ref{projLy},
$L(y\cdot\mu)$ is projective in $\op{Coker}(L(y\cdot\mu)\otimes E)$ and
$\overline{N}$ is projective in $\op{Coker}(\overline{N}\otimes E)$ by
Proposition~\ref{smallN}. In particular, all modules of the form
$L(y\cdot\mu)\otimes E$ and $\overline{N}\otimes E$, where $E$ is
finite-dimensional, are projective in the corresponding categories.
It follows that both categories  $\op{Coker}(L(y\cdot\mu)\otimes E)$ and
$\op{Coker}(\overline{N}\otimes E)$ have enough projectives. Now the
condition  (L\ref{lll3}) follows for instance from \cite[Section~5]{Au}.
\end{proof}

Lemma~\ref{s9.2-lem1} allows us to consider the category
$\mathcal{O}\{\mathfrak{p},\op{Coker}(L(y\cdot\mu)\otimes E)\}$ and
the category
$\mathcal{O}\{\mathfrak{p},\op{Coker}(\overline{N}\otimes E)\}$.
Both categories have a block decomposition with respect to central
characters. By \cite[Theorem~6.1]{Ma}, these blocks are equivalent
to module categories over some standardly stratified algebras (it is
easy to see that these algebras are even weakly properly stratified
in the sense of \cite{Fr}). Denote by
$\mathcal{O}\{\mathfrak{p},\op{Coker}(L(y\cdot\mu)\otimes
E)\}_{\op{int}}$ and
$\mathcal{O}\{\mathfrak{p},\op{Coker}(\overline{N}\otimes
E)\}_{\op{int}}$ the direct sums of all blocks corresponding to
integral central characters. The main result of this section is the
following statement:

\begin{theorem}\label{s9.2-thm2}
There is a blockwise equivalence of categories
\begin{eqnarray*}
\xi:\quad\mathcal{O}\{\mathfrak{p},\op{Coker}(\overline{N}\otimes E)\}_{\op{int}} &\cong&
\mathcal{O}\{\mathfrak{p}, \op{Coker}(L(y\cdot\mu)\otimes E)\}_{\op{int}},
\end{eqnarray*}
which sends proper standard modules to proper standard modules.
\end{theorem}

\begin{proof}
To construct a blockwise equivalence it is enough to verify the
assumptions of \cite[Theorem~5]{KM2}. The module $L(y\cdot\mu)$ is
projective in $\mathscr{C}:=\op{Coker}(L(y\cdot\mu)\otimes E)$
(Lemma~\ref{projLy}). Hence the induced module
$\Delta(\mathfrak{p},L(y\cdot\mu))$ is both standard and proper
standard in $\mathcal{O}\{\mathfrak{p}, \mathscr{C}\}$ for
any linear functional on the center of $\ma'$ which extends the
$\mathfrak{a}$-action on $L(y\cdot\mu)$. We pick the linear
functional such that the module $\Delta(\mathfrak{p},L(y\cdot\mu))$
is projective in some regular block of
$\mathcal{O}\{\mathfrak{p},\mathscr{C}\}_{\op{int}}$. It is easy to
see that all projective modules in
$\mathcal{O}\{\mathfrak{p},\mathscr{C}\}_{\op{int}}$ can be obtained
by translating $\Delta(\mathfrak{p},L(y\cdot\mu))$. In particular,
we have
\begin{displaymath}
\mathcal{O}\{\mathfrak{p},\mathscr{C}\}_{\op{int}}\cong
\op{Coker}(\Delta(\mathfrak{p},L(y\cdot\mu))\otimes E).
\end{displaymath}
Analogously
\begin{displaymath}
\mathcal{O}(\mathfrak{p}, \op{Coker}(\overline{N}\otimes E))_{\op{int}}\cong
\op{Coker}(\Delta(\mathfrak{p},\overline{N})\otimes E).
\end{displaymath}
for the same linear functional. To be able to apply
\cite[Theorem~5]{KM2} we just have to verify that Kostant's problem
has a positive answer for the modules
$\Delta(\mathfrak{p},\overline{N})$ and
$\Delta(\mathfrak{p},L(y\cdot\mu))$. This will be done in the
following Lemmas \ref{s9.2-lem3} and \ref{s9.2-lem4}. Hence there is
an equivalence of categories $\xi$, but it is left to show that is
sends proper standard objects to such. The partial ordering on the
simple modules in $\mathcal{O}\{\mathfrak{p},
\mathscr{C}\}_{\op{int}}$ induces a partial ordering on the simple
modules in $\mathcal{O}\{\mathfrak{p},
\op{Coker}(\overline{N}\otimes E)\}_{\op{int}}$, which defines a
stratified structure. Since proper standard modules have a
categorical definition, they will be sent to proper standard modules
by any blockwise equivalence.
\end{proof}

For a finite dimensional $\mg$-module $E$ we denote by $\tilde{E}$
its underlying $\ma$-module. Let $\tilde{E}_0$ be the direct sum of
finite dimensional $\ma$-submodules of $\tilde{E}$ where the center
of the reductive Lie algebra $\ma'$ acts trivially.

\begin{lemma}\label{s9.2-lem3a}
Kostant's problem has a positive answer for $\Delta(\mathfrak{p},L(\mu))$.
\end{lemma}

\begin{proof}
The module $\Delta(\mathfrak{p},L(\mu))$ is a quotient of the dominant Verma module and therefore Kostant's problem is affirmative
by \cite[6.9 (10)]{Ja2}.
\end{proof}

\begin{lemma}\label{s9.2-lem3}
Kostant's problem has a positive answer for $\Delta(\mathfrak{p},L(y\cdot\mu))$.
\end{lemma}

\begin{proof}
For any simple finite-dimensional $\mathfrak{g}$-module $E$ we have
\small
\begin{eqnarray}
\label{multi}
&&[U(\mathfrak{g})/\mathrm{Ann}_{U(\mg)}(\Delta(\mathfrak{p},L(\mu))):E]=
\dim\mathrm{Hom}_{\mathfrak{g}}
(\Delta(\mathfrak{p},L(\mu))\otimes E,\Delta(\mathfrak{p},L(\mu))
\end{eqnarray}
\normalsize
by Lemma~\ref{s9.2-lem3a} and \cite[6.8(3)]{Ja2}. Since for $\zeta\in\{\mu,
y\cdot\mu\}$, the module $\Delta(\mathfrak{p},L(\zeta))$ is a
projective standard module in its corresponding
$\op{Coker}$-category, the standard adjointness gives
\begin{displaymath}\label{firsteq}
\begin{array}{rcl}
\mathrm{Hom}_{\mathfrak{g}}
(\Delta(\mathfrak{p},L(\zeta))\otimes E,\Delta(\mathfrak{p},L(\zeta)))&=&
\mathrm{Hom}_{\mathfrak{g}}
(\Delta(\mathfrak{p},L(\zeta)),\Delta(\mathfrak{p},L(\zeta))\otimes E^*)\\
&=&\mathrm{Hom}_{\mathfrak{a}}(L(\zeta),L(\zeta)\otimes E_0^*)\\
&=&\mathrm{Hom}_{\mathfrak{a}}(L(\zeta)\otimes E_0,L(\zeta)).
\end{array}
\end{displaymath}
The latter is however independent of the choice of $\zeta$ by Theorem~\ref{cor2-061107}, and therefore
\small
\begin{equation}\label{s9.2-eq11}
\mathrm{Hom}_{\mathfrak{g}}
(\Delta(\mathfrak{p},L(\mu))\otimes E,\Delta(\mathfrak{p},L(\mu)))=
\mathrm{Hom}_{\mathfrak{g}}(\Delta(\mathfrak{p},L(y\cdot\mu))\otimes E,
\Delta(\mathfrak{p},L(y\cdot\mu )).
\end{equation}
\normalsize
The modules $L(y\cdot\mu)$ and $L(\mu)$ have the same
annihilator (by Theorem~\ref{cor2-061107} again), therefore the
modules $\Delta(\mathfrak{p},L(y\cdot\mu))$ and
$\Delta(\mathfrak{p},L(\mu))$ have the same annihilator by
\cite[Proposition~5.1.7]{Di}. Together with \eqref{s9.2-eq11} and
Lemma~\ref{s9.2-lem3a} we deduce that  Kostant's problem has a
positive answer for $\Delta(\mathfrak{p},L(y\cdot\mu))$.
\end{proof}

\begin{lemma}\label{s9.2-lem4}
Kostant's problem has a positive answer for $\Delta(\mathfrak{p},\overline{N})$.
\end{lemma}

\begin{proof}
Since $\Delta(\mathfrak{p},\overline{N})$ is a projective standard module in
the corresponding $\op{Coker}$-category, as in \eqref{firsteq} we have
\begin{displaymath}
\mathrm{Hom}_{\mathfrak{g}}
(\Delta(\mathfrak{p},\overline{N}),\Delta(\mathfrak{p},\overline{N})\otimes E)=
\mathrm{Hom}_{\mathfrak{a}}(\overline{N},\overline{N}\otimes E_0).
\end{displaymath}
Recall that $\overline{N}$ and $L(\mu)$ have the same annihilator
(Proposition~\ref{smallN}), and Kostant's map is surjective in both
cases (Theorem~\ref{cor2-061107} and Proposition~\ref{smallN}).
Together with \eqref{multi} we have
\begin{displaymath}
\mathrm{Hom}_{\mathfrak{g}}
(\Delta(\mathfrak{p},\overline{N}),\Delta(\mathfrak{p},\overline{N})
\otimes E)\cong \mathrm{Hom}_{\mathfrak{g}}(\Delta(\mathfrak{p},L(y\cdot\mu)),
\Delta(\mathfrak{p},L(y\cdot\mu))\otimes E).
\end{displaymath}
Now, $\Delta(\mathfrak{p},L(y\cdot\mu))$ and
$\Delta(\mathfrak{p},\overline{N})$ have the same annihilator
(Proposition~\ref{smallN} and \cite[Proposition 5.1.7]{Di}. So, the latter
equality and the fact that Kostant's problem has a positive  answer
for  $\Delta(\mathfrak{p},L(y\cdot\mu))$ imply that  Kostant's
problem has a positive  answer for
$\Delta(\mathfrak{p},\overline{N})$. This completes the proof.
\end{proof}

\subsection{The rough structure of
generalized Verma modules: main results}\label{s9.5}

The equivalence $\xi$ from Theorem~\ref{s9.2-thm2} induces a
bijection between the sets of the isomorphism classes of
indecomposable projective modules in the categories
\begin{displaymath}
\mathscr{Y}_{\overline{N}}=\mathcal{O}\{\mathfrak{p},\op{Coker}
(\overline{N}\otimes E)\}_{\op{int}}\quad\text{ and }\quad
\mathscr{Y}_{L(y\cdot\mu)}=\mathcal{O}\{\mathfrak{p},
\op{Coker}(L(y\cdot\mu)\otimes E)\}_{\op{int}}.
\end{displaymath}
Therefore $\xi$ also induces a bijection
\begin{displaymath}
\overline{\xi}:\quad\op{Irr}(\mathscr{Y}_{L(y\cdot\mu)})\rightarrow
\op{Irr}(\mathscr{Y}_{\overline{N}})
\end{displaymath}
between the sets of isomorphism classes of simple objects in
$\mathscr{Y}_{L(y\cdot\mu)}$ and $\mathscr{Y}_{\overline{N}}$ respectively.
This induces moreover a bijection
\begin{displaymath}
\hat{\xi}:\quad\op{Irr}^{\mathfrak{g}}(\mathscr{Y}_{\overline{N}})\rightarrow
\op{Irr}^{\mathfrak{g}}(\mathscr{Y}_{L(y\cdot\mu)})
\end{displaymath}
between the sets of isomorphism classes of the simple quotients, as
$\mg$-modules, of the modules from $\op{Irr}(\mathscr{Y}_{\overline{N}})$
and $\op{Irr}(\mathscr{Y}_{L(y\cdot\mu)})$ respectively. Each module $X\in
\op{Irr}^{\mathfrak{g}}(\mathscr{Y}_{\overline{N}})$ or
$\op{Irr}^{\mathfrak{g}}(\mathscr{Y}_{L(y\cdot\mu)})$ has the form
$L(\mathfrak{p},V_X)$ for a uniquely defined simple $\ma'$-module $V_X$.

As a consequence of Theorem~\ref{s9.2-thm2} we obtain the following result:

\begin{theorem}\label{s9.5-cor1}
For $X,Y\in \op{Irr}^{\mathfrak{g}}(\mathscr{Y}_{\overline{N}})$ we have the
following multiplicity formula in the category of $\mg$-modules:
\begin{displaymath}
[\Delta(\mathfrak{p},V_X):L(\mathfrak{p},V_Y)]=
[\Delta(\mathfrak{p},V_{\hat{\xi}(X)}):
L(\mathfrak{p},V_{\hat{\xi}(Y)})].
\end{displaymath}
\end{theorem}

\begin{proof}
Let $P(X)\in\mathscr{Y}_{\overline{N}}$ be an indecomposable projective,
whose head (as a $\mg$-module) is isomorphic to $X$. Then
$[\Delta(\mathfrak{p},V_X):L(\mathfrak{p},V_Y)]$ is just the dimension of
the homomorphism space from $P(X)$ to the proper standard module in
$\mathscr{Y}_{\overline{N}}$ corresponding to $X$ (see \cite[Section~5]{KM2}).
Exactly the same holds if we replace $X$ by $\hat{\xi}(X)$ and work with
the category $\mathscr{Y}_{L(y\cdot\mu)}$ instead of
$\mathscr{Y}_{\overline{N}}$. Since $\xi$ is an equivalence of categories
(Theorem~\ref{s9.2-thm2}) sending proper standard objects to proper standard
objects, the claim follows.
\end{proof}

\begin{remark}[Additional remarks to
Theorem~\ref{s9.5-cor1}]\label{finalremarks}
{\rm
Theorem~\ref{s9.5-cor1} describes only multiplicities of certain simple
subquotients of $\Delta(\mathfrak{p},V_X)$, namely multiplicities of those
simple subquotients, which occur as heads of indecomposable projectives in
$\mathscr{Y}_{\overline{N}}$. Following \cite{KM2} we call this the {\em rough
structure} of $\Delta(\mathfrak{p},V_X)$. The theorem reduces the question
about the rough structure of the module $\Delta(\mathfrak{p},V_X)$ to the
analogous question for the module $\Delta(\mathfrak{p},V_{\hat{\xi}(X)})$.
The latter module is an object of $\mathcal{O}$ and hence the problem can
be solved inductively using  the Kazhdan-Lusztig combinatorics.
}
\end{remark}

Let $L$ be as in Subsection~\ref{s9.3}. Then the module
$\Delta(\mathfrak{p},L)$ has generalized trivial integral central character,
and $L(\mathfrak{p},L)$ is the simple top of some indecomposable projective
module, $P$ say,  in $\mathscr{Y}_{\overline{N}}$. Let $\mathscr{X}$ be the
block of $\mathscr{Y}_{\overline{N}}$ corresponding to the trivial central
character. It contains $P$ by construction. By Theorem~\ref{s9.2-thm2} and
Subsection~\ref{s5.4}, simple modules in $\mathscr{X}$ are (bijectively)
indexed by $(x,w)\in\mathds{I}(\mathbf{R}')$. Therefore, there is a
pair $(x,w)$ for each $\Delta(\mathfrak{p},L)$ in $\mathscr{X}$.
Theorem~\ref{s9.5-cor1} allows us to formulate the following
irreducibility criterion for generalized Verma modules:

\begin{theorem}\label{s9.5-cor2}
Let $(x,w)$ be the pair associated with $\Delta(\mathfrak{p},L)$.
Then the module  $\Delta(\mathfrak{p},L)$ is irreducible if and only
if $w=\overline{w}$.
\end{theorem}

\begin{proof}
Theorem~\ref{s9.5-cor1} reduces this to the category
$\mathcal{O}\{\mathfrak{p},\mathscr{A}^{\mathbf{R}'}\}$ from
Subsection~\ref{s5.4}. For the category
$\mathcal{O}\{\mathfrak{p},\mathscr{A}^{\mathbf{R}'}\}$ the statement follows
from the proof of Theorem~\ref{thm53}.
\end{proof}

\begin{remark}[Unnecessary restrictions]\hfill\\
\label{simplicity}
{\rm
\begin{enumerate}[(i)]
\item The restriction of integrability for the central character is not
really essential and can be taken away using methods proposed by
Soergel in \cite[Bemerkung~1]{Sperv} on the reduction of the
Kazhdan-Lusztig conjecture to the integral case.
\item In this paper we only worked with the trivial central
character to avoid even more notation. The singular case follows by
translation to the regular case, using our results there and
translating back (invoking the fact that the composition of these
translation functors is just a multiple of the identity).
\end{enumerate}
}
\end{remark}

\end{document}